\UseRawInputEncoding

\documentclass[12pt,reqno]{article}

\usepackage{amsfonts}
\usepackage{mathrsfs}

\usepackage{amsmath,amsthm,amsfonts,amssymb,latexsym,mathrsfs,color}

\usepackage[colorlinks=true,
linkcolor=blue, filecolor=webbrown, citecolor=webred ]{hyperref}
\definecolor{webblue}{rgb}{0,.5,0}
\definecolor{webred}{rgb}{0,.5,0}
\definecolor{webbrown}{rgb}{.6,0,0}

\newtheorem{thm}{Theorem}[section]

\newtheorem{cla}{Claim}

\theoremstyle{definition}
\newtheorem{definition}[thm]{Definition}

\newtheorem{rem}{Remark}

\numberwithin{equation}{section}

\numberwithin{equation}{section}

\newcommand{\Z}{\mathbb Z}
\newcommand{\Q}{\mathbb Q}

\title{Coefficientwise Hankel-total positivity from a universal cumulant lifting\thanks{This work was supported in part by the National Key R\&D Program of China (No.~2024YFA1013903) and the National Natural Science Foundation of China (No.~12471310).}}
\author{ Bao-Xuan Zhu\thanks{
{\it Email address:} {bxzhu@jsnu.edu.cn (B.-X. Zhu)}}}
\date{\footnotesize
School of Mathematics and Statistics, Jiangsu Normal University,
Xuzhou 221116, P.R. China}

\usepackage[T1]{fontenc}
\usepackage{lmodern,mathtools,microtype,needspace,booktabs,titlesec,array}
\theoremstyle{plain}
\newtheorem{theorem}[thm]{Theorem}
\newtheorem{lemma}[thm]{Lemma}
\newtheorem{proposition}[thm]{Proposition}
\newtheorem{corollary}[thm]{Corollary}
\newtheorem{conjecture}[thm]{Conjecture}
\newtheorem{claim}[cla]{Claim}
\theoremstyle{remark}

\titleformat{\part}[block]{\normalfont\Large\bfseries\raggedright}
  {Part~\thepart.}{0.65em}{}
\titleformat{\section}{\normalfont\large\bfseries\raggedright}
  {\thesection.}{0.65em}{}
\titleformat{\subsection}{\normalfont\normalsize\bfseries\raggedright}
  {\thesubsection.}{0.65em}{}
\titleformat{\subsubsection}{\normalfont\normalsize\itshape\raggedright}
  {\thesubsubsection.}{0.65em}{}
\titlespacing*{\part}{0pt}{2.2ex plus 1ex minus .2ex}{1.5ex plus .2ex}
\titlespacing*{\section}{0pt}{2.4ex plus 1ex minus .2ex}{1.1ex plus .2ex}
\titlespacing*{\subsection}{0pt}{1.8ex plus .8ex minus .2ex}{.7ex plus .2ex}
\renewcommand{\thepart}{\arabic{part}}
\renewenvironment{abstract}
  {\small\begin{center}\normalfont\bfseries\abstractname\end{center}\noindent\ignorespaces}
  {\par}
\allowdisplaybreaks[2]
\makeatletter
\renewcommand{\@maketitle}{%
  \newpage\null\vskip .5em
  \begin{center}
  {\LARGE\bfseries\@title\par}\vskip 1.2em
  {\large\@author\par}\vskip .8em
  {\@date\par}
  \end{center}\par\vskip 1em}
\makeatother
\hypersetup{
 pdftitle={Coefficientwise Hankel-total positivity from a universal cumulant lifting},
 pdfauthor={Bao-Xuan Zhu},
 pdfsubject={Functional digraphs; Stieltjes and Hamburger cumulants; coefficientwise Hankel total positivity},
 bookmarksdepth=3
}
\newcommand{\R}{\mathbb R}
\newcommand{\Nzero}{\mathbb N_0}
\newcommand{\E}{\mathbb E}
\newcommand{\by}{\boldsymbol y}
\newcommand{\bq}{\boldsymbol q}
\newcommand{\bs}{\boldsymbol s}
\newcommand{\bx}{\boldsymbol x}
\newcommand{\bJ}{\boldsymbol J}

\newcommand{\dif}{\mathrm d}
\newcommand{\cP}{\mathcal P}
\newcommand{\Dshift}{\mathscr D}

\newcommand{\N}{\mathbb N_0}

\newcommand{\TP}{\mathrm{TP}}

\DeclareMathOperator{\diag}{diag}
\DeclareMathOperator{\adj}{adj}

\definecolor{addedblue}{rgb}{0,0.16,0.67}

\definecolor{revisionpurple}{rgb}{0.43,0.08,0.56}

\begin{document}
\maketitle

\begin{abstract}
Using probabilistic methods, we establish universal cumulant lifting criteria for Stieltjes and
Hamburger moment sequences. In the Stieltjes case, we characterize
coefficientwise Hankel total positivity after convolution with every
initial moment sequence by requiring each once-shifted cumulant
sequence to be Stieltjes. We also prove equivalence with the
numerical moment property at every nonnegative intensity.
In the Hamburger case, we use twice-shifted cumulants to obtain
coefficientwise-positive leading Hankel determinants, a degree bound,
an explicit highest-degree component, and drift invariance.
We work with formal series and require neither exponential moments
nor moment determinacy. We establish
composition, binomial-convolution, and substitution results, and
apply them to implicit equations and autonomous differential equations, tree and forest
families, and lognormal compound-Poisson moments in a unified approach. From known
moment sequences, we obtain general coefficientwise Hankel-totally positive
families containing the classical Bessel, Lah, and Laguerre
polynomials. We not only give new proofs of coefficientwise Hankel-total positivity of many combinatorial polynomials, but also prove coefficientwise Hankel total positivity for the component polynomials of functional
digraphs confirming a conjecture of Sokal and establish fixed-shape Schl\"afli-Gessel-Seo polynomials positivity supporting another conjecture of Sokal.

\end{abstract}

\clearpage
{\renewcommand{\baselinestretch}{0.90}\normalsize
\makeatletter
\patchcmd{\l@section}{1.0em}{0.25em}{}{\PackageError{manuscript}{Cannot adjust contents spacing}{}}
\makeatother
\tableofcontents}
\clearpage

\section{Introduction}\label{sec:introduction}

Total positivity of matrices has been deeply studied in various branches of mathematics. For example, Brenti \cite{Bre89} developed total positivity to study unimodality and log-concavity of sequences. Lusztig extended total
positivity to reductive groups~\cite{Lusztig1994}. Fomin and
Zelevinsky~\cite{FominZelevinsky1999} studied totally nonnegative
parts of double Bruhat cells and their parametrizations. Rietsch
related totally positive Toeplitz matrices to the quantum cohomology
of partial flag varieties~\cite{Rietsch2003}. Kodama and
Williams~\cite{KodamaWilliams2014} connected the totally nonnegative
Grassmannian with regular soliton solutions of the
Kadomtsev--Petviashvili equation. Belton, Guillot, Khare and Putinar~\cite{BGKP2016} obtained a
classification for entrywise polynomial preservers of positive
semidefiniteness in fixed dimension. In their work on
moment-sequence transforms~\cite{BGKP2022}, they classified
entrywise maps preserving moment sequences and totally
nonnegative Hankel matrices, and also treated multivariable
transforms. We also refer the reader to monographs
\cite{Kar68,Pinkus2010} for more details.

Recent years, there many important progress on coefficientwise Hankel-total positivity. For example, in \cite{PSZ2023},  P\'{e}tr\'{e}olle, Sokal and Zhu developed the
theory for coefficientwise Hankel-total positivity based on lattice paths and branched
continued fractions. See Zhu \cite{Zhu2024,Zhu2026Riordan} for methods from production matrices and Riordan arrays in combinatorics. There are still many difficult open problems on coefficientwise total positivity. We also refer the reader to
\cite{DDPS23,PSZ2023,Sokal_totalpos,PS19,PSZ2023,Zhu20,Zhu21,Zhu212,Zhu2024,Zhu2026Riordan} for more details. In this present paper, we will develop the probabilistic method to prove coefficientwise Hankel-total positivity. We address the functional-digraph problem within this development by characterizing the exponential cumulant families for which a
universal coefficientwise lifting property holds.

\subsection{Functional digraphs and Sokal's conjecture}

A \emph{functional digraph} on \([n]=\{1,\ldots,n\}\) is the directed
graph of a map from \([n]\) to itself: every vertex has out-degree one.
Each weakly connected component consists of a directed cycle with rooted
trees attached to its cyclic vertices. Let \(\psi_{n,k}\) be the number of such
digraphs with \(k\) components, and put
\begin{equation}\label{eq:display-001}
 \psi_n(y)=\sum_{k=0}^n\psi_{n,k}y^k,\qquad \psi_0(y)=1.
\end{equation}
We call \(\psi_n(y)\) the \emph{functional digraph component polynomial}.
For rooted labeled trees on \(n\) vertices, we write
\(T(t)\) for the corresponding exponential generating function. It is known that
\begin{equation}\label{eq:tree-digraph-egf}
 T(t)=\sum_{n\ge1}n^{n-1}\frac{t^n}{n!},\qquad T(t)=te^{T(t)}.
\end{equation}
By the exponential
formula~\cite[Corollary~5.1.6]{Stanley}, we obtain
\begin{equation}\label{eq:univariate-digraph-egf}
 \sum_{n\ge0}\psi_n(y)\frac{t^n}{n!}=(1-T(t))^{-y}.
\end{equation}
Sokal studied these functional digraph
component polynomials and proposed the next conjecture~\cite[Conjecture~6.9(b)]{Sokal2023}.

\begin{conjecture}[\cite{Sokal2023}]
\label{conj:sokal-digraph}
The functional digraph component polynomial sequence \((\psi_n(y))_{n\ge0}\)
is coefficientwise Hankel-totally positive in \(y\).
\end{conjecture}
Chen and Sokal~\cite{ChenSokal2024} studied related row-generating
polynomials associated with the entries \(\binom nk n^{n-k}\), using a model of partial functional digraphs.
To describe the progress on this conjecture, we assign an additional
weight \(q\) to each cyclic vertex and denote the resulting polynomials by
\(\psi_n(q,y)\). Their generating functions are
\begin{align}
 \sum_{n\ge0}\psi_n(q,y)\frac{t^n}{n!}
    &=(1-qT(t))^{-y},\label{eq:scalar-egf}\\
 C_q(t):=-\log(1-qT(t))
    &=\sum_{n\ge1}c_n(q)\frac{t^n}{n!}.\label{eq:connected-egf}
\end{align}
By \eqref{eq:coefficient-extraction-notation}, we extract the coefficient
of \(y\) in \(\psi_n(q,y)\) to select precisely the digraphs with
one weakly connected component. Hence
\(c_n(q)=[y]\psi_n(q,y)\).
At \(q=1\), we recover \(\psi_n(1,y)=\psi_n(y)\).
Conjecture~\ref{conj:sokal-digraph} was also stated as
Conjecture~6.16(ii) in~\cite{Zhu2024} and the following result was proved in \cite[Proposition~6.18]{Zhu2024}.

\begin{thm}[\cite{Zhu2024}]
\label{prop:connected-stieltjes}
 For fixed \(q,y\geq0\), the sequence \((\psi_n(q,y))_{n\ge0}\) is a numerical Stieltjes moment sequence. In particular, for each fixed \(q\ge0\), \((c_{n+1}(q))_{n\ge0}\) is a Stieltjes moment sequence.
\end{thm}

The result concerns fixed nonnegative parameter values. Coefficientwise Hankel-total positivity is stronger than Stieltjes moment property. In Conjecture~\ref{conj:sokal-digraph}, we require every coefficient
of each Hankel minor to be nonnegative. We extend Theorem \ref{prop:connected-stieltjes} to the following result for coefficientwise Hankel-total positivity, which in particular implies Conjecture~\ref{conj:sokal-digraph}.

\begin{theorem}[Functional-digraph Hankel-total positivity]\label{thm:digraph-scalar}
For every fixed \(q\ge0\), the
sequence \((\psi_n(q,y))_{n\ge0}\) is coefficientwise Hankel-totally
positive in \(y\). In particular, \((\psi_n(y))_{n\ge0}\) is coefficientwise Hankel-totally
positive in \(y\).
\end{theorem}

We will extend Theorem \ref{thm:digraph-scalar} by assigning independent variables to different component colors.

\begin{definition}
\label{def:colored-digraph}
Fix \(m\ge1\) and \(\bq=(q_1,\ldots,q_m)\in[0,\infty)^m\).
Assign one color \(\ell\in\{1,\ldots,m\}\) to each component of a
functional digraph. Give that component weight \(y_\ell\) and each of
its cyclic vertices weight \(q_\ell\). The weight of a colored digraph
is the product of its component and vertex weights. Let
\(\Psi_n(\bq;\by)\) be the sum of these weights over all digraphs on
\([n]\) and all component colorings, where \(\by=(y_1,\ldots,y_m)\).
We call \(\Psi_n(\bq;\by)\) the
\emph{colored functional digraph polynomial}.
\end{definition}

By the exponential formula, we obtain the generating function
\begin{equation}\label{eq:colored-egf}
 \sum_{n\ge0}\Psi_n(\bq;\by)\frac{t^n}{n!}
   =\prod_{\ell=1}^m(1-q_\ell T(t))^{-y_\ell}
   =\exp\!\left(\sum_{\ell=1}^m y_\ell C_{q_\ell}(t)\right).
\end{equation}

Taking \(m=1\), \(q_1=q\), and \(y_1=y\), we obtain
\(\Psi_n((q);(y))=\psi_n(q,y)\). We extend Theorem~\ref{thm:digraph-scalar} to the following joint positivity for the colored model.

\Needspace{11\baselineskip}
\begin{theorem}[Colored functional-digraph Hankel-total positivity]\label{thm:digraph}
For every fixed \(\bq\in[0,\infty)^m\),
\((\Psi_n(\bq;\by))_{n\ge0}\) is coefficientwise
Hankel-totally positive jointly in \(y_1,\ldots,y_m\).
\end{theorem}

Theorems~\ref{thm:digraph-scalar} and~\ref{thm:digraph} follow
from the Stieltjes lifting criterion below; their proof is given in
Section~\ref{sec:digraph}.

\subsection{A universal lifting criterion for coefficientwise Hankel-total positivity}


We will use Stieltjes moment property to present a new universal lifting criterion for coefficientwise
Hankel-total positivity. To place the lifting criterion in context, we first recall the role of
classical cumulants in probability and statistics. For a real random
variable \(X\) with finite absolute moments of every order, we define
the moment and cumulant formal power series by
\begin{equation*}\tag{1.A}\label{intro:eq:formal-cumulants}
 A_X(t)=\sum_{n\ge0}\E[X^n]\frac{t^n}{n!},\qquad
 K_X(t)=\log A_X(t)=\sum_{n\ge1}\kappa_n(X)\frac{t^n}{n!}.
\end{equation*}
We use the normalization \(A_X(0)=1\), and require no convergence
of either formal power series. In particular,
\(\kappa_1(X)=\E[X]\) and
\(\kappa_2(X)=\E[X^2]-(\E[X])^2\); thus the first two cumulants
are the mean and variance. For independent random variables \(X,Y\)
with all absolute moments finite, we have
\(K_{X+Y}(t)=K_X(t)+K_Y(t)\), as a formal identity, and hence
\(\kappa_n(X+Y)=\kappa_n(X)+\kappa_n(Y)\).
For the classical theory and its formulation through partition
lattices, see Speed~\cite{Speed1983}.

In enumerative combinatorics, we use the same moment--cumulant
relation to pass between arbitrary labeled structures and their
connected components. By summing products of cumulants over set
partitions, we recover the moments and give the independently
colored version in \eqref{eq:bell-partitions}.
By the exponential formula~\cite[Corollary~5.1.6]{Stanley}, we may
therefore regard the connected-component coefficients \(c_n(q)\)
in \eqref{eq:connected-egf} as cumulants of the functional-digraph
family. In this way, we connect Sokal's conjecture with the
positivity properties of moment--cumulant transformations.

In the classical moment problem, for a real sequence, it is a
\emph{Stieltjes moment sequence}, that is,
\begin{equation}\label{eq:display-003}
 a_n=\int_0^\infty x^n\,d\mu(x)\quad(n\ge0)
\end{equation}
for a nonnegative measure with all moments finite, where the zero measure is
allowed. For Stieltjes moment sequences, one
requires support in \([0,\infty)\); for Hamburger moments, one
allows support on \(\R\). One also studies these questions through
the linear functional \(\mathcal L(x^n)=a_n\): for a Stieltjes
representation, one requires
\(\mathcal L(p(x)^2)\ge0\) and
\(\mathcal L(xp(x)^2)\ge0\) for every real polynomial \(p(x)\).
We recall the equivalent matrix criteria in
Lemma~\ref{lem:moment-criteria}; see
Schm\"udgen~\cite[Chapter~3]{Schmudgen2017}.

In probability, one encounters the shifted cumulant condition
through compound-Poisson distributions. Let \(d,\lambda\ge0\),
let \(N_y\) have Poisson distribution with mean \(\lambda y\),
and let \(J_1,J_2,\ldots\) be independent copies of a nonnegative
random variable \(J\), independent of \(N_y\), with all moments
finite. For \(y\ge0\), set
\(X_y=dy+\sum_{h=1}^{N_y}J_h\), where an empty sum equals zero.
By the standard compound-Poisson cumulant formula
\cite[Chapter~15]{LastPenrose2017}, we have
\begin{equation*}\tag{1.C}\label{intro:eq:compound-cumulants}
 \kappa_n(X_y)=y\bigl(d\,\mathbf1_{\{n=1\}}
                    +\lambda\E[J^n]\bigr)\qquad(n\ge1).
\end{equation*}
Here \(\mathbf1_{\{n=1\}}\) equals one for \(n=1\) and zero
otherwise. Writing \(\delta_0\) for unit mass at zero and
\(P_J\) for the distribution of \(J\), we obtain
\(\kappa_{n+1}(X_1)=\int_0^\infty x^n\,\nu(\dif x)\)
with the finite measure
\(\nu(\dif x)=d\,\delta_0(\dif x)+\lambda xP_J(\dif x)\).
Thus we can interpret the once-shifted cumulants as Stieltjes
moments, including the deterministic drift as an atom at zero.

More generally, probabilists construct nondecreasing L\'evy
processes, called \emph{subordinators}, from Bernstein functions.
For a Bernstein function \(\phi(s)\) with \(\phi(0)=0\), one
obtains a subordinator \((S_y)_{y\ge0}\) with
\begin{equation*}\tag{1.D}\label{intro:eq:subordinator-laplace}
 \E[e^{-sS_y}]=e^{-y\phi(s)}\qquad(s,y\ge0).
\end{equation*}
Here one calls \(\phi(s)\) a Bernstein function when it is
nonnegative, smooth for \(s>0\), and satisfies
\((-1)^{r-1}\phi^{(r)}(s)\ge0\) for \(r\ge1\).
For the L\'evy representation and the convolution-semigroup
correspondence, see Schilling, Song and
Vondra\v cek~\cite[Theorems~3.2 and~5.2]{SSV2012}.
We recall the corresponding cumulant representation, under the
additional requirement of finite moments of every order, in
Lemma~\ref{lem:representation}. Through this representation, we
also explain the shift by one in the lifting criterion.

For fixed nonnegative indeterminates, one can use independent sums
to combine such moment families. In this paper, we retain the
intensities as independent indeterminates and prove nonnegativity
of every coefficient in every Hankel minor. We also allow an
arbitrary initial Stieltjes moment sequence, without normalization
of its zeroth moment. In Theorem~\ref{thm:lifting}, we characterize
this universal coefficientwise preservation by the shifted
cumulant condition. We work with formal power series throughout;
in particular, we do not require an exponential moment or
uniqueness of a representing measure.


Let \(K(t)=\sum_{n\ge1}\kappa_n t^n/n!\) be any real formal power
series, and write \(e^{yK(t)}=\sum_{n\ge0}B_n(y)t^n/n!\). To formulate the universal lifting property, we state the following equivalence.

\Needspace{20\baselineskip}
\begin{theorem}[One-parameter lifting criterion]\label{thm:scalar-lifting}
The following statements are equivalent:
\begin{enumerate}
\renewcommand{\labelenumi}{\textup{(\roman{enumi})}}
\item The shifted cumulant sequence
\((\kappa_{n+1})_{n\ge0}\) is a Stieltjes moment sequence.
\item The polynomial sequence \((B_n(y))_{n\ge0}\) is coefficientwise
Hankel-totally positive in \(y\).
\item For every real Stieltjes moment sequence
\(\boldsymbol a=(a_n)_{n\ge0}\), the sequence
\begin{equation}\label{eq:scalar-lifting-convolution}
 F_n(y)=\sum_{k=0}^n\binom nk a_{n-k}B_k(y)\qquad(n\ge0)
\end{equation}
is coefficientwise Hankel-totally positive in \(y\).
\end{enumerate}
\end{theorem}

The last assertion is the \emph{universal lifting} property: arbitrary
initial Stieltjes moments may be convolved with the component family.
Closure of real Stieltjes moment
sequences under binomial convolution is classical (see \cite{PS64},
\cite[Corollary~4.3]{Wang-Zhu} for instance). The lifting assertion extends it to coefficientwise Hankel-total positivity.
We obtain this theorem from the multivariate result below with
\(m=1\).

We now retain independent intensities for different formal power series.
Fix \(m\ge1\).  For \(\ell=1,\ldots,m\), let
\(\kappa_1^{(\ell)},\kappa_2^{(\ell)},\ldots\) be real numbers,
and define
\begin{equation}\label{eq:bell-egf}
 \sum_{n\ge0}M_n(\by)\frac{t^n}{n!}
 =\exp\!\left(\sum_{\ell=1}^{m}y_\ell K_\ell(t)\right),
 \qquad
 K_\ell(t)=\sum_{n\ge1}\kappa_n^{(\ell)}\frac{t^n}{n!}.
\end{equation}
For a real sequence \(\boldsymbol a=(a_n)_{n\ge0}\), define
\begin{equation}\label{eq:universal-operator}
 F_n(\by):=(\mathcal T_{\by}\boldsymbol a)_n
 =\sum_{k=0}^{n}\binom nk a_{n-k}M_k(\by).
\end{equation}
If \(A(t)=\sum_{n\ge0}a_nt^n/n!\),  then (\ref{eq:universal-operator}) is equivalent to the formal identity
\begin{equation}\label{eq:universal-egf}
 \sum_{n\ge0}F_n(\by)\frac{t^n}{n!}
 =A(t)\exp\!\left(\sum_{\ell=1}^{m}y_\ell K_\ell(t)\right).
\end{equation}
In particular, \(F_0(\by)=a_0\).

We characterize the admissible cumulants as follows, which for
\(m=1\) implies Theorem~\ref{thm:scalar-lifting}.

\Needspace{19\baselineskip}
\begin{theorem}[Universal multivariate Stieltjes--cumulant equivalence]
\label{thm:lifting}
For the polynomial sequence \((M_n(\by))_{n\ge0}\) defined by \eqref{eq:bell-egf}, the following
conditions are equivalent:
\begin{enumerate}
\renewcommand{\labelenumi}{\textup{(\roman{enumi})}}
\item For each \(\ell=1,\ldots,m\), the shifted sequence
\((\kappa_{n+1}^{(\ell)})_{n\ge0}\) is a Stieltjes moment sequence.
\item \((M_n(\by))_{n\ge0}\) is coefficientwise Hankel-totally positive
jointly in \(y_1,\ldots,y_m\).
\item For every fixed \(\by\in[0,\infty)^m\),
\((M_n(\by))_{n\ge0}\) is a Stieltjes moment sequence.
\item For every real Stieltjes moment sequence \(\boldsymbol a\),
\((F_n(\by))_{n\ge0}\) in
\eqref{eq:universal-operator} is coefficientwise Hankel-totally positive
jointly in \(y_1,\ldots,y_m\).
\item For every real Stieltjes moment sequence \(\boldsymbol a\) and
every fixed \(\by\in[0,\infty)^m\), \((F_n(\by))_{n\ge0}\) is a
Stieltjes moment sequence.
\end{enumerate}
\end{theorem}

Under any of these conditions, we have \(\kappa_n^{(\ell)}\ge0\)
for all \(n,\ell\), so we need no separate sign assumption on the
cumulants. In condition~\textup{(iv)}, we allow every initial real
Stieltjes moment sequence and assert coefficientwise positivity in \(\by\).

We will prove Theorem \ref{thm:lifting} in Section~\ref{sec:lifting}. For sufficiency of Theorem \ref{thm:lifting}, the proof combines the Poisson square-integral
identity of Last and Penrose~\cite{LastPenrose2011}, Schur
alternants, and a partitioned Mecke expansion. The resulting
kernels are independent of the indeterminates. Necessity follows
by extracting the lowest-degree terms of shifted Hankel minors.

\subsection{A Hamburger extension on the real line}\label{sec:hamburger}

A real sequence \((a_n)_{n\ge0}\) is a \emph{Hamburger moment sequence}
if \(a_n=\int_{-\infty}^\infty x^n\,\mu(\dif x)\) for a nonnegative measure
\(\mu\) on \(\R\) with all absolute moments finite.
The zero measure is allowed.
The Hamburger moment theorem says that this is equivalent to positive
semidefiniteness of every finite Hankel matrix
\([a_{i+j}]_{i,j=0}^{d}\); see \cite[Chapter~3]{Schmudgen2017}.

For arbitrary real support, the universal coefficientwise
criterion, stated in the following as in Theorem~\ref{thm:lifting}, involves twice-shifted cumulants and leading principal
determinants. In addition, Remark~\ref{rem:principal-minors} explains why its
coefficientwise assertion does not extend to arbitrary principal
minors.

\begin{theorem}[Universal Hamburger lifting]\label{thm:hamburger}
Let \(M_n\) and \(F_n\) be defined by \eqref{eq:bell-egf} and
\eqref{eq:universal-operator} with real cumulants and an arbitrary
initial real sequence \(\boldsymbol a=(a_n)_{n\ge0}\), and put
\begin{equation}\label{eq:display-017}
 H_r^a(\by)=\det[F_{i+j}(\by)]_{i,j=0}^{r-1}.
\end{equation}
The following conditions are equivalent:
\begin{enumerate}
\renewcommand{\labelenumi}{\textup{(\roman{enumi})}}
\item For every \(\ell\), the twice-shifted sequence
\((\kappa_{n+2}^{(\ell)})_{n\ge0}\) is a Hamburger moment sequence.
\item For every initial Hamburger moment sequence and \(r\ge1\),
\(H_r^a(\by)\in\R_{\ge0}[\by]\).
\item For every initial Hamburger moment sequence and every
\(\by\in[0,\infty)^m\), \((F_n(\by))_{n\ge0}\) is a Hamburger
moment sequence.
\item For every \(\by\in[0,\infty)^m\),
\((M_n(\by))_{n\ge0}\) is a Hamburger moment sequence.
\end{enumerate}
Under these conditions, fix an initial Hamburger moment sequence
and let \(D=\binom r2\) and
\(v(\by)=\sum_\ell y_\ell\kappa_2^{(\ell)}\).
Then
\begin{equation}\label{eq:hamburger-top-degree}
 \deg_{\by}H_r^a\le D,\qquad
 [H_r^a]_D
 =a_0^r\!\left(\prod_{j=0}^{r-1}j!\right)v(\by)^D,
\end{equation}
where \([P]_D\) denotes the homogeneous part of total degree \(D\).
The determinant is independent of all first cumulants
\(\kappa_1^{(\ell)}\).
\end{theorem}

Using the common moment and Poisson tools from Section~\ref{sec:lifting}, we will prove Theorem~\ref{thm:hamburger} in
Section~\ref{sec:hamburger-proof}. To clarify the scope of Theorem~\ref{thm:hamburger}, we record
four remarks on drift invariance, Gaussian components, principal
minors, and composition.

\begin{rem}[Invariance under common polynomial drifts]\label{rem:common-drifts}
The fact that these determinants do not depend on the first
cumulants has an algebraic explanation.
If any polynomial \(c(\by)\in\R[\by]\) is used to replace
\(F_n\) by its common binomial translate
\begin{equation}\label{eq:display-027}
 \widetilde F_n(\by)
 =\sum_{k=0}^{n}\binom nk c(\by)^{n-k}F_k(\by),
\end{equation}
then \(\widetilde H_r^a:=\det[\widetilde F_{i+j}]_{i,j=0}^{r-1}
=H_r^a\). This is the parameter form of
Layman's binomial invariance~\cite[Theorem~1]{Layman2001}.
Indeed, writing \(\mathsf H_r=[F_{i+j}]_{i,j=0}^{r-1}\), one has
\(\widetilde{\mathsf H}_r=P_r\mathsf H_rP_r^{\mathsf T}\), where
\(P_r\) is the \(r\times r\) lower triangular matrix with
\((P_r)_{i,j}=\binom ij c(\by)^{i-j}\) for \(0\le j\le i<r\)
and zero entries for \(j>i\). Its diagonal entries are one.
Hence the determinant assertion also survives arbitrary common polynomial
drifts, including drifts with negative coefficients.
\end{rem}

\begin{rem}[Gaussian cumulants and moment determinants]\label{rem:gaussian-determinants}
For the pure exponential family with
\begin{equation}\label{eq:display-028}
 K_\ell(t)=b_\ell t+\frac{\sigma_\ell^2t^2}{2},
 \qquad b_\ell\in\R,\quad \sigma_\ell^2\ge0,
\end{equation}
the representing measure in \textup{(i)} is
\(\nu_\ell=\sigma_\ell^2\delta_0\).
Take the initial measure \(\delta_0\).
By Gaussian scaling and \eqref{eq:hamburger-top-degree}, we obtain
the exact formula
\begin{equation}\label{eq:gaussian-determinants}
 \det[M_{i+j}(\by)]_{i,j=0}^{r-1}
 =\left(\prod_{j=0}^{r-1}j!\right)
       \left(\sum_\ell\sigma_\ell^2y_\ell\right)^{\binom r2}.
\end{equation}
This is also the classical Gaussian moment determinant, obtained from
the monic Hermite squared norms; see
\cite[Table~18.3.1 and Section~18.2 (ix)]{DLMF-Norms}.
This family lies beyond the Stieltjes theorem: when the mean is zero
and the variance is positive, the once-shifted cumulant sequence
starts with \(0,\sigma^2\), which cannot be Stieltjes.
\end{rem}

\begin{rem}[Leading determinants versus arbitrary principal minors]\label{rem:principal-minors}
The coefficientwise conclusion in Theorem~\ref{thm:hamburger}
cannot be extended to every principal Hankel minor.
For example, let
\begin{equation}\label{eq:display-029}
 K(t)=e^t-1-10t.
\end{equation}
Then \(\kappa_1=-9\) and \(\kappa_n=1\) for \(n\ge2\), so
\((\kappa_{n+2})\) is represented by \(\delta_1\).
Nevertheless, the pure-family fourth moment is
\begin{equation}\label{eq:hamburger-counterexample}
 M_4(y)=y-33y^2+486y^3+6561y^4.
\end{equation}
It is the principal \(1\times1\) minor at Hankel index \(2\), and it
has a negative coefficient.  Numerically it remains the fourth moment
of \(N_y-10y\), where \(N_y\) is Poisson with mean \(y\).
Thus leading principal determinants and arbitrary principal minors
must be distinguished in the real-line extension.
\end{rem}

\begin{rem}[Failure of composition closure on the real line]\label{rem:hamburger-composition}
The composition closure in Theorem~\ref{thm:composition} also
uses the nonnegative-support hypothesis.  The exponent \(K(t)=t^2/2\)
satisfies the Hamburger criterion.  Its self-composition
\(K(K(t))=t^4/8\) has second cumulant zero and fourth cumulant
\(3\).  A nonnegative representing measure for its twice-shifted
sequence would have mass zero but a positive second moment, which
is impossible.
\end{rem}

\subsection{Notation and concepts}

We use the following notation for coefficient extraction and positivity.
For a formal power series \(f(x)=\sum_{j\ge0}f_jx^j\), we write
\begin{equation}\label{eq:coefficient-extraction-notation}
 [x^k]f(x)=f_k\qquad(k\ge0).
\end{equation}
Thus, for an exponential generating function
\(F(t)=\sum_{n\ge0}a_nt^n/n!\), we have
\(a_n=n![t^n]F(t)\). For several variables, we use
\([\boldsymbol z^{\boldsymbol\alpha}]\) to extract the coefficient of
\(z_1^{\alpha_1}\cdots z_d^{\alpha_d}\).
We write \(\R,\Q,\Z\) for the real numbers, rational numbers and
integers, respectively, and put \(\R_{\ge0}=[0,\infty)\) and
\(\N=\Nzero=\{0,1,2,\ldots\}\).
For a multi-index \(\boldsymbol\alpha=(\alpha_1,\ldots,\alpha_d)\in\Nzero^d\),
we use
\begin{equation}\label{eq:multi-index-conventions}
 |\boldsymbol\alpha|=\sum_{j=1}^d\alpha_j,\qquad
 \boldsymbol\alpha!=\prod_{j=1}^d\alpha_j!,\qquad
 \boldsymbol z^{\boldsymbol\alpha}=\prod_{j=1}^d z_j^{\alpha_j}.
\end{equation}
We write \(\R[\boldsymbol z]\)
for the polynomial ring and \(\R[[\boldsymbol z]]\) for the ring of
formal power series, and use
\begin{equation}\label{eq:nonnegative-polynomial-ring}
 \R_{\ge0}[\boldsymbol z]
 =\left\{\sum_{\boldsymbol\alpha}c_{\boldsymbol\alpha}
       \boldsymbol z^{\boldsymbol\alpha}:
       c_{\boldsymbol\alpha}\ge0,\quad
       \text{only finitely many }c_{\boldsymbol\alpha}\ne0\right\}.
\end{equation}
For polynomials or formal power series \(P,Q\) in the same
indeterminates, we write \(P\succeq_{\mathrm{coeff}}Q\) when every
coefficient of \(P-Q\) is nonnegative. In particular,
\(P\succeq_{\mathrm{coeff}}0\) for a polynomial means membership in
\eqref{eq:nonnegative-polynomial-ring}.
When we assert coefficientwise positivity, we require this condition
for every minor under consideration.

Throughout, a matrix is \emph{\textbf{totally positive}} (\(\TP\)) if its every
minor is nonnegative and \emph{\textbf{coefficientwise totally positive}} if its every
minor is a nonnegative coefficient polynomial. For a sequence of polynomials \(a_n(\boldsymbol z)\),
it is called \emph{\textbf{coefficientwise Hankel-totally positive}} if the Hankel matrix \([a_{i+j}(\boldsymbol z)]_{i,j\ge0}\) is coefficientwise totally positive. In other words, for every \(r\ge1\) and strictly
increasing lists \(I=(i_1<\cdots<i_r)\), \(J=(j_1<\cdots<j_r)\),
\begin{equation}\label{eq:display-002}
 \det[a_{i_a+j_b}(\boldsymbol z)]_{a,b=1}^r
       \in\R_{\ge0}[\boldsymbol z].
\end{equation}
For a real sequence, Hankel total positivity is equivalent to being a
Stieltjes moment sequence. Coefficientwise positivity implies positivity at every nonnegative
parameter value $\boldsymbol z$. For a general polynomial, the converse fails:
\((u-1)^2\ge0\) for every \(u\ge0\), although its coefficient of
\(u\) is negative. The distinction is therefore essential when we
pass from numerical moment sequences to combinatorial polynomials.

We write \(S_n\) for the permutations of \([1,n]\), with
\(\operatorname{sgn}(\sigma)\in\{-1,1\}\) their signs.
For a set \(A\), we use \(\mathbf1_A\) for its indicator, and
\(\delta_x\) for the unit point mass at \(x\).
In a finite-dimensional coordinate space, we write
\(\boldsymbol1=(1,\ldots,1)\) and use \(\boldsymbol e_j\) for the
vector whose \(j\)-th entry is one and whose other entries are zero.

For a matrix \(M\) and increasing index lists
\(I=(i_1,\ldots,i_r)\), \(J=(j_1,\ldots,j_r)\), let
\begin{equation}\label{eq:matrix-submatrix-notation}
 M[I,J]:=[M_{i_a,j_b}]_{a,b=1}^r.
\end{equation}
We use \(M^{\mathsf T}\) for the transpose, \(I\) for an identity
matrix of the required size, and \(E_{i,j}\) for the matrix with
its only nonzero entry, equal to one, in position \((i,j)\).
We use \(\diag\) for a diagonal or block-diagonal matrix and
\(\adj(M)\) for the classical adjugate, the transpose of the cofactor
matrix. Unless a finite range is specified, we index infinite matrices
by \(\Nzero\). Empty products and determinants of size zero have value one.

For \(a\in\R\) and \(n\in\Nzero\), our factorial conventions are
\begin{equation}\label{eq:factorial-conventions}
 (a)_n:=\prod_{h=0}^{n-1}(a+h),\qquad
 a^{\underline n}:=\prod_{h=0}^{n-1}(a-h),\qquad
 \binom an:=\frac{a^{\underline n}}{n!}.
\end{equation}

For a positive measure \(\mu\), we use interval limits for
one-dimensional integrals. Our endpoint convention is
\begin{equation}\label{eq:measure-endpoint-convention}
 \int_a^b f(x)\,\mu(\dif x):=\int_{[a,b]}f(x)\,\mu(\dif x),\qquad
 \int_{a^+}^b f(x)\,\mu(\dif x):=\int_{(a,b]}f(x)\,\mu(\dif x).
\end{equation}
At an infinite endpoint we integrate over the corresponding half-line.
Thus \(\int_0^\infty f(x)\,\mu(\dif x)\) includes a possible atom at
zero. We extend measures originally defined on \((0,\infty)\) or
\(\R\setminus\{0\}\) by assigning zero mass to the omitted points.
We retain a set below the integral sign for a multidimensional domain
or a general measurable space.

\subsection{Organization of the paper}

In Section~\ref{sec:lifting}, we collect the common moment and Poisson
tools and prove the Stieltjes lifting criterion in
Theorem~\ref{thm:lifting}. In Section~\ref{sec:hamburger-proof}, we
prove the Hamburger lifting criterion by using the real-line
Vandermonde expansion. In Section~\ref{sec:universal}, we establish
binomial-convolution, composition, and substitution results, discuss
logarithms of ordinary Stieltjes series, and collect moment inputs
for the applications.

In Section~\ref{sec:digraph}, we apply the lifting criteria to
functional digraphs and their extensions. In
Section~\ref{sec:equations}, we derive general results for implicit
equations and autonomous differential equations from
P\'olya-frequency series. In Section~\ref{sec:tree-applications},
we use known moment representations for tree and forest sequences
to verify the lifting hypotheses and obtain families that are
coefficientwise Hankel-totally positive in the lifting variables.
We treat the Ward sequence and several related tree families in
this way.

In Section~\ref{sec:applications}, we obtain further general families
from polynomial moment inputs, including extensions of the Bessel,
Lah, and Laguerre polynomials and a fixed-shape result for the
Schl\"afli-Gessel-Seo polynomials. We also consider set partitions,
lognormal moments, Poisson differences, and matrix consequences.
For each application, we distinguish the fixed parameters used to
verify the moment condition from the variables in which we prove
coefficientwise positivity.  We not only give new proofs of coefficientwise Hankel-total positivity of many combinatorial polynomials in a unified approach, but also prove coefficientwise Hankel-total positivity for the component polynomials of functional
digraphs confirming a conjecture of Sokal and establish fixed-shape Schl\"afli-Gessel-Seo polynomials positivity supporting another coefficientwise Hankel-total positivity conjecture of Sokal. In Section~\ref{sec:scope}, we discuss
the remaining questions about positivity jointly in the structural
parameters.

\section{Proof of the Stieltjes lifting theorem}\label{sec:lifting}

In this section, we will prove Theorem~\ref{thm:lifting}. We first collect the moment
criteria, alternant identities, and Poisson formulas used. For the Stieltjes sufficiency of Theorem~\ref{thm:lifting}, we construct square-integral
expansions with variable-independent kernels and integrate them
against an arbitrary initial moment measure. The lowest-degree
terms of shifted Hankel minors then give necessity.

By the classical exponential formula of
Stanley~\cite[Corollary~5.1.6]{Stanley}, we may express each moment
polynomial as a finite sum over set partitions. Let \(\cP_0=\{\varnothing\}\) and \(\cP_n\) denote
the set of partitions of \([n]\). With block weight
\(\sum_\ell y_\ell\kappa_{|B|}^{(\ell)}\), we obtain
\begin{equation}\label{eq:bell-partitions}
 M_n(\by)=
 \sum_{\pi\in\cP_n}\prod_{B\in\pi}
       \left(\sum_{\ell=1}^{m}y_\ell\kappa_{|B|}^{(\ell)}\right).
\end{equation}
Thus \(M_0=1\) and \(\deg_{\by}M_n\le n\).

\subsection{Common tools}

We begin with the classical Hamburger and Stieltjes moment criteria (see Schm\"udgen~\cite[Chapter~3]{Schmudgen2017},
Pinkus~\cite{Pinkus2010} and Karlin~\cite{Kar68} for instance).

\begin{lemma}[Classical moment criteria]\label{lem:moment-criteria}
Let \((a_n)_{n\ge0}\) be a real sequence, and, for \(h\in\{0,1\}\), define
\begin{equation}\label{eq:finite-moment-matrices}
 H_d^{(h)}(\boldsymbol a):=[a_{i+j+h}]_{i,j=0}^{d}
 \qquad(d\ge0).
\end{equation}
\begin{enumerate}
\renewcommand{\labelenumi}{\textup{(\roman{enumi})}}
\item The sequence is a Hamburger moment sequence if and only if
\(H_d^{(0)}(\boldsymbol a)\) is positive semidefinite for every \(d\ge0\).
\item The sequence is a Stieltjes moment sequence if and only if both
\(H_d^{(0)}(\boldsymbol a)\) and \(H_d^{(1)}(\boldsymbol a)\) are positive
semidefinite for every \(d\ge0\). Equivalently, its infinite Hankel
matrix \([a_{i+j}]_{i,j\ge0}\) is totally positive.
\item If a representing Stieltjes measure has infinitely many positive
support points, then every Hankel minor with strictly increasing row
and column index lists is strictly positive.
\end{enumerate}
\end{lemma}

To convert a determinant of integrals into an integral of two
determinants, we use Andr\'eief's integration identity~\cite{Andreief1886}.
We integrate over the full Cartesian product, so we include the
normalizing factor \(1/r!\).

\begin{lemma}[Andr\'eief's integration identity]\label{lem:andreief-source}
Let \(\mu\) be a positive measure on a measurable space \(X\), and let
\(f_1(t),\ldots,f_r(t)\) and \(g_1(t),\ldots,g_r(t)\) be measurable
functions on \(X\) such that
\(\int_X |f_i(t)g_j(t)|\,d\mu(t)<\infty\) for all \(i,j\). Then we have
\begin{equation}\label{eq:andreief-identity}
 \det\left[\int_X f_i(t)g_j(t)\,d\mu(t)\right]_{i,j=1}^r
 =\frac1{r!}\int_{X^r}
      \det[f_i(t_j)]_{i,j=1}^r
      \det[g_i(t_j)]_{i,j=1}^r\prod_{j=1}^r d\mu(t_j).
\end{equation}
\end{lemma}

We next recall the notation for Schur polynomials. For a partition
\(\lambda=(\lambda_1\ge\cdots\ge\lambda_r\ge0)\), we draw
\(\lambda_i\) boxes in row \(i\). Let
\(\operatorname{SSYT}_r(\lambda)\) be the set of fillings of these boxes
by \(1,\ldots,r\) that are weakly increasing along rows and strictly
increasing down columns. Let \(m_j(P)\) be the number of entries equal
to \(j\) in a filling \(P\), and we define
\begin{equation}\label{eq:schur-tableau-definition}
 s_\lambda(t_1,\ldots,t_r)
 =\sum_{P\in\operatorname{SSYT}_r(\lambda)}
                  \prod_{j=1}^r t_j^{m_j(P)}.
\end{equation}
For the empty partition, we set \(s_{\varnothing}=1\).
From \eqref{eq:schur-tableau-definition}, we have nonnegative integer
coefficients. To factor the generalized Vandermonde determinants,
we use the alternant formula in
Macdonald~\cite[Chapter~I, Sections~3 and~5]{Macdonald1995}, or
Stanley~\cite[Sections~7.10 and~7.15]{Stanley}.

\begin{lemma}[Schur alternants and coefficient positivity]\label{lem:schur-source}
For integers \(0\le n_1<\cdots<n_r\), set
\(\lambda_i=n_{r+1-i}-(r-i)\) with \(1\le i\le r\). Then we have
\begin{equation}\label{eq:schur-alternant}
 \det[t_j^{n_i}]_{i,j=1}^r
 =\Delta(\boldsymbol t)s_\lambda(\boldsymbol t),
 \qquad
 \Delta(\boldsymbol t)=\prod_{1\le i<j\le r}(t_j-t_i).
\end{equation}
In particular, the determinant in \eqref{eq:schur-alternant} is positive
when \(0<t_1<\cdots<t_r\).
\end{lemma}

For the subsequent matrix factorizations, we record three standard
closure rules. We cite the finite matrix results in the corresponding
parts.

\begin{lemma}[Classical matrix closure rules]\label{lem:matrix-tools}
\begin{enumerate} The following hold.
\renewcommand{\labelenumi}{\textup{(\roman{enumi})}}
\item \textup{(Finite products.)}
The product of coefficientwise totally positive matrices is coefficientwise
totally positive whenever the Cauchy--Binet sums for its minors are finite \cite[Proposition~1.4]{Pinkus2010}.
\item \textup{(Bidiagonal matrices.)}
An upper or lower bidiagonal matrix with coefficientwise nonnegative
entries is coefficientwise totally positive \cite[Lemma~9.1]{PSZ2023}.
\end{enumerate}
\end{lemma}

\subsection{Stieltjes representations and compound-Poisson moments}

We will use the Poisson exponential formula
\cite[Theorem~3.9]{LastPenrose2017} and the Fock isometry
\cite[Theorem~1.1]{LastPenrose2011}. We allow atomic intensity
measures. For the exponential formula, we use bounded integrands.
For the Fock isometry, we use polynomial functions of finite Poisson
counts with bounded jumps; consequently, all required polynomial
moments and second moments are finite.

\begin{lemma}[Poisson exponential formula \cite{LastPenrose2017} and Fock isometry \cite{LastPenrose2011}]
\label{lem:poisson-source}
Let \(\eta\) be a Poisson random measure on a measurable space
\(E\), with finite intensity measure \(\Lambda\).
\begin{enumerate}
\renewcommand{\labelenumi}{\textup{(\roman{enumi})}}
\item For a bounded measurable function \(f:E\to\R\),
\begin{equation}\label{eq:poisson-exponential-source}
 \E\exp\left(\int_E f(z)\,\eta(\dif z)\right)
       =\exp\left(\int_E(e^{f(z)}-1)\,\Lambda(\dif z)\right).
\end{equation}
\item For a square-integrable real function \(F\) of the point measure \(\eta\), define
\(D_zF(\eta)=F(\eta+\delta_z)-F(\eta)\). Then we have
\begin{equation}\label{eq:fock-isometry-source}
 \E\bigl[(F(\eta))^2\bigr]
 =\sum_{k\ge0}\frac1{k!}
  \int_{E^k}\bigl(\E[D_{z_1}\cdots D_{z_k}F(\eta)]\bigr)^2
                         \,\prod_{h=1}^k\Lambda(\dif z_h).
\end{equation}
For \(k=0\), the summand is \((\E F)^2\).
\end{enumerate}
\end{lemma}

To handle products of Poisson integrals, we use the multivariate
Mecke equation~\cite[Theorem~4.4]{LastPenrose2017}.
We distinguish point occurrences from their locations: two distinct
occurrences may be situated at the same location. By distinguishing occurrences from locations, we can work with
atomic intensity measures throughout.

\begin{lemma}[Factorial Mecke identity \cite{LastPenrose2017}]\label{lem:mecke-source}
For a Poisson random measure \(\eta\) on a measurable space \(E\)
with intensity measure \(\Lambda\), let
\(\eta^{(k)}\) count ordered \(k\)-tuples of distinct point occurrences.
For every nonnegative measurable function \(H(\eta,\boldsymbol z)\), we have
\begin{equation}\label{eq:factorial-mecke-source}
 \E\int_{E^k} H(\eta,\boldsymbol z)\,d\eta^{(k)}(\boldsymbol z)
 =\int_{E^k} \E H\left(\eta+\sum_{i=1}^k\delta_{z_i},
                       \boldsymbol z\right)\,\prod_{h=1}^k\Lambda(\dif z_h).
\end{equation}
\end{lemma}

Before applying the Poisson exponential formula
\eqref{eq:poisson-exponential-source}, the Fock isometry
\eqref{eq:fock-isometry-source}, and the factorial Mecke identity
\eqref{eq:factorial-mecke-source} to the cumulant families, we first
separate the mass at zero in each representing measure.
For \(x\ge0\), we write \(1\wedge x:=\min\{1,x\}\). We obtain the following
moment form of the classical subordinator representation; for the
corresponding Laplace exponents, see
\cite[Theorems~3.2 and~5.2]{SSV2012}. We include the elementary
measure conversion because we shall use the individual formal moments.

\Needspace{12\baselineskip}
\begin{lemma}[L\'evy cumulants in moment form \cite{SSV2012}]\label{lem:representation}
Assume that the shifted sequence
\((\kappa_{n+1}^{(\ell)})_{n\ge0}\) is a Stieltjes moment sequence for each \(\ell=1,\ldots,m\) . Then there exist
\(b_\ell\ge0\) and nonnegative measures \(\rho_\ell\) on \((0,\infty)\)
such that, for every \(\ell\),
\begin{equation}\label{eq:levy-cumulants}
 \kappa_n^{(\ell)}
 =b_\ell\mathbf1_{\{n=1\}}
    +\int_0^\infty x^n\,\rho_\ell(\dif x)
 \qquad(n\ge1).
\end{equation}
All moments in \eqref{eq:levy-cumulants} are finite, and
\begin{equation}\label{eq:levy-integrability}
 \int_0^\infty (1\wedge x)\,\rho_\ell(\dif x)<\infty.
\end{equation}
\end{lemma}

\begin{proof}
For each \(\ell\), choose a nonnegative representing measure
\(\nu_\ell\) on \([0,\infty)\) such that
\begin{equation}\label{eq:shifted-cumulant-measure}
 \kappa_{n+1}^{(\ell)}
 =\int_0^\infty x^n\,\nu_\ell(\dif x)
 \qquad(n\ge0).
\end{equation}
Taking \(n=0\) in \eqref{eq:shifted-cumulant-measure}, we obtain
\(\nu_\ell([0,\infty))=\kappa_1^{(\ell)}<\infty\). Define
\begin{equation}\label{eq:drift-and-jump-measure}
 b_\ell=\nu_\ell(\{0\}),\qquad
 \rho_\ell(\dif x)=x^{-1}\nu_\ell(\dif x)\quad(x>0).
\end{equation}
For \(n=1\), we recover \eqref{eq:levy-cumulants} by adding the mass
at zero to \(\int_0^\infty  x\,d\rho_\ell=\nu_\ell((0,\infty))\).
On the other hand, for \(n\ge2\), the atom at zero contributes nothing to the
\((n-1)\)-st moment of \(\nu_\ell\). Hence
\(\int_0^\infty  x^n\,d\rho_\ell=\int_0^\infty  x^{n-1}\,d\nu_\ell
=\kappa_n^{(\ell)}\). Finally, by \eqref{eq:drift-and-jump-measure}, we have
\begin{equation}\label{eq:levy-integrability-calculation}
 \int_0^1 x\,\rho_\ell(\dif x) =\nu_\ell((0,1])<\infty,\qquad \int_{1^+}^\infty \rho_\ell(\dif x) =\int_{1^+}^\infty x^{-1}\nu_\ell(\dif x) \le\nu_\ell((1,\infty))<\infty.
\end{equation}
Adding such two finite integrals in \eqref{eq:levy-integrability-calculation},
we obtain \eqref{eq:levy-integrability}.
\end{proof}

Until the truncation step, we assume that each \(\rho_\ell\) is finite
and compactly supported in \((0,\infty)\). For the probabilistic
calculation we fix \(\by\in(0,\infty)^m\) and define
\begin{equation}\label{eq:combined-intensity}
 \rho_{\by}=\sum_{\ell=1}^{m}y_\ell\rho_\ell,
 \qquad b(\by)=\sum_{\ell=1}^{m}b_\ell y_\ell.
\end{equation}
We write \(N_{\by}\) for a Poisson random measure of 
\(\rho_{\by}\), and set
\begin{equation}\label{eq:compound-poisson-variable}
 J_{\by}=\int_0^\infty x\,N_{\by}(\dif x),
 \qquad X_{\by}=b(\by)+J_{\by}.
\end{equation}

\begin{lemma}[Compound-Poisson moment representation]
\label{lem:compound-poisson-moments}
Let \(b_\ell\ge0\) and \(\rho_\ell\) in
Lemma~\ref{lem:representation}, and suppose that each
\(\rho_\ell\) is finite and compactly supported in \((0,\infty)\).
For fixed \(\by\in(0,\infty)^m\), let \(X_{\by}\) be defined by
\eqref{eq:combined-intensity}--\eqref{eq:compound-poisson-variable}.
Then \(X_{\by}\ge0\), its moment generating function is finite
for every real \(t\), and
\(M_n(\by)=\E X_{\by}^{\,n}\) for every \(n\ge0\).
In particular, \((M_n(\by))_{n\ge0}\) is a Stieltjes moment sequence.
\end{lemma}

\begin{proof}
For clarity, we will derive the moment generating function explicitly.
Put \(\lambda_{\by}=\rho_{\by}((0,\infty))\). If
\(\lambda_{\by}>0\), we may write \(J_{\by}=\sum_{j=1}^{P}Y_j\),
where \(P\) has the Poisson distribution of mean \(\lambda_{\by}\),
and the \(Y_j\) are independent with common law
\(\rho_{\by}/\lambda_{\by}\), independently of \(P\).
Conditioning on \(P\), and then summing its probability mass function,
we obtain
\begin{equation}\label{eq:compound-poisson-conditioning}
 \E\,e^{tJ_{\by}} =\sum_{k\ge0}e^{-\lambda_{\by}} \frac{\lambda_{\by}^k}{k!} \left(\frac1{\lambda_{\by}} \int_0^\infty e^{tx}\,\rho_{\by}(\dif x)\right)^k =\exp\left(\int_0^\infty (e^{tx}-1) \,\rho_{\by}(\dif x)\right).
\end{equation}
If \(\lambda_{\by}=0\), we take \(J_{\by}=0\), and both sides of
\eqref{eq:compound-poisson-conditioning} are one.
For bounded jump sizes, we may expand \(e^{tx}-1\) under the finite
integral. By \eqref{eq:levy-cumulants} and
\eqref{eq:combined-intensity}, we have
\begin{equation}\label{eq:compound-poisson-cumulant-expansion}
 b(\by)t+\int_0^\infty (e^{tx}-1)\,\rho_{\by}(\dif x) =\sum_{\ell=1}^m y_\ell \left(b_\ell t+ \sum_{n\ge1}\frac{t^n}{n!} \int_0^\infty x^n\,\rho_\ell(\dif x)\right) =\sum_{\ell=1}^m y_\ell K_\ell(t).
\end{equation}
Multiplying \eqref{eq:compound-poisson-conditioning} by
\(e^{b(\by)t}\), and using
\eqref{eq:compound-poisson-cumulant-expansion}, we obtain
\begin{equation}\label{eq:poisson-mgf}
 \E \,e^{tX_{\by}}=\E\, e^{b(\by)t+tJ_{\by}}
 =\exp\left(b(\by)t+
          \int_0^\infty (e^{tx}-1)\,\rho_{\by}(\dif x)\right)
 =\exp\left(\sum_{\ell=1}^m y_\ell K_\ell(t)\right).
\end{equation}
We may differentiate at zero because the moment generating function
is finite for every real \(t\). Comparing \eqref{eq:poisson-mgf}
with \eqref{eq:bell-egf}, we conclude that
\begin{equation}\label{eq:poisson-moments}
 M_n(\by)=\E\, X_{\by}^{\,n}\qquad(n\ge0).
\end{equation}
\end{proof}

We will use only finite polynomial identities in \(\by\). Once we
have proved such an identity on $\by>0$, we may therefore
extend it to all boundary specializations as an identity in the polynomial ring \(\R[\by]\).

\subsection{A multivariate shifted Vandermonde square}

For \(\bx=(x_1,\ldots,x_r)\), we continue to use
\(\Delta(\bx)=\prod_{i<j}(x_j-x_i)\), as in
\eqref{eq:schur-alternant}. Take independent copies
\(J_{\by}^{(1)},\ldots,J_{\by}^{(r)}\) of \(J_{\by}\), and write
\(\bJ_{\by}=(J_{\by}^{(1)},\ldots,J_{\by}^{(r)})\).
We shall obtain nonnegative coefficients by expressing the
second moment of the Vandermonde as a finite sum of integrals of squares.

\begin{lemma}[Multivariate Poisson--Vandermonde positivity]\label{lem:vandermonde}
For every fixed \(\bs\in\R^r\), define
\begin{equation}\label{eq:V-definition}
 V_{\bs}(\by):=\E\!\left[
                         \bigl(\Delta(\bJ_{\by}+\bs)\bigr)^2\right].
\end{equation}
Then we have \(V_{\bs}(\by)\in\R_{\ge0}[\by]\) and
\(\deg_{\by}V_{\bs}\le D:=\binom r2\).
Every coefficient is a polynomial in \(\bs\) of degree at most \(2D\).
The same assertions hold for finite compactly supported nonnegative
jump measures on \(\R\setminus\{0\}\).
\end{lemma}

\begin{proof}
To remove the  indeterminates from the kernels in the Fock expansion,
we first compute the mean of the Vandermonde. For this purpose,
we establish the following mean identity.

\begin{claim}\label{claim:vandermonde-mean}
For every \(\bs\in\R^r\), we have
\begin{equation}\label{eq:vandermonde-mean}
 \E[\Delta(\bJ_{\by}+\bs)]=\Delta(\bs).
\end{equation}
\end{claim}
\begin{proof}
Choose \(B<\infty\) that bounds the absolute jump sizes, and let
\(P_j\) be the Poisson count in coordinate \(j\). Then
\(|J_{\by}^{(j)}|\le BP_j\). For \(k\ge0\), define
\begin{equation}\label{eq:shifted-moment-polynomial}
 p_k(u):=\E(J_{\by}+u)^k
       =\sum_{h=0}^k\binom kh(\E\, J_{\by}^{\,h})u^{k-h}.
\end{equation}
Since a Poisson count has finite moments
of every order, we may integrate every polynomial above. In particular, \(p_k(u)\) is monic of degree \(k\), because
\(\E\, J_{\by}^{0}=1\).

Expanding the determinant and then using independence in each product,
we obtain
\begin{equation}\label{eq:independent-vandermonde-expectation}
 \E\,\Delta(\bJ_{\by}+\bs) =\sum_{\sigma\in S_r}\operatorname{sgn}(\sigma) \prod_{i=1}^r\E(s_{\sigma(i)}+ J_{\by}^{(\sigma(i))})^{i-1} =\det[p_{i-1}(s_j)]_{i,j=1}^r.
\end{equation}
Here we use independence because the indices \(\sigma(i)\) are distinct.

Write \(A_{\by}(\bs)=[p_{i-1}(s_j)]_{i,j=1}^r\) and
\(V(\bs)=[s_j^{i-1}]_{i,j=1}^r\). By
\eqref{eq:shifted-moment-polynomial}, we have
\begin{equation}\label{eq:monic-change-of-basis}
 A_{\by}(\bs)=L_{\by}V(\bs),\qquad
 (L_{\by})_{i,h}=
 \begin{cases}
 \binom{i-1}{h-1}\E J_{\by}^{\,i-h},&1\le h\le i\le r,\\
 0,&1\le i<h\le r.
 \end{cases}
\end{equation}
Since \(L_{\by}\) is lower triangular with diagonal entries equal to one,
we obtain
\begin{equation}\label{eq:monic-vandermonde-reduction}
 \det A_{\by}(\bs)=\det L_{\by}\det V(\bs)
                 =\Delta(\bs).
\end{equation}
Combining \eqref{eq:independent-vandermonde-expectation} and
\eqref{eq:monic-vandermonde-reduction} yields
\eqref{eq:vandermonde-mean}.
\end{proof}

Using the mean identity, we now express the second moment as
an integral of squares with kernels independent of \(\by\).
We use this expansion to establish coefficientwise nonnegativity
and the degree bound in the following claim.

\begin{claim}\label{claim:vandermonde-square}
For every \(\bs\in\R^r\) and \(\boldsymbol\alpha\in\Nzero^m\), we have
\begin{equation}\label{eq:claim-vandermonde-coefficients}
 [\by^{\boldsymbol\alpha}]V_{\bs}(\by)\ge0\quad(|\boldsymbol\alpha|\le D),\qquad [\by^{\boldsymbol\alpha}]V_{\bs}(\by)=0\quad(|\boldsymbol\alpha|>D).
\end{equation}
\end{claim}
\begin{proof}
\emph{Step 1: One point measure for all coordinates.}
Take \(E=\{1,\ldots,r\}\times(0,\infty)\). For Borel
\(A\subseteq(0,\infty)\), define
\begin{equation}\label{eq:coordinate-intensity-measures}
 \Lambda_\ell(\{j\}\times A)=\rho_\ell(A),\qquad
 \Lambda_{\by}=\sum_{\ell=1}^m y_\ell\Lambda_\ell.
\end{equation}
Let \(\eta\) be a Poisson random measure with intensity \(\Lambda_{\by}\).
The restrictions of \(\eta\) to the coordinate copies are independent;
hence we may realize \(J_{\by}^{(j)}\) as the jump sum on copy \(j\).
For a function \(f:\R^r\to\R\) and \(z=(j,x)\in E\), put
\begin{equation}\label{eq:translation-difference}
 \Dshift_{(j,x)}f(\bs)=f(\bs+x\boldsymbol e_j)-f(\bs).
\end{equation}

\smallskip\noindent
\emph{Step 2: Expected difference kernels.}
Set \(F(\eta)=\Delta(\bJ_{\by}+\bs)\). On adding an occurrence
\((j,x)\), we increase only the \(j\)-th jump sum by \(x\).
Therefore, iterating \eqref{eq:translation-difference} and then taking
expectations, we obtain
\begin{equation}\label{eq:mean-difference-kernel}
 \E[D_{z_1}\cdots D_{z_k}F(\eta)] =\Dshift_{z_1}\cdots\Dshift_{z_k} \E[\Delta(\bJ_{\by}+\bs)] =\Dshift_{z_1}\cdots\Dshift_{z_k}\Delta(\bs).
\end{equation}
For the first equality, we interchange expectation with a finite signed
sum of integrable polynomials. For the second, we use Claim~\ref{claim:vandermonde-mean} at every translated shift. We have thus
obtained kernels that do not depend on \(\by\).

\smallskip\noindent
\emph{Step 3: A finite expansion into integrals of squares.}
For fixed jump vectors, a nonzero translation difference reduces the
degree in \(\bs\) by at least one. Since \(\Delta\) has degree \(D\),
the last expression in \eqref{eq:mean-difference-kernel} vanishes for
\(k>D\). Write \(P_j\) for the Poisson count on copy \(j\), and
choose \(B\) that bounds the absolute jump sizes. Then, for a constant
\(C_{\bs,B}\),
\begin{equation}\label{eq:vandermonde-square-integrability}
 |F(\eta)|^2\le C_{\bs,B}
                   \left(1+\sum_{j=1}^r P_j\right)^{2D},
 \qquad \E|F(\eta)|^2<\infty.
\end{equation}
Applying the Fock isometry \eqref{eq:fock-isometry-source} to
\(F(\eta)=\Delta(\bJ_{\by}+\bs)\), and then using \eqref{eq:mean-difference-kernel}, we have
\begin{equation}\label{eq:V-square-integrals}
 V_{\bs}(\by)
 =\sum_{k=0}^{D}\frac1{k!}\int_{E^k}
       \left(\Dshift_{z_1}\cdots\Dshift_{z_k}\Delta(\bs)\right)^2
       \prod_{h=1}^{k}\Lambda_{\by}(\dif z_h).
\end{equation}
For \(k=0\), we interpret the integral as \(\Delta(\bs)^2\).

\smallskip\noindent
\emph{Step 4: Coefficient extraction and polynomial dependence on the shifts.}
Let \(|\boldsymbol\alpha|=k\le D\) and
\(\boldsymbol\alpha!=\prod_{\ell=1}^m\alpha_\ell!\).
Choose an ordered color list \((c_1,\ldots,c_k)\) containing
\(\alpha_\ell\) copies of \(\ell\), with all copies of color \(1\)
first, then color \(2\), and so on. We integrate \(z_h\) against
\(\Lambda_{c_h}\).
Applying the expansion \eqref{eq:coordinate-intensity-measures}
to \eqref{eq:V-square-integrals}, we have
\begin{equation}\label{eq:V-coefficient-squares}
 [\by^{\boldsymbol\alpha}]V_{\bs}(\by)
 =\frac1{\boldsymbol\alpha!}
   \int_{E^k}
    \left(\Dshift_{z_1}\cdots\Dshift_{z_k}\Delta(\bs)\right)^2
    \prod_{h=1}^k\Lambda_{c_h}(\dif z_h)
 \ge0.
\end{equation}
Indeed, the difference operators commute, so the kernel is symmetric
in \(z_1,\ldots,z_k\). The \(k!/\boldsymbol\alpha!\) assignments of
the colors give equal integrals; after multiplication by \(1/k!\),
we obtain \(1/\boldsymbol\alpha!\).
No term with \(|\boldsymbol\alpha|>D\) occurs in
\eqref{eq:V-square-integrals}. Finally, the squared kernel has degree
at most \(2(D-k)\) in \(\bs\). Since the measures are finite with
bounded jump sizes, we may integrate its finitely many polynomial
coefficients. Thus every coefficient of \(V_{\bs}\) is a polynomial
in \(\bs\) of degree at most \(2D\).
\end{proof}

Combining such two claims, we prove the lemma for nonnegative jump
locations. We used no sign restriction on the shifts or on the jump
locations in the determinant and square-integral calculations.
Applying the same calculations to
\(E=\{1,\ldots,r\}\times(\R\setminus\{0\})\), we have the real-line
extension as well.
\end{proof}

\subsection{A partitioned form of the Mecke formula}

For the next step, we need a power of a Poisson integral multiplied by
a function of the same process. We record the following standard
consequence of the multivariate Mecke equation
\cite[Theorem~4.4]{LastPenrose2017}, recalled in
Lemma~\ref{lem:mecke-source}. We obtain it by grouping the factors in
the power according to set partitions. We include the complete
reduction to specify the treatment of repeated locations and the
absence of any additional factorial.

\Needspace{12\baselineskip}
\begin{lemma}[Partitioned consequence of the Mecke formula]\label{lem:mecke}
Let \(N\) be a Poisson random measure on \((0,\infty)\) with finite,
compactly supported measure \(\theta\), and put
\(J=\int_0^\infty x\,N(\dif x)\). For \(p\in\Nzero\) and a nonnegative measurable function \(G\)
on the space of finite point measures on \((0,\infty)\), we have
\begin{equation}\label{eq:partitioned-mecke}
 \E\,[J^pG(N)]
 =\sum_{\pi\in\cP_p}
   \int_{(0,\infty)^{|\pi|}}
    \left(\prod_{B\in\pi}x_B^{|B|}\right)
    \E\, G\!\left(N+\sum_{B\in\pi}\delta_{x_B}\right)
    \prod_{B\in\pi}\theta(\dif x_B).
\end{equation}
For \(p=0\), we use the
single empty-partition term. The identity holds in \([0,\infty]\).
\end{lemma}

\begin{proof}
\emph{Step 1: Grouping the factors by distinct occurrences.}
We first prove a deterministic identity for a finite point measure.
Write \(\eta=\sum_{j=1}^M\delta_{x_j}\), where the indices
\(j\) distinguish occurrences even when some locations \(x_j\) coincide.
For \(p\ge1\), we expand the power and obtain
\begin{equation}\label{eq:ordinary-point-tuple-expansion}
 \left(\int_0^\infty  x\,\eta(\dif x)\right)^p
   =\sum_{a_1,\ldots,a_p=1}^M x_{a_1}\cdots x_{a_p}.
\end{equation}
For a tuple \((a_1,\ldots,a_p)\), place two positions \(u,v\in[p]\)
in the same block exactly when \(a_u=a_v\). In this way we associate
one partition \(\pi\in\cP_p\) with the tuple. Order the blocks by
their least elements, say \(B_1,\ldots,B_k\). For a fixed partition,
we must choose \(k\) distinct occurrence indices \(j_1,\ldots,j_k\),
one for each ordered block. The product in
\eqref{eq:ordinary-point-tuple-expansion} is then
\(\prod_{h=1}^k x_{j_h}^{|B_h|}\). Conversely, from the partition and the ordered distinct indices,
we recover exactly one original tuple.
Therefore, using the factorial measure of distinct occurrences, we deduce
\begin{equation}\label{eq:power-factorial-partition}
 \left(\int_0^\infty x\,\eta(\dif x)\right)^p
 =\sum_{\pi\in\cP_p}
   \int_{(0,\infty)^{|\pi|}}\left(\prod_{B\in\pi}x_B^{|B|}\right)
                      \eta^{(|\pi|)}(\dif\boldsymbol x_\pi).
\end{equation}

\smallskip\noindent
\emph{Step 2: The factorial Mecke identity for one partition.}
Because \(\theta\) is finite, the Poisson random measure \(N\) is
finite almost surely. Applying \eqref{eq:power-factorial-partition}
to \(N\), and then multiplying by \(G(N)\), we may take expectations
term by term. For a partition \(\pi=\{B_1,\ldots,B_k\}\), with the
blocks in their canonical order, write \(\bx=(x_1,\ldots,x_k)\) and set
\begin{equation}\label{eq:partitioned-mecke-integrand}
 w_\pi(\bx)=\prod_{h=1}^k x_h^{|B_h|},\qquad
 g_k(\bx)=\E\, G\left(N+\sum_{h=1}^k\delta_{x_h}\right).
\end{equation}
Thus \(w_\pi(\bx)\) is the product of the block weights, and
\(g_k(\bx)\) is the expected value of \(G\) after inserting one
occurrence at each of \(x_1,\ldots,x_k\). Applying the factorial
Mecke identity \eqref{eq:factorial-mecke-source} to
\(H_\pi(\eta,\bx)=G(\eta)w_\pi(\bx)\), we obtain
\begin{equation}\label{eq:partitioned-mecke-one-term}
 \E\int_{(0,\infty)^k}G(N)w_\pi(\bx)\,N^{(k)}(\dif\bx)
 =\int_{(0,\infty)^k}w_\pi(\bx)g_k(\bx)\,\prod_{h=1}^k\theta(\dif x_h).
\end{equation}
\smallskip\noindent
\emph{Step 3: Summation and the empty partition.}
Summing \eqref{eq:partitioned-mecke-one-term} over all partitions,
we recover \eqref{eq:partitioned-mecke}. We require no factor \(k!\):
we fixed one block order, and we already count all ordered distinct
occurrence choices through the factorial measure. Moreover, since
the integrands are nonnegative, we may interchange expectation
and integration by Tonelli's theorem without assuming that
\(\E\, G(N)<\infty\). For \(p=0\), both sides are
\(\E\, G(N)\) under the empty-product convention.

For example, when \(p=2\), we have exactly two partitions. Therefore,
using \(g_1(x)=\E\, G(N+\delta_x)\) and
\(g_2(x,y)=\E\, G(N+\delta_x+\delta_y)\), we obtain
\begin{equation}\label{eq:partitioned-mecke-p2}
 \E\,[J^2G(N)]
 =\int_0^\infty x^2g_1(x)\,\theta(\dif x)
 +\int_{(0,\infty)^2}xy\,g_2(x,y)\,\theta(\dif x)\theta(\dif y).
\end{equation}
In the second integral we also include \(x=y\) when the intensity
has atoms; in that case we add two occurrences at that location.
For \(\theta=\lambda\delta_c\) and \(G(N)=1\), we obtain
\(\E\, J^2=c^2(\lambda+\lambda^2)\), as required by the Poisson
second moment. This completes the proof.
\end{proof}

\subsection{Universal sufficiency of Theorem \ref{thm:lifting}}\label{sec:stieltjes-suff}

We now prove preservation for an arbitrary initial Stieltjes moment
sequence. We first reduce a Hankel minor to weighted Poisson moment;
we then prove coefficientwise nonnegativity and justify the passage
from compactly supported jumps to general jumps.

\begin{proof}[Proof of Theorem~\ref{thm:lifting}:
\textup{(i)}\(\Rightarrow\)\textup{(iv)}]
Assume \textup{(i)}, and fix an initial Stieltjes moment sequence
\(\boldsymbol a=(a_n)_{n\ge0}\). Choose a finite nonnegative measure
\(\mu\) on \([0,\infty)\) such that
\begin{equation}\label{eq:display-005}
 a_n=\int_0^\infty z^n\,\mu(\dif z)<\infty
 \qquad(n\ge0).
\end{equation}
If \(a_0=0\), then \(\mu\) has zero total mass, so every
\(a_n\) and every \(F_n(\by)\) is zero. The conclusion is
immediate in that case. We may therefore assume \(a_0>0\),
without normalizing its value. This preserves the full range of
initial moment sequences in the statement.

Applying Lemma~\ref{lem:representation} to each shifted cumulant
sequence, we have the drift and jump representation
\eqref{eq:levy-cumulants}. We first assume that every \(\rho_\ell\)
is finite and compactly supported in \((0,\infty)\).

\smallskip\noindent
\emph{Step 1: An integral representation of a fixed minor.}
Fix strictly increasing lists \(I=(i_1,\ldots,i_r)\) and
\(J=(j_1,\ldots,j_r)\) of nonnegative integers, and let
\begin{equation}\label{eq:initial-minor-definition}
 D^\mu_{I,J}(\by):=\det[F_{i_a+j_b}(\by)]_{a,b=1}^r.
\end{equation}
For numerical \(\by\in(0,\infty)^m\), we use \(X_{\by}\) and
\(J_{\by}\) from \eqref{eq:compound-poisson-variable}.
Applying the binomial theorem to \((z+X_{\by})^n\), and using
\eqref{eq:poisson-moments} and \eqref{eq:universal-operator}, we have
\begin{equation}\label{eq:initial-moment-integral}
 \begin{aligned}
 \int_0^\infty \E(z+X_{\by})^n\,\mu(\dif z)
 &=\sum_{k=0}^n\binom nk
       \left(\int_0^\infty z^{n-k}\,\mu(\dif z)\right)
       \E X_{\by}^{\,k}\\
 &=\sum_{k=0}^n\binom nk a_{n-k}M_k(\by)=F_n(\by).
 \end{aligned}
\end{equation}
Every term is finite; we require no convergence assumption on
the formal power series \(A(t)\).

Take independent copies \(J_{\by}^{(1)},\ldots,J_{\by}^{(r)}\), and write
\begin{equation}\label{eq:display-006}
 \boldsymbol W_{\by}(\boldsymbol z)
 =\bJ_{\by}+\boldsymbol z+b(\by)\boldsymbol1.
\end{equation}
The measure whose moments occur in \eqref{eq:initial-moment-integral}
is the convolution of \(\mu\) with the law of \(X_{\by}\) and it has
all moments finite. Applying Andr\'eief's identity
\eqref{eq:andreief-identity} to that measure and the monomials with
exponents \(I,J\), and then using the Schur alternant formula
\eqref{eq:schur-alternant}, we obtain
\begin{equation}\label{eq:initial-andreief}
 D^\mu_{I,J}(\by)
 =\frac1{r!}\int_{[0,\infty)^r}
  \E\!\left[
    \Delta(\bJ_{\by}+\boldsymbol z)^2
    s_\lambda(\boldsymbol W_{\by}(\boldsymbol z))
    s_\eta(\boldsymbol W_{\by}(\boldsymbol z))
  \right]\,\prod_{j=1}^r\mu(\dif z_j),
\end{equation}
where \(\lambda_a=i_{r+1-a}-(r-a)\) and
\(\eta_a=j_{r+1-a}-(r-a)\). To remove the drift from the Vandermonde,
we used the identity
\begin{equation}\label{eq:common-drift-vandermonde}
 \Delta(\bx+c\boldsymbol1)
 =\prod_{a<b}((x_b+c)-(x_a+c))=\Delta(\bx).
\end{equation}

\smallskip\noindent
\emph{Step 2: Schur and multinomial expansions.}
Using the tableau definition \eqref{eq:schur-tableau-definition}, we write
\begin{equation}\label{eq:schur-product-expansion}
 s_\lambda(\boldsymbol w)s_\eta(\boldsymbol w)
 =\sum_{\boldsymbol\alpha\in\Nzero^r}
       c_{\boldsymbol\alpha}\boldsymbol w^{\boldsymbol\alpha},
 \qquad c_{\boldsymbol\alpha}\ge0,
\end{equation}
with finitely many nonzero coefficients. For multi-indices in
\(\Nzero^r\), put \(\boldsymbol\alpha!=\prod_j\alpha_j!\).
Applying the multinomial formula to each factor, we have
\begin{equation}\label{eq:initial-coordinate-multinomial}
 \prod_{j=1}^r(J_{\by}^{(j)}+z_j+b(\by))^{\alpha_j}
 =\sum_{\boldsymbol\beta+\boldsymbol\gamma+
                  \boldsymbol\delta=\boldsymbol\alpha}
  \frac{\boldsymbol\alpha!}
       {\boldsymbol\beta!\boldsymbol\gamma!\boldsymbol\delta!}
  \left(\prod_{j=1}^r(J_{\by}^{(j)})^{\beta_j}\right)
  \boldsymbol z^{\boldsymbol\gamma}b(\by)^{|\boldsymbol\delta|}.
\end{equation}
Substituting these expansions into \eqref{eq:initial-andreief}, and
taking the deterministic factors outside the expectation, we obtain
\begin{equation}\label{eq:detailed-schur-reduction}
 \begin{aligned}
 D^\mu_{I,J}(\by)
 &=\frac1{r!}\sum_{\boldsymbol\alpha}c_{\boldsymbol\alpha}
   \sum_{\boldsymbol\beta+\boldsymbol\gamma+
                  \boldsymbol\delta=\boldsymbol\alpha}
   \frac{\boldsymbol\alpha!}
        {\boldsymbol\beta!\boldsymbol\gamma!\boldsymbol\delta!}
        b(\by)^{|\boldsymbol\delta|}\\
   &\qquad{}\times\int_{[0,\infty)^r}\boldsymbol z^{\boldsymbol\gamma}
   \E\!\left[\Delta(\bJ_{\by}+\boldsymbol z)^2
            \prod_{j=1}^r(J_{\by}^{(j)})^{\beta_j}\right]
       \prod_{j=1}^r\mu(\dif z_j).
 \end{aligned}
\end{equation}
Thus we need to treat finitely many expressions of the form
\begin{equation}\label{eq:schur-reduction-basic-term}
 b(\by)^h\int_{[0,\infty)^r}\boldsymbol z^{\boldsymbol\gamma}
   \E\!\left[\Delta(\bJ_{\by}+\boldsymbol z)^2
            \prod_{j=1}^r(J_{\by}^{(j)})^{\beta_j}\right]
       \prod_{j=1}^r\mu(\dif z_j),
 \qquad h\in\Nzero.
\end{equation}
To separate the positivity argument from the integrability argument,
we first keep \(\boldsymbol z\) fixed and prove coefficientwise
nonnegativity. We then justify integration against
\(\mu\) in each initial coordinate.

\begin{claim}\label{claim:initial-mecke-positive}
For \(\boldsymbol z\in[0,\infty)^r\) and
\(\boldsymbol\beta\in\Nzero^r\), we have
\begin{equation}\label{eq:claim-initial-coefficients}
 P_{\boldsymbol\beta,\boldsymbol z}(\by)
 :=\E\!\left[\Delta(\bJ_{\by}+\boldsymbol z)^2
            \prod_{j=1}^r(J_{\by}^{(j)})^{\beta_j}\right]
 \in\R_{\ge0}[\by].
\end{equation}
\end{claim}
\begin{proof}
\emph{Step 1: Coordinate processes and tagged partition blocks.}
Fix \(\boldsymbol z\) and \(\boldsymbol\beta\), and first take
numerical \(\by\in(0,\infty)^m\). Note
\begin{equation}\label{eq:claim3-coordinate-processes}
 J_{\by}^{(j)}=\int_0^\infty x\,N_{\by}^{(j)}(\dif x),
 \qquad \rho_{\by}=\sum_{\ell=1}^m y_\ell\rho_\ell,
\end{equation}
using independent Poisson random measures \(N_{\by}^{(j)}\) of common
intensity \(\rho_{\by}\). Let
\(\boldsymbol\pi=(\pi_1,\ldots,\pi_r)\), with
\(\pi_j\in\cP_{\beta_j}\). To distinguish blocks in different
coordinates, use the tagged set
\begin{equation}\label{eq:claim3-tagged-blocks}
 \mathcal B_{\boldsymbol\pi}
 :=\{(j,B):1\le j\le r,\ B\in\pi_j\},\qquad
 k_{\boldsymbol\pi}:=\sum_{j=1}^r|\pi_j|.
\end{equation}
For one integration variable \(x_{j,B}\) per tagged block, put
\begin{equation}\label{eq:claim3-block-weights}
 \begin{aligned}
 w_{\boldsymbol\pi}(\boldsymbol x)
   &:=\prod_{(j,B)\in\mathcal B_{\boldsymbol\pi}}x_{j,B}^{|B|},\\
 s_{\boldsymbol\pi,j}(\boldsymbol x)
   &:=\sum_{B\in\pi_j}x_{j,B},\qquad
 \boldsymbol s_{\boldsymbol\pi}(\boldsymbol x)
   :=(s_{\boldsymbol\pi,1}(\boldsymbol x),\ldots,
      s_{\boldsymbol\pi,r}(\boldsymbol x)).
 \end{aligned}
\end{equation}
On inserting the occurrences indexed by \(\pi_j\), we increase the
\(j\)-th jump sum by \(s_{\boldsymbol\pi,j}(\boldsymbol x)\).

\smallskip\noindent
\emph{Step 2: Successive uses of the partitioned Mecke identity.}
For fixed realizations of all other coordinates, the Vandermonde square
is a nonnegative measurable function of \(N_{\by}^{(j)}\).
Applying Lemma~\ref{lem:mecke} to that function and
\((J_{\by}^{(j)})^{\beta_j}\), successively for \(j=1,\ldots,r\),
we have
\begin{equation}\label{eq:claim3-mecke-expansion}
 \begin{aligned}
 P_{\boldsymbol\beta,\boldsymbol z}(\by)
 &=\sum_{\boldsymbol\pi}
   \int_{(0,\infty)^{k_{\boldsymbol\pi}}}
   w_{\boldsymbol\pi}(\boldsymbol x)
   \E\!\left[\Delta(\bJ_{\by}+\boldsymbol z+
                      \boldsymbol s_{\boldsymbol\pi}(\boldsymbol x))^2\right]
   \prod_{(j,B)\in\mathcal B_{\boldsymbol\pi}}
                     \rho_{\by}(\dif x_{j,B})\\
 &=\sum_{\boldsymbol\pi}\int_{(0,\infty)^{k_{\boldsymbol\pi}}}
   w_{\boldsymbol\pi}(\boldsymbol x)
   V_{\boldsymbol z+\boldsymbol s_{\boldsymbol\pi}(\boldsymbol x)}(\by)
   \prod_{(j,B)\in\mathcal B_{\boldsymbol\pi}}
                    \rho_{\by}(\dif x_{j,B}).
 \end{aligned}
\end{equation}
We may interchange the successive integrations by Tonelli's theorem,
because all the integrands in this step are nonnegative.

\smallskip\noindent
\emph{Step 3: An explicit nonnegative coefficient formula.}
For a color assignment
\(c:\mathcal B_{\boldsymbol\pi}\to[m]\), define
\begin{equation}\label{eq:claim3-color-counts}
 n_\ell(c):=|\{(j,B)\in\mathcal B_{\boldsymbol\pi}:c(j,B)=\ell\}|,
 \qquad \boldsymbol n(c):=(n_1(c),\ldots,n_m(c)).
\end{equation}
Expanding the product of measure functions, we obtain
\begin{equation}\label{eq:claim3-colour-expansion}
 \prod_{(j,B)\in\mathcal B_{\boldsymbol\pi}}
                     \rho_{\by}(\dif x_{j,B})
 =\sum_{c:\mathcal B_{\boldsymbol\pi}\to[m]}
    \by^{\boldsymbol n(c)}
    \prod_{(j,B)\in\mathcal B_{\boldsymbol\pi}}
                     \rho_{c(j,B)}(\dif x_{j,B}).
\end{equation}
By Claim~\ref{claim:vandermonde-square} and its proof, we may write
\begin{equation}\label{eq:claim3-shift-coefficients}
 V_{\boldsymbol u}(\by)
 =\sum_{|\boldsymbol\xi|\le D}
             v_{\boldsymbol\xi}(\boldsymbol u)\by^{\boldsymbol\xi},
 \qquad v_{\boldsymbol\xi}(\boldsymbol u)\ge0
 \quad(\boldsymbol u\in[0,\infty)^r).
\end{equation}
Each \(v_{\boldsymbol\xi}(\boldsymbol u)\) is a polynomial in
\(\boldsymbol u\) of degree at most \(2D\). For fixed \(\boldsymbol z\),
the inserted jumps range over a compact set, so these polynomials and
the weights \(w_{\boldsymbol\pi}\) are bounded on the integration
domains. Since every \(\rho_\ell\) is finite, all coefficient
integrals are finite. We may therefore extract coefficients in the
finite polynomial expansions. Combining
\eqref{eq:claim3-mecke-expansion}--\eqref{eq:claim3-shift-coefficients},
we obtain, for \(\boldsymbol\nu\in\Nzero^m\),
\begin{equation}\label{eq:claim3-coefficient-certificate}
 \begin{aligned}
 [\by^{\boldsymbol\nu}]P_{\boldsymbol\beta,\boldsymbol z}(\by)
 &=\sum_{\boldsymbol\pi}
   \sum_{\substack{c:\mathcal B_{\boldsymbol\pi}\to[m]\\
                    \boldsymbol n(c)\le\boldsymbol\nu}}
   \int_{(0,\infty)^{k_{\boldsymbol\pi}}}
        w_{\boldsymbol\pi}(\boldsymbol x)\,
        v_{\boldsymbol\nu-\boldsymbol n(c)}
          (\boldsymbol z+\boldsymbol s_{\boldsymbol\pi}(\boldsymbol x))
   \prod_{(j,B)\in\mathcal B_{\boldsymbol\pi}}
                         \rho_{c(j,B)}(\dif x_{j,B})\ge0.
 \end{aligned}
\end{equation}
Here we compare multi-indices componentwise and set
\(v_{\boldsymbol\xi}=0\) for \(|\boldsymbol\xi|>D\).
Every integrand and every measure in the last expression is
nonnegative, which proves the coefficientwise assertion.
There are only finitely many partitions and colors, and the degree is
at most \(D+|\boldsymbol\beta|\). Thus the expansions are polynomial
identities; their validity on \((0,\infty)^m\) also gives validity
on the boundary. When \(k_{\boldsymbol\pi}=0\), we use the empty-product
convention, with \(\boldsymbol n(c)=0\) and
\(\boldsymbol s_{\boldsymbol\pi}=0\).
\end{proof}

It remains to justify integration against the initial measure and to
remove the restriction on the jump measures. For general measures
\(\rho_\ell\) in \eqref{eq:levy-cumulants}, define, for integers \(L\ge1\),
\begin{equation}\label{eq:display-007}
 \rho_{\ell,L}(\dif x) :=\mathbf1_{[1/L,L]}(x)\rho_\ell(\dif x),\qquad \kappa_{n,L}^{(\ell)} :=b_\ell\mathbf1_{\{n=1\}} +\int_{1/L}^L x^n\,\rho_\ell(\dif x) \qquad(n\ge1).
\end{equation}
Let \(M_{n,L}\) be given by \eqref{eq:bell-partitions} with
\(\kappa_n^{(\ell)}\) replaced by \(\kappa_{n,L}^{(\ell)}\), and put
\begin{equation}\label{eq:claim4-truncated-polynomials}
 F_{n,L}(\by):=\sum_{k=0}^n\binom nk a_{n-k}M_{k,L}(\by),\qquad
 D^\mu_{I,J,L}(\by):=\det[F_{i_a+j_b,L}(\by)]_{a,b=1}^r.
\end{equation}
We now justify the coefficientwise limit as the jump truncation
is removed. In the following claim, we also verify finiteness and
nonnegativity after integration against the initial measure.

\begin{claim}\label{claim:initial-integrability}
For every \(\boldsymbol\alpha\in\Nzero^m\), we have
\begin{equation}\label{eq:claim-initial-limit}
 0\le[\by^{\boldsymbol\alpha}]D^\mu_{I,J}(\by)
 =\lim_{L\to\infty}
       [\by^{\boldsymbol\alpha}]D^\mu_{I,J,L}(\by)<\infty.
\end{equation}
\end{claim}
\begin{proof}
\emph{Step 1: Polynomial bounds and integration against \(\mu\).}
First assume that the jump measures are finite and compactly supported.
For \(\boldsymbol\beta\in\Nzero^r\) and
\(\boldsymbol\nu\in\Nzero^m\), write
\begin{equation}\label{eq:claim4-weighted-coefficients}
 p_{\boldsymbol\beta,\boldsymbol\nu}(\boldsymbol z)
 :=[\by^{\boldsymbol\nu}]P_{\boldsymbol\beta,\boldsymbol z}(\by)
 \ge0\qquad(\boldsymbol z\in[0,\infty)^r).
\end{equation}
In \eqref{eq:claim3-coefficient-certificate}, we integrate translates
of the polynomials \(v_{\boldsymbol\xi}\), whose degrees are at most
\(2D\), over bounded jump domains with finite measures.
Consequently, \(p_{\boldsymbol\beta,\boldsymbol\nu}(\boldsymbol z)\)
is a polynomial in \(\boldsymbol z\) of degree at most \(2D\).
For a constant \(C\) depending on the fixed indices and jump measures,
we therefore have
\begin{equation}\label{eq:claim4-polynomial-bound}
 0\le\boldsymbol z^{\boldsymbol\gamma}
       p_{\boldsymbol\beta,\boldsymbol\nu}(\boldsymbol z)
 \le C(1+z_1+\cdots+z_r)^{2D+|\boldsymbol\gamma|}
 \quad(\boldsymbol z\in[0,\infty)^r).
\end{equation}
To verify integrability of this bound, let \(q\in\Nzero\).
Applying the multinomial formula to \((1+z_1+\cdots+z_r)^q\), we have
\begin{equation}\label{eq:claim4-moment-integrability}
 \begin{aligned}
 \int_{[0,\infty)^r}(1+z_1+\cdots+z_r)^q
                                  \,\prod_{j=1}^r\mu(\dif z_j)&=\sum_{\delta_0+\cdots+\delta_r=q}
       \frac{q!}{\delta_0!\cdots\delta_r!}
       \prod_{j=1}^r\left(\int_0^\infty z_j^{\delta_j}
                                                \,\mu(\dif z_j)\right)\\
 &=\sum_{\delta_0+\cdots+\delta_r=q}
       \frac{q!}{\delta_0!\cdots\delta_r!}\prod_{j=1}^r a_{\delta_j}
       <\infty.
 \end{aligned}
\end{equation}
Using \(q=2D+|\boldsymbol\gamma|\), we obtain
\begin{equation}\label{eq:claim4-integrable-coefficients}
 0\le\int_{[0,\infty)^r}\boldsymbol z^{\boldsymbol\gamma}
       p_{\boldsymbol\beta,\boldsymbol\nu}(\boldsymbol z)
                    \,\prod_{j=1}^r\mu(\dif z_j)<\infty.
\end{equation}
Thus we may integrate the finite polynomial expansion coefficient by
coefficient. In particular, no compact-support assumption on \(\mu\)
is needed.

\smallskip\noindent
\emph{Step 2: Nonnegative coefficients of the minor in the compact case.}
Recall \(b(\by)=\sum_{\ell=1}^m b_\ell y_\ell\), with \(b_\ell\ge0\).
For \(h\in\Nzero\), applying the multinomial formula to \(b(\by)^h\),
we have
\begin{equation}\label{eq:claim4-drift-coefficients}
 b(\by)^h=\sum_{|\boldsymbol\omega|=h}
     \frac{h!}{\boldsymbol\omega!}
     \left(\prod_{\ell=1}^m b_\ell^{\omega_\ell}\right)
       \by^{\boldsymbol\omega}.
\end{equation}
Here and below we interpret a zeroth power as one, including \(0^0\)
in these polynomial coefficient formulas. Combining
\eqref{eq:detailed-schur-reduction},
\eqref{eq:claim4-integrable-coefficients}, and
\eqref{eq:claim4-drift-coefficients}, we obtain the explicit formula
\begin{equation}\label{eq:claim4-minor-coefficients}
 \begin{aligned}
 [\by^{\boldsymbol\alpha}]D^\mu_{I,J}(\by)
 &=\frac1{r!}\sum_{\boldsymbol\tau}c_{\boldsymbol\tau}
   \sum_{\boldsymbol\beta+\boldsymbol\gamma+
                   \boldsymbol\delta=\boldsymbol\tau}
   \frac{\boldsymbol\tau!}
        {\boldsymbol\beta!\boldsymbol\gamma!\boldsymbol\delta!}\times\sum_{\substack{|\boldsymbol\omega|=|\boldsymbol\delta|\\
                             \boldsymbol\omega\le\boldsymbol\alpha}}
   \frac{|\boldsymbol\delta|!}{\boldsymbol\omega!}
   \left(\prod_{\ell=1}^m b_\ell^{\omega_\ell}\right)\\
 &\quad\times\int_{[0,\infty)^r}\boldsymbol z^{\boldsymbol\gamma}
       p_{\boldsymbol\beta,\boldsymbol\alpha-\boldsymbol\omega}
                            (\boldsymbol z)
                    \,\prod_{j=1}^r\mu(\dif z_j)\ge0.
 \end{aligned}
\end{equation}
We used \(\boldsymbol\tau,\boldsymbol\beta,\boldsymbol\gamma,
\boldsymbol\delta\in\Nzero^r\) and
\(\boldsymbol\alpha,\boldsymbol\omega\in\Nzero^m\) to distinguish
the two kinds of multi-indices. Every sum is finite and every factor
is nonnegative. We have therefore proved the desired assertion for
finite compactly supported jump measures.

\smallskip\noindent
\emph{Step 3: Finite compact truncations of general jump measures.}
Return to the general representation \eqref{eq:levy-cumulants}.
For \(\rho_{\ell,L}\) in \eqref{eq:display-007}, we have
\begin{equation}\label{eq:display-008}
 \rho_{\ell,L}((0,\infty)) =\rho_\ell([1/L,L]) \le L\int_0^1 x\,\rho_\ell(\dif x) +\rho_\ell((1,\infty))<\infty.
\end{equation}
We used \(1\le L\,x\) on \([1/L,1]\) and
\eqref{eq:levy-integrability-calculation}. The support lies in
\([1/L,L]\), so the truncated jump measures satisfy the assumptions
of Steps 1 and 2. More explicitly, the nonnegative measure
\begin{equation}\label{eq:claim4-truncated-representation}
 \nu_{\ell,L}(\dif x):=b_\ell\,\delta_0(\dif x)
                              +x\,\rho_{\ell,L}(\dif x)
\end{equation}
has finite mass and compact support, and
\begin{equation}\label{eq:claim4-truncated-shifted-moments}
 \int_0^\infty x^n\,\nu_{\ell,L}(\dif x)
 =b_\ell\mathbf1_{\{n=0\}}
        +\int_{1/L}^L x^{n+1}\,\rho_\ell(\dif x)
 =\kappa_{n+1,L}^{(\ell)}\qquad(n\ge0).
\end{equation}
Applying the compact-case conclusion \eqref{eq:claim4-minor-coefficients}
to these truncated result, we have
\begin{equation}\label{eq:claim4-truncated-positivity}
 [\by^{\boldsymbol\alpha}]D^\mu_{I,J,L}(\by)\ge0
 \qquad(L\ge1).
\end{equation}

\smallskip\noindent
\emph{Step 4: Coefficientwise convergence of moments and minors.}
Since \([1/L,L]\) increases to \((0,\infty)\), applying the monotone
convergence theorem to the nonnegative functions
\(x^n\mathbf1_{[1/L,L]}(x)\), we have
\begin{equation}\label{eq:claim4-cumulant-convergence}
 \kappa_{n,L}^{(\ell)}
 \longrightarrow b_\ell\mathbf1_{\{n=1\}}
                   +\int_0^\infty x^n\,\rho_\ell(\dif x)
 =\kappa_n^{(\ell)}<\infty.
\end{equation}
For each fixed \(n\), the set-partition formula
\eqref{eq:bell-partitions} expresses every coefficient of \(M_{n,L}\)
as a finite polynomial in the cumulants with indices at most \(n\).
Using that formula and \eqref{eq:claim4-truncated-polynomials}, we obtain
\begin{equation}\label{eq:claim4-moment-coefficient-convergence}
 [\by^{\boldsymbol\alpha}]M_{n,L}(\by) \longrightarrow[\by^{\boldsymbol\alpha}]M_n(\by),\qquad [\by^{\boldsymbol\alpha}]F_{n,L}(\by) \longrightarrow[\by^{\boldsymbol\alpha}]F_n(\by).
\end{equation}
To make the last passage to a determinant explicit, expand by permutations
and then extract the coefficient of \(\by^{\boldsymbol\alpha}\):
\begin{equation}\label{eq:claim4-determinant-coefficient-limit}
 \begin{aligned}
 [\by^{\boldsymbol\alpha}]D^\mu_{I,J,L}
 &=\sum_{\sigma\in S_r}\operatorname{sgn}(\sigma)
   \sum_{\boldsymbol\alpha_1+\cdots+
                         \boldsymbol\alpha_r=\boldsymbol\alpha}
      \prod_{a=1}^r
       [\by^{\boldsymbol\alpha_a}]F_{i_a+j_{\sigma(a)},L}\\
 &\longrightarrow
   \sum_{\sigma\in S_r}\operatorname{sgn}(\sigma)
   \sum_{\boldsymbol\alpha_1+\cdots+
                         \boldsymbol\alpha_r=\boldsymbol\alpha}
      \prod_{a=1}^r
       [\by^{\boldsymbol\alpha_a}]F_{i_a+j_{\sigma(a)}}\\
 &=[\by^{\boldsymbol\alpha}]D^\mu_{I,J}.
 \end{aligned}
\end{equation}
All sums in this calculation are finite. Moreover,
\(\deg F_{n,L}\le n\) and
\(\deg D^\mu_{I,J,L}\le\sum_a i_a+\sum_a j_a\), independently of
\(L\). Thus we obtain a finite limiting coefficient, and from
\eqref{eq:claim4-truncated-positivity} we retain its nonnegativity.

For the fixed index lists \(I,J\), put \(N=i_r+j_r\).
Only the initial moments \(a_0,\ldots,a_N\) and the cumulants
\(\kappa_n^{(\ell)}\) with \(1\le n\le N\) occur in this
calculation. Each coefficient of the minor is a polynomial in
these finitely many real numbers. Thus we use monotone convergence
for the nonnegative cumulant integrals, followed by continuity of
finite polynomial expressions for the minor coefficients.
We require no monotonicity of determinant coefficients, no
interchange of an infinite formal power series with an integral,
and no convergence of a moment generating function. This proves
\eqref{eq:claim-initial-limit}.
\end{proof}
Since \(I,J\), the minor order \(r\), and the initial Stieltjes
moment sequence are arbitrary, Claim~\ref{claim:initial-integrability}
yields \textup{(iv)}.
\end{proof}

\subsection{Necessity and completion of the equivalence of Theorem~\ref{thm:lifting}}
\label{sec:stieltjes-necessity}

We next recover the shifted cumulant minors from the lowest nonzero
terms of the moment minors. This will prove necessity using only
numerical Stieltjes moment sequences in \textup{(iii)}.

\begin{proof}[Completion of the proof of Theorem~\ref{thm:lifting}]
We have proved \textup{(i)}\(\Rightarrow\)\textup{(iv)}.
For the initial measure \(\delta_0\), the moments are
\(a_0=1\) and \(a_n=0\) for \(n\ge1\), so
\eqref{eq:universal-operator} gives \(F_n=M_n\). Hence
\textup{(iv)}\(\Rightarrow\)\textup{(ii)} and
\textup{(v)}\(\Rightarrow\)\textup{(iii)}.
Under \textup{(iv)} or \textup{(ii)}, evaluating every nonnegative
polynomial minor at \(\by\in[0,\infty)^m\) gives a numerical
nonnegative minor. Applying Lemma~\ref{lem:moment-criteria}\textup{(ii)}
to the resulting numerical Hankel matrices, we have
\textup{(iv)}\(\Rightarrow\)\textup{(v)} and
\textup{(ii)}\(\Rightarrow\)\textup{(iii)}.
It remains to prove \textup{(iii)}\(\Rightarrow\)\textup{(i)}.

\smallskip\noindent
\emph{Step 1: Restriction to one color and the one-block term.}
Assume that \textup{(iii)} holds, fix \(\ell\in\{1,\ldots,m\}\), and let
\(\boldsymbol e_\ell\) be the \(\ell\)-th coordinate vector. Define
\begin{equation}\label{eq:display-009}
 M_n^{[\ell]}(u):=M_n(u\boldsymbol e_\ell),\qquad
 \sum_{n\ge0}M_n^{[\ell]}(u)\frac{t^n}{n!}=e^{uK_\ell(t)}.
\end{equation}
Substituting \(\by=u\boldsymbol e_\ell\) into the finite partition
formula \eqref{eq:bell-partitions}, we obtain, for \(n\ge1\),
\begin{equation}\label{eq:one-color-partition-expansion}
 \begin{aligned}
 M_n^{[\ell]}(u)
 &=\sum_{\pi\in\cP_n}u^{|\pi|}
                       \prod_{B\in\pi}\kappa_{|B|}^{(\ell)}\\
 &=u\,\kappa_n^{(\ell)}+
   \sum_{\substack{\pi\in\cP_n\\|\pi|\ge2}}
                u^{|\pi|}\prod_{B\in\pi}\kappa_{|B|}^{(\ell)}.
 \end{aligned}
\end{equation}
There is exactly one partition with one block, namely \(\{[n]\}\).
Every other partition contributes a power of \(u\) at least two.
Thus we may write
\begin{equation}\label{eq:linear-cumulant}
 M_n^{[\ell]}(u)=u\,\kappa_n^{(\ell)}+u^2R_{n,\ell}(u),
 \qquad
 R_{n,\ell}(u):=
 \sum_{\substack{\pi\in\cP_n\\|\pi|\ge2}}
        u^{|\pi|-2}\prod_{B\in\pi}\kappa_{|B|}^{(\ell)}\in\R[u].
\end{equation}
For \(n=1\), this is an empty sum and \(R_{1,\ell}=0\).
We assume no sign condition on the coefficients of \(R_{n,\ell}\).

\smallskip\noindent
\emph{Step 2: The lowest-degree term of a shifted minor.}
Take arbitrary strictly increasing nonnegative integer lists
\(I=(i_1,\ldots,i_r)\) and \(J=(j_1,\ldots,j_r)\), and set
\begin{equation}\label{eq:display-010}
 Q_{I,J}^{[\ell]}(u):=
             \det[M_{i_a+j_b+1}^{[\ell]}(u)]_{a,b=1}^r.
\end{equation}
This is a minor of the ordinary moment Hankel matrix with row list
\(I\) and column list \((j_1+1,\ldots,j_r+1)\).
Under \textup{(iii)}, \((M_n^{[\ell]}(u))_{n\ge0}\) is a Stieltjes
moment sequence for every \(u>0\). Applying
Lemma~\ref{lem:moment-criteria}\textup{(ii)} to this sequence, we have
\begin{equation}\label{eq:necessity-moment-minor-positive}
 Q_{I,J}^{[\ell]}(u)\ge0\qquad(u>0).
\end{equation}
Every index \(i_a+j_b+1\) is positive, so we may use
\eqref{eq:linear-cumulant} for every entry.
Factoring \(u\) from each row, and then expanding the determinant
by permutations, we obtain
\begin{equation}\label{eq:necessity-factored-determinant}
 \begin{aligned}
 u^{-r}Q_{I,J}^{[\ell]}(u)
 &=\det\!\left[\kappa_{i_a+j_b+1}^{(\ell)}
                 +uR_{i_a+j_b+1,\ell}(u)\right]_{a,b=1}^r\\
 &=\sum_{\sigma\in S_r}\operatorname{sgn}(\sigma)
     \prod_{a=1}^r\left(\kappa_{i_a+j_{\sigma(a)}+1}^{(\ell)}
                 +u\,R_{i_a+j_{\sigma(a)}+1,\ell}(u)\right).
 \end{aligned}
\end{equation}
In each product, choosing only cumulant factors gives the constant
term. Every other choice contributes at least one factor of \(u\).
In consequence, for a polynomial \(S_{I,J,\ell}(u)\), we have
\begin{equation}\label{eq:leading-shifted-minor}
 Q_{I,J}^{[\ell]}(u)
 =u^r\det[\kappa_{i_a+j_b+1}^{(\ell)}]_{a,b=1}^r
      +u^{r+1}S_{I,J,\ell}(u).
\end{equation}
Dividing by the positive number \(u^r\), and then taking
\(u\to0^+\) in \eqref{eq:necessity-factored-determinant}, we obtain
\begin{equation}\label{eq:display-011}
 \det[\kappa_{i_a+j_b+1}^{(\ell)}]_{a,b=1}^r
 =\lim_{u\to0^+}u^{-r}Q_{I,J}^{[\ell]}(u)\ge0.
\end{equation}
We use nonnegativity of the numerical values for \(u>0\), together
with the fact that all terms of degree below \(r\) vanish; we do not
assume coefficientwise nonnegativity in \textup{(iii)}.

\smallskip\noindent
\emph{Step 3: The Stieltjes criterion for the shifted cumulants.}
Put \(c_n^{(\ell)}:=\kappa_{n+1}^{(\ell)}\) for \(n\ge0\).
For every choice of \(I,J\), we have
\begin{equation}\label{eq:necessity-shifted-hankel}
 \det[c_{i_a+j_b}^{(\ell)}]_{a,b=1}^r
 =\det[\kappa_{i_a+j_b+1}^{(\ell)}]_{a,b=1}^r\ge0.
\end{equation}
Thus every minor of the infinite Hankel matrix
\([c_{i+j}^{(\ell)}]_{i,j\ge0}\) is nonnegative.
Applying the total-positivity characterization in
Lemma~\ref{lem:moment-criteria} \textup{(ii)} to
\((c_n^{(\ell)})_{n\ge0}\), we have a nonnegative measure
\(\nu_\ell\) on \([0,\infty)\) such that
\begin{equation}\label{eq:necessity-recovered-measure}
 \kappa_{n+1}^{(\ell)}
 =\int_0^\infty x^n\,\nu_\ell(\dif x)<\infty
 \qquad(n\ge0).
\end{equation}
Since \(\ell\) was arbitrary, we obtain \textup{(i)}.
Combining the implications in Step 1 with
\textup{(iii)}\(\Rightarrow\)\textup{(i)}, we prove the equivalence
of all five conditions. We used only finite polynomial identities
and limits of finite determinants. This completes the proof.
\end{proof}

\section{Proof of the Hamburger lifting theorem}
\label{sec:hamburger-proof}

In this section, we will prove Theorem~\ref{thm:hamburger} stated in
Section~\ref{sec:hamburger}. The proof uses the Hamburger moment
criterion in Lemma~\ref{lem:moment-criteria} and the Poisson
Vandermonde expansion of Lemma~\ref{lem:vandermonde}.

\begin{proof}[Proof of Theorem~\ref{thm:hamburger}]
We first derive the necessary moment condition from
\textup{(iv)}. For sufficiency, we treat compactly supported jumps,
compute the highest-degree component, and then approximate arbitrary
twice-shifted moment measures. Finally, we use Gaussian initial
measures to deduce \textup{(iv)} from \textup{(ii)}, thereby
completing the equivalence.

\Needspace{12\baselineskip}
To recover the twice-shifted cumulants, we test the moment
function on squares of polynomials divisible by the coordinate
variable and let a single intensity tend to zero.
In the following claim, we make this limiting argument explicit
and recover the required nonnegative quadratic forms.

\begin{claim}\label{claim:hamburger-necessity}
Under condition~\textup{(iv)}, for each \(\ell\in\{1,\ldots,m\}\), we have
\begin{equation}\label{eq:claim-hamburger-quadratic}
 \sum_{i,j=0}^d c_ic_j\kappa_{i+j+2}^{(\ell)}\ge0
 \qquad\bigl(d\ge0,\ (c_0,\ldots,c_d)\in\R^{d+1}\bigr).
\end{equation}
In consequence, \((\kappa_{n+2}^{(\ell)})_{n\ge0}\) is a Hamburger
moment sequence.
\end{claim}
\begin{proof}
For the implication \textup{(iv)}\(\Rightarrow\)\textup{(i)},
fix \(\ell\) and any polynomial \(p(x)=\sum_{i=0}^{d}c_ix^i\).
For every \(u>0\), applying the Hamburger moment property to
the square of the product of the coordinate variable and this
polynomial, we obtain
\begin{equation}\label{eq:display-018}
 \frac1u\sum_{i,j=0}^{d}c_ic_jM_{i+j+2}(u\boldsymbol e_\ell)
 \ge0,
\end{equation}
since the numerator is the integral of \(x^2p(x)^2\) against a nonnegative measure.
For every \(n\ge1\),
\(M_n(u\boldsymbol e_\ell)=u\,\kappa_n^{(\ell)}+O(u^2)\).
Only finitely many moments occur in the quadratic form.
Thus, letting \(u\downarrow0\) and passing to the limit term
by term, we obtain
\begin{equation}\label{eq:display-019}
 \sum_{i,j=0}^{d}c_ic_j\kappa_{i+j+2}^{(\ell)}\ge0.
\end{equation}
Applying the Hamburger moment theorem to these nonnegative
quadratic forms, we obtain the representing measure required in
\textup{(i)}. Since the quadratic forms involve only cumulants
of order at least two, we impose no condition on the first cumulant.
\end{proof}
\Needspace{8\baselineskip}
We now turn to sufficiency. We begin with finite compactly
supported jump measures, for which we can construct the output
measure directly and then obtain coefficientwise positivity from
an expansion in squares. We carry out this construction of a representing measure and
the subsequent square-integral expansion in the following claim.

\begin{claim}\label{claim:hamburger-compact}
For finite compactly supported jump measures on \(\R\setminus\{0\}\),
\begin{equation}\label{eq:claim-hamburger-compact}
 H_r^a(\by)\in\R_{\ge0}[\by],\qquad
 \deg_{\by}H_r^a\le D,\qquad D=\binom r2.
\end{equation}
Moreover, for each \(\by\in[0,\infty)^m\), there is a nonnegative
measure \(\mu_{\by}\) on \(\R\) such that
\begin{equation}\label{eq:claim-hamburger-measure}
 F_n(\by)=\int_{-\infty}^\infty x^n\,\mu_{\by}(\dif x)\qquad(n\ge0).
\end{equation}
\end{claim}
\begin{proof}
To prove \textup{(i)}\(\Rightarrow\)\textup{(ii)} and
\textup{(iii)} in this case, take finite compactly supported jump
measures \(\rho_1,\ldots,\rho_m\) on \(\R\setminus\{0\}\).
Let \(b_1,\ldots,b_m\) be arbitrary real drifts, and set
\begin{equation}\label{eq:display-020}
 K_\ell(t)=b_\ell t+
           \int_{-\infty}^\infty (e^{tx}-1)\,\rho_\ell(\dif x).
\end{equation}
For a numerical \(\by\ge\boldsymbol0\), let \(J_{\by}\) be the
compound-Poisson sum with intensity
\(\rho_{\by}=\sum_\ell y_\ell\rho_\ell\), and put
\(b(\by)=\sum_\ell y_\ell b_\ell\).
If \(\mu\) represents the initial Hamburger sequence, then
\begin{equation}\label{eq:display-021}
 F_n(\by)=\int_{-\infty}^\infty \E(z+J_{\by}+b(\by))^n\,\mu(\dif z).
\end{equation}
By integrating the law of the translated compound-Poisson
variable against the initial measure, we obtain a finite
nonnegative measure on \(\R\) with these moments. This proves
\textup{(iii)} for compactly supported jumps.

To prove coefficientwise positivity, we next express each
leading determinant as an integral of a squared Vandermonde.
Let \(\bJ_{\by}\) consist of \(r\) independent copies of
\(J_{\by}\).  By Andr\'eief's identity \eqref{eq:andreief-identity} and the
common-translation invariance of the Vandermonde, we obtain
\begin{equation}\label{eq:real-andreief}
 H_r^a(\by)
 =\frac1{r!}\int_{\R^r}
       \E\left[\left(\Delta(\bs+\bJ_{\by})\right)^2\right]
       \,\prod_{j=1}^r\mu(\dif s_j).
\end{equation}
We have thus removed the common drift before expanding in the intensities.
Applying the real-line version of Lemma~\ref{lem:vandermonde} to
\(\Delta(\bs+\bJ_{\by})\), with
\(E=\{1,\ldots,r\}\times(\R\setminus\{0\})\) and
\(\Lambda_\ell(\{b\}\times dx)=\rho_\ell(dx)\), and then
integrating each initial coordinate against \(\mu\), we obtain
\begin{equation}\label{eq:real-fock}
 H_r^a(\by)=\frac1{r!}\sum_{k=0}^{D}\frac1{k!}
 \int_{\R^r}\int_{E^k}
 \bigl(\Dshift_{z_1}\cdots\Dshift_{z_k}
                         \Delta(\bs)\bigr)^2
 \prod_{h=1}^{k}\Lambda_{\by}(\dif z_h)
 \,\prod_{j=1}^r\mu(\dif s_j).
\end{equation}
Expanding the product of measure functions, we express each
coefficient as an integral of a square against nonnegative measures.
The difference kernels are independent of the indeterminates and
polynomial in the initial positions and bounded jump sizes.
By the absolute moment hypothesis on \(\mu\), every coefficient
integral is finite. We therefore obtain \textup{(ii)}, with
degree at most \(D\), as required.
\end{proof}
\Needspace{6\baselineskip}
Having established the degree bound in the compact-jump case,
we now compute the highest-degree component. At this degree,
the differences exhaust the degree of the Vandermonde, so the
result depends only on the second jump moments.
We now evaluate this component by counting the nonzero derivatives
of the Vandermonde, as described in the following claim.

\begin{claim}\label{claim:hamburger-top}
For finite compactly supported jump measures, we have
\begin{equation}\label{eq:claim-hamburger-highest}
 [H_r^a]_D(\by)
 =a_0^r\left(\prod_{j=0}^{r-1}j!\right)
      \left(\sum_{\ell=1}^m y_\ell\kappa_2^{(\ell)}\right)^D.
\end{equation}
\end{claim}
\begin{proof}
In the square-integral expansion, only the term \(k=D\)
contributes to the degree-\(D\) part. Write
\(\partial_j=\partial/\partial s_j\).
Since \(\Delta(\bs)\) is homogeneous of degree \(D\), each
successive difference lowers the degree by at least one.
Consequently, only the linear term in each jump can contribute
after all differences have been taken, and we obtain
\begin{equation}\label{eq:display-022}
 \Dshift_{(b_1,x_1)}\cdots\Dshift_{(b_D,x_D)}\Delta
 =\left(\prod_{h=1}^{D}x_h\right)
                     \partial_{b_1}\cdots\partial_{b_D}\Delta.
\end{equation}
On the right side, we obtain a constant derivative. To sum its
square over the coordinate choices, note that in the monomial
expansion of \(\Delta\), the exponent vectors are precisely the
\(r!\) permutations of \((0,1,\ldots,r-1)\).
A nonzero derivative occurs only when the coordinate
multiplicities agree with one of these vectors. For each such
vector, the absolute value of the derivative is
\(\prod_{j=0}^{r-1}j!\), and the number of ordered derivative
lists producing it is \(D!/\prod_{j=0}^{r-1}j!\).
In consequence,
\begin{equation}\label{eq:display-023}
 \sum_{b_1,\ldots,b_D=1}^{r}
   (\partial_{b_1}\cdots\partial_{b_D}\Delta)^2
 =r!D!\prod_{j=0}^{r-1}j!.
\end{equation}
For each jump-square factor, integration gives
\(\int_{-\infty}^\infty x^2\rho_{\by}(\dif x)=v(\by)\).
Since the remaining kernel is constant in the initial positions,
integration against the product initial measure gives \(a_0^r\).
Substituting these evaluations into \eqref{eq:real-fock}, we obtain
\eqref{eq:hamburger-top-degree} in the compact-jump case.
\end{proof}

\Needspace{12\baselineskip}
It remains to remove the compact-jump assumption and to include
a possible mass at zero in the twice-shifted moment measure.
We approximate that mass by symmetric small jumps and adjust
the drift to preserve the first cumulant.
To justify this approximation, we verify convergence of the
cumulants and of every fixed determinant coefficient in the
following claim.

\begin{claim}\label{claim:hamburger-limit}
Under condition~\textup{(i)}, there are drift-plus-compact-jump
exponents \(K_{\ell,N}(t)\) such that
\begin{equation}\label{eq:claim-hamburger-limit}
 \begin{aligned}
 K_{\ell,N}^{(n)}(0)&\longrightarrow\kappa_n^{(\ell)}
                         &&(n\ge1,\ 1\le\ell\le m),\\
 [\by^{\boldsymbol\alpha}]H_r^{a,N}(\by)
     &\longrightarrow[\by^{\boldsymbol\alpha}]H_r^a(\by)
                         &&(r\ge1,\ \boldsymbol\alpha\in\Nzero^m).
 \end{aligned}
\end{equation}
Here \(N\to\infty\), and \(H_r^{a,N}(\by)\) denotes the determinant
formed from \(K_{1,N}(t),\ldots,K_{m,N}(t)\). In particular, we retain
\eqref{eq:claim-hamburger-compact} and
\eqref{eq:claim-hamburger-highest} in the limit.
\end{claim}
\begin{proof}
To treat arbitrary data in \textup{(i)}, we choose nonnegative measures
\(\nu_\ell\) on \(\R\), with all absolute moments finite, such that
\begin{equation}\label{eq:display-024}
 \kappa_{n+2}^{(\ell)}=\int_{-\infty}^\infty x^n\,\nu_\ell(\dif x)
 \qquad(n\ge0),
 \qquad g_\ell=\nu_\ell(\{0\}).
\end{equation}
To obtain finite jump measures, we truncate away from zero
and infinity before dividing by the square of the jump size.
We represent the atom at zero by two symmetric jumps.
For integers \(N\ge1\), define
\begin{equation}\label{eq:real-jump-approximation}
 \begin{aligned}
 \rho_{\ell,N}(\dif x)
 &=\mathbf1_{\{1/N\le|x|\le N\}}x^{-2}\nu_\ell(\dif x)
   +\frac{g_\ell N^2}{2}(\delta_{1/N}+\delta_{-1/N})(\dif x),\\
 b_{\ell,N}&=\kappa_1^{(\ell)}
       -\int_{-\infty}^\infty x\,\rho_{\ell,N}(\dif x).
 \end{aligned}
\end{equation}

The jump measure in \eqref{eq:real-jump-approximation} is finite:
its total mass is at most
\(N^2\nu_\ell(\R)+g_\ell N^2\), since
\(x^{-2}\le N^2\) on the truncated domain. It also has bounded
support. Hence its first absolute moment is finite and
\(b_{\ell,N}\) is a well-defined real number.

Let \(K_{\ell,N}\) be the associated drift-plus-jump exponent.
Its first cumulant is exactly \(\kappa_1^{(\ell)}\), while for
\(n\ge2\) its \(n\)-th cumulant is
\begin{equation}\label{eq:display-025}
 \int_{-N}^{-1/N}x^{n-2}\,\nu_\ell(\dif x)
 +\int_{1/N}^N x^{n-2}\,\nu_\ell(\dif x)
 +\frac{g_\ell}{2}N^{2-n}(1+(-1)^n).
\end{equation}
For the second cumulant, the final term reproduces the atom at
zero; for every higher cumulant, it tends to zero. The two
truncated integrals converge by finiteness of the corresponding
absolute moments. We therefore obtain convergence to
\(\kappa_n^{(\ell)}\), including the possible Gaussian contribution.

Since each coefficient of a moment polynomial depends on only
finitely many cumulants, we obtain coefficientwise convergence of
every fixed moment polynomial \(F_{n,N}(\by)\) to
\(F_n(\by)\). We then pass to the limit in each determinant,
which is a finite polynomial expression in these moments, with
uniformly bounded degree. Thus we retain coefficientwise
positivity, the degree bound, and the formula for the degree-\(D\)
component.

For each numerical \(\by\ge\boldsymbol0\), we also have positive
semidefinite finite Hankel matrices for \((F_{n,N}(\by))\).
Taking entrywise limits, we retain positive semidefiniteness.
Applying the Hamburger moment theorem to the limiting sequence,
we obtain \textup{(iii)} without an analytic assumption.
Finally, since every approximating determinant is independent
of the chosen first cumulants, we retain common-drift invariance
in the limit. This proves all the asserted consequences of
\textup{(i)}.
\end{proof}

We can now complete the equivalence. Choosing the initial
measure \(\delta_0\) in condition~\textup{(iii)}, we obtain
\textup{(iv)}. It remains to deduce \textup{(iv)} directly from
\textup{(ii)}. Assume \textup{(ii)} and fix
\(\by\ge\boldsymbol0\).
To use Sylvester's criterion, we first make all leading
determinants strictly positive by choosing a Gaussian initial
measure. For \(\varepsilon>0\), take a centered Gaussian
probability measure of variance \(\varepsilon\), and denote its
moment sequence by \(\boldsymbol a^\varepsilon=(a_n^\varepsilon)_{n\ge0}\).
For this initial measure, every leading Hankel determinant is
strictly positive, since the integral of the square of a nonzero
polynomial is positive. By \textup{(ii)}, the output determinants
have nonnegative coefficients and hence are bounded below by
their constant terms. Thus, for every \(r\), we obtain
\begin{equation}\label{eq:display-026}
 H_r^{a^{\varepsilon}}(\by)
 \ge H_r^{a^{\varepsilon}}(\boldsymbol0)>0.
\end{equation}

Fix \(d\ge0\) and form the symmetric matrix
\(\mathsf H_d^\varepsilon(\by)
=[F_{i+j}^{a^\varepsilon}(\by)]_{i,j=0}^{d}\).
For each \(1\le r\le d+1\), its leading principal minor is
strictly positive by \eqref{eq:display-026}.
Applying Sylvester's criterion~\cite[Section~7.2]{HornJohnson2013}
to \(\mathsf H_d^\varepsilon(\by)\), we obtain positive
definiteness.

As \(\varepsilon\downarrow0\), the Gaussian initial moments
converge to those of \(\delta_0\). Using the finite convolution
formula \eqref{eq:universal-operator}, we therefore obtain
\(F_n^{a^\varepsilon}(\by)\longrightarrow M_n(\by)\) for every
fixed \(n\). In particular, for every real vector
\(\boldsymbol c=(c_0,\ldots,c_d)\), we have
\(\sum_{i,j=0}^d c_ic_jM_{i+j}(\by)
=\lim_{\varepsilon\downarrow0}
  \boldsymbol c^{\mathsf T}\mathsf H_d^\varepsilon(\by)
  \boldsymbol c\ge0\).
Since \(d\) and \(\boldsymbol c\) are arbitrary, applying the
Hamburger moment criterion in
Lemma~\ref{lem:moment-criteria}\textup{(i)} to
\((M_n(\by))_{n\ge0}\), we obtain \textup{(iv)}.
This completes the proof.

\end{proof}

\section{Composition, convolution, and substitution}\label{sec:universal}

We now study the operators associated with
Theorem~\ref{thm:lifting}. We first prove a binomial-convolution
identity for two independently lifted families and recall a
classical logarithm criterion. We then establish composition and
substitution closure for the admissible exponents.

Define the class of admissible formal exponents by
\begin{equation}\label{eq:admissible-class}
 \mathcal C=
 \left\{K(t)=\sum_{n\ge1}\kappa_n\frac{t^n}{n!}\in t\R[[t]]:
       (\kappa_{n+1})_{n\ge0}\text{ is a Stieltjes moment sequence}
 \right\}.
\end{equation}

Recall the operator \(\mathcal T_{\by}\) defined in
\eqref{eq:universal-operator}. The finite sum in \eqref{eq:universal-operator} also defines these
operators on polynomial-valued sequences.
They satisfy the formal identities
\begin{equation}\label{eq:operator-semigroup}
 \mathcal T_{\by}\mathcal T_{\boldsymbol u}
 =\mathcal T_{\by+\boldsymbol u},
 \qquad
 \mathcal T_{\boldsymbol0}=\operatorname{Id}.
\end{equation}
Indeed, their exponential generating functions are multiplied by
\(\exp(\sum_\ell y_\ell K_\ell)\) and
\(\exp(\sum_\ell u_\ell K_\ell)\), respectively.
All coefficient sums are finite.
By Theorem~\ref{thm:lifting}, after substituting
\(y_\ell\mapsto y_\ell+u_\ell\), we obtain joint coefficientwise positivity
for the output in both parameter vectors.
More generally, substituting arbitrary polynomials
\(y_\ell=p_\ell(\boldsymbol w)\in\R_{\ge0}[\boldsymbol w]\)
preserves the coefficientwise conclusion.

\subsection{Binomial convolution of lifted families}

Wang and Zhu~\cite[Corollary~4.3]{Wang-Zhu} recalled the classical
fact that binomial convolution preserves real Stieltjes moment
sequences. By joint lifting for the exponential families specified, we will extend it to the coefficientwise Hankel-total positivity in the following
Proposition~\ref{prop:lifted-binomial-convolution}. In fact, at fixed nonnegative indeterminates, we can apply it to \(F(\by)\) and \(G(\boldsymbol z)\) to obtain a
numerical Stieltjes sequence.

\begin{proposition}[Binomial convolution of two lifted families]
\label{prop:lifted-binomial-convolution}
Let \(A,B\) be initial Stieltjes moment series, and let
\(K_1,\ldots,K_m,L_1,\ldots,L_p\in\mathcal C\). Define
\(F_n(\by)\) and \(G_n(\boldsymbol z)\) by
\begin{equation}\label{eq:two-lifted-families}
 \sum_{n\ge0}F_n(\by)\frac{t^n}{n!}
 =A(t)e^{\sum_i y_iK_i(t)},\qquad
 \sum_{n\ge0}G_n(\boldsymbol z)\frac{t^n}{n!}
 =B(t)e^{\sum_j z_jL_j(t)}.
\end{equation}
Then the binomial convolution
\begin{equation}\label{eq:lifted-binomial-convolution}
 H_n(\by,\boldsymbol z)
 =\sum_{k=0}^n\binom nk F_k(\by)G_{n-k}(\boldsymbol z)
 \qquad(n\ge0)
\end{equation}
is jointly coefficientwise Hankel-totally positive in
\(\by,\boldsymbol z\).
Write \(\boldsymbol a\star\boldsymbol b\) for binomial convolution
and \(\mathcal T^K\) for the lifting operator with exponent list
\(K\). We have
\begin{equation}\label{eq:convolution-operator-intertwining}
 (\mathcal T^K_{\by}\boldsymbol a)\star
 (\mathcal T^L_{\boldsymbol z}\boldsymbol b)
 =\mathcal T^{(K,L)}_{(\by,\boldsymbol z)}
       (\boldsymbol a\star\boldsymbol b).
\end{equation}
For the same exponent list, the right side is
\(\mathcal T^K_{\by+\boldsymbol z}
(\boldsymbol a\star\boldsymbol b)\).
The coefficientwise conclusion is preserved by identifying
variables.
\end{proposition}

\begin{proof}
Let \(\mu,\nu\) represent the initial sequences.
Integrating \((x+s)^n\) against \(\mu(\dif x)\nu(\dif s)\),
we obtain the moments of their binomial convolution. Thus
\(A(t)B(t)\) is an initial Stieltjes moment series, including
arbitrary total masses and the zero-measure cases. Multiplying
the two generating functions, we obtain
\(A(t)B(t)\exp(\sum_i y_iK_i+\sum_j z_jL_j)\).
Applying Theorem~\ref{thm:lifting} to this initial series and
the concatenated exponent list, we derive the joint coefficientwise
assertion. Extracting \(n![t^n]\), we obtain
\eqref{eq:convolution-operator-intertwining}. For a common list,
we combine its intensities by addition. Finally, substituting
the same indeterminate for selected intensity variables, we retain
nonnegative coefficients.
\end{proof}

\subsection{Logarithms of ordinary Stieltjes series}

We use the following classical consequence of the infinite
divisibility of mixtures of exponential distributions
\cite{Goldie1967}. If \(X\) has Stieltjes moments \(b_n\), and
\(E\) is an independent unit exponential variable, then
\(\E(XE)^n=n!b_n\). Applying Goldie's theorem to \(XE\) and
the L\'evy representation to its nonnegative infinitely divisible
law, we obtain the shifted cumulant moment property below.
We include an analytic proof using
\cite[Theorems~3.2, 7.3 and Corollary~3.8]{SSV2012}, so that
we also cover ordinary formal series with zero convergence radius.
\begin{proposition}[A logarithm criterion]
\label{powerapp:prop:ordinary-log}
Let \((b_n)_{n\ge0}\) be a real Stieltjes moment sequence with
\(b_0=1\). Define the ordinary formal power series \(B(t)\) and
the exponential coefficients \(\kappa_n\) of its logarithm by
\begin{equation}\label{powerapp:eq:ordinary-log}
 B(t)=\sum_{n\ge0}b_nt^n,\qquad
 \log B(t)=\sum_{n\ge1}\kappa_n\frac{t^n}{n!}.
\end{equation}
Then \((\kappa_{n+1})_{n\ge0}\) is a Stieltjes moment sequence;
equivalently, \(\log B\in\mathcal C\).
\end{proposition}

\begin{proof}
Choose a probability measure \(\mu\) on \([0,\infty)\) with moments
\(b_n\), and, for \(s\ge0\), let
\begin{equation}\label{powerapp:eq:stieltjes-laplace}
 L(s)=\int_0^\infty\frac{\mu(\dif x)}{1+sx},\qquad
 h(s)=\frac1{L(s)}-1,\qquad \phi(s)=-\log L(s).
\end{equation}
We have \(0<L(s)\le1\) and \(L(0)=1\). For \(x>0\), we can
write \((1+sx)^{-1}=x^{-1}(s+x^{-1})^{-1}\) and a mass of \(\mu\)
at zero contributes a nonnegative constant. Thus \(L(s)\) is a
Stieltjes function in the classical terminology of
\cite[Chapter~2]{SSV2012}. By the reciprocal characterization in
\cite[Theorem~7.3]{SSV2012}, its reciprocal is a Bernstein function.
Since \(L(s)^{-1}\ge1\), we may subtract its value at zero
and conclude that \(h(s)\) is also a Bernstein function.
We also know that \(\log(1+s)\) is a Bernstein function. Applying the
composition rule in \cite[Corollary~3.8]{SSV2012} to \(h(s)\) and
\(\log(1+s)\), we obtain that \(\phi(s)=\log(1+h(s))\) is a
Bernstein function with \(\phi(0)=0\).

To identify all derivatives at zero, we differentiate the first
integral in \eqref{powerapp:eq:stieltjes-laplace}. For every
\(n\ge0\), using the integrable bound \(n!x^n\) and applying
dominated convergence to the differentiated integrals, we have
\begin{equation}\label{powerapp:eq:stieltjes-jets}
 L^{(n)}(s)=(-1)^n n!\int_0^\infty
                   \frac{x^n}{(1+sx)^{n+1}}\,\mu(\dif x),\qquad
 L^{(n)}(0+)=(-1)^n n!b_n.
\end{equation}
In consequence, by the finite coefficient formulas for a logarithm,
\(\phi^{(n+1)}(0+)=(-1)^n\kappa_{n+1}\).
At each order, we use only finitely many derivatives in this
comparison, so no convergence assumption on \(B(t)\) is needed.

By the classical L\'evy representation
\cite[Theorem~3.2]{SSV2012}, there are a drift \(d\ge0\) and a
positive measure \(\rho\) on \((0,\infty)\) such that
\(\phi(s)=ds+\int_0^\infty(1-e^{-sx})\,\rho(\dif x)\), with
\(\int_0^\infty\min\{1,x\}\,\rho(\dif x)<\infty\).
For \(s>0\), we may differentiate under this integral. Then, by
monotone convergence as \(s\downarrow0\), we obtain
\begin{equation}\label{powerapp:eq:ordinary-log-moments}
 \kappa_{n+1}=(-1)^n\phi^{(n+1)}(0+)
 =d\,\mathbf1_{\{n=0\}}+
             \int_0^\infty x^{n+1}\,\rho(\dif x)\qquad(n\ge0).
\end{equation}
By \eqref{powerapp:eq:stieltjes-jets}, all derivatives at zero
are finite, and hence so are the moment integrals above. We may
therefore use the finite positive measure
\(\nu(\dif x)=d\,\delta_0(\dif x)+x\rho(\dif x)\) to represent
\((\kappa_{n+1})_{n\ge0}\). This completes the proof..
\end{proof}

\subsection{Nonnegative mixtures and composition}

We begin with closure under nonnegative mixtures. The proof
uses quadratic forms and requires no measurable choice of
representing measures.

\begin{proposition}[Nonnegative mixtures]\label{prop:mixtures}
Let \(f_n(s)\) be measurable functions on \([0,\infty)\), and suppose
that \((f_n(s))_{n\ge0}\) is a Stieltjes moment sequence for every \(s\).
Let \(\sigma\) be a nonnegative measure on \([0,\infty)\) such that
\(\int_0^\infty  f_n(s)\,\sigma(\dif s)<\infty\) for every \(n\ge0\).
Then
\begin{equation}\label{eq:display-012}
 g_n=\int_0^\infty  f_n(s)\,\sigma(\dif s)
\end{equation}
is a Stieltjes moment sequence.
\end{proposition}

\begin{proof}
For \(h=0,1\) and any finite real vector \((c_0,\ldots,c_d)\), we have
\begin{equation}\label{eq:display-013}
 \sum_{i,j=0}^d c_ic_jg_{i+j+h}
 =\int_0^\infty \sum_{i,j=0}^d c_ic_jf_{i+j+h}(s)\,d\sigma(s)\ge0.
\end{equation}
The integrand is nonnegative by the Stieltjes property at each \(s\).
Since every term is integrable by hypothesis, we may integrate the
finite sums term by term. Applying Lemma~\ref{lem:moment-criteria}
to the resulting quadratic forms, we obtain the Stieltjes property.
\end{proof}

For \(K(t),L(t)\in t\R[[t]]\), with
\(K(t)=\sum_{k\ge1}\kappa_k t^k/k!\), we use \(\circ\) for
composition by formal substitution:
\begin{equation}\label{eq:formal-composition-definition}
 (K\circ L)(t):=K(L(t))
       =\sum_{k\ge1}\frac{\kappa_k}{k!}L(t)^k.
\end{equation}
Because \(L(0)=0\), only finitely many powers contribute to each
coefficient of \(t\); thus we have a well-defined formal power series.
Using the mixture property together with this definition, we establish
closure under composition.

\begin{theorem}[Composition closure]\label{thm:composition}
If \(K,L\in\mathcal C\), then \(K\circ L\in\mathcal C\).
\end{theorem}

\begin{proof}
Fix a real number \(u\ge0\), and write
\begin{equation}\label{eq:display-014}
 e^{uK(w)}=\sum_{k\ge0}B_k(u)\frac{w^k}{k!}.
\end{equation}
By Theorem~\ref{thm:lifting}, there is a nonnegative measure
\(\mu_u\) on \([0,\infty)\), with all moments finite, such that
\(B_k(u)=\int_0^\infty  s^k\,\mu_u(\dif s)\).
Also define the moment polynomials associated with \(L\) by
\begin{equation}\label{eq:display-015}
 e^{sL(t)}=\sum_{n\ge0}P_n(s)\frac{t^n}{n!}.
\end{equation}
For every \(s\ge0\), \((P_n(s))_{n\ge0}\) is a Stieltjes moment
sequence, again by Theorem~\ref{thm:lifting}.
Moreover, \(P_n(s)\) is a polynomial of degree at most \(n\).
Since \(L(0)=0\), only powers \(L(t)^k\) with \(k\le n\)
contribute to a fixed coefficient, giving
\begin{equation}\label{eq:composition-mixture}
 n![t^n]e^{uK(L(t))}
 =\int_0^\infty P_n(s)\,\mu_u(\dif s).
\end{equation}
Since \(P_n(s)\) is a polynomial and \(\mu_u\) has all moments
finite, the integral is finite.
By Proposition~\ref{prop:mixtures}, we conclude that the sequence on the
left-hand side is a Stieltjes moment sequence.
Since this holds for every \(u\ge0\), applying the implication \textup{(iii)}\(\Rightarrow\)\textup{(i)}
of Theorem~\ref{thm:lifting} to the formal exponent \(K(L(t))\),
we obtain \(K\circ L\in\mathcal C\).
We use only finite coefficient calculations and polynomial
integrals. Hence we can use this argument even when the formal power
series have zero radius of convergence.
\end{proof}

\subsection{Universal substitution--convolution preservers}

The next theorem characterizes formal substitutions of the
initial series that preserve universal coefficientwise Hankel-total positivity.

\begin{theorem}[Substitution and convolution]\label{thm:substitution}
Fix \(L\in t\R[[t]]\) and \(K_1,\ldots,K_m\in t\R[[t]]\).
For an initial real sequence \(\boldsymbol a=(a_n)_{n\ge0}\) with
formal exponential generating function
\(A(t)=\sum_{n\ge0}a_nt^n/n!\), define \(Q_n(\by)\) by
\begin{equation}\label{eq:substitution-operator}
 \sum_{n\ge0}Q_n(\by)\frac{t^n}{n!}
 =A(L(t))\exp\!\left(\sum_{\ell=1}^{m}y_\ell K_\ell(t)\right).
\end{equation}
The following conditions are equivalent:
\begin{enumerate}
\renewcommand{\labelenumi}{\textup{(\roman{enumi})}}
\item \(L,K_1,\ldots,K_m\in\mathcal C\).
\item For every initial real Stieltjes moment sequence,
\((Q_n(\by))_{n\ge0}\) is jointly coefficientwise Hankel-totally
positive in \(\by\).
\item For every initial real Stieltjes moment sequence and every
\(\by\in[0,\infty)^m\), \((Q_n(\by))_{n\ge0}\) is a Stieltjes
moment sequence.
\end{enumerate}
\end{theorem}

\begin{proof}
Suppose \textup{(i)} holds, and let \(\mu\) represent the initial
sequence.  Define \(P_n(s)\) by
\(e^{sL(t)}=\sum_{n\ge0}P_n(s)t^n/n!\).
For each \(s\ge0\), we know from Theorem~\ref{thm:lifting} that
\((P_n(s))_{n\ge0}\) is a Stieltjes moment sequence.
Since \(P_n(s)\) is a polynomial of degree at most \(n\), we obtain
by formal substitution the finite coefficient identity
\begin{equation}\label{eq:display-016}
 b_n:=n![t^n]A(L(t))
     =\int_0^\infty P_n(s)\,\mu(\dif s).
\end{equation}
The integral is finite. The numerical preservation step is the
row-moment criterion of Wang and Zhu~\cite[Theorem~4.2]{Wang-Zhu},
or equivalently Proposition~\ref{prop:mixtures}. Thus \((b_n)\) is a
Stieltjes moment sequence.
Applying Theorem~\ref{thm:lifting} to the initial Stieltjes sequence
\((b_n)_{n\ge0}\) and the exponents \(K_1(t),\ldots,K_m(t)\),
we obtain \textup{(ii)}.  By nonnegative specialization, we have
\textup{(ii)}\(\Rightarrow\)\textup{(iii)}.

Conversely, assume that \textup{(iii)} holds. Taking the initial measure \(\delta_0\), we obtain the exponential
family of Theorem~\ref{thm:lifting}, so
\(K_1,\ldots,K_m\in\mathcal C\).
Set \(\by=\boldsymbol0\), and use the initial measure \(\delta_s\)
for every real \(s\ge0\).  Its initial series is \(A(t)=e^{st}\),
so the coefficients of \(e^{sL(t)}\) form a Stieltjes moment sequence
for every \(s\ge0\).
By the converse direction of Theorem~\ref{thm:lifting}, we now obtain
\(L\in\mathcal C\). We complete the proof.
\end{proof}

For formal power series \(B(t),L(t)\) with \(L(0)=0\), write
\(\mathcal S_{B,L}A=e^{B(t)}A(L(t))\).
These operators satisfy
\begin{equation}\label{eq:substitution-composition-law}
 \mathcal S_{B,L}\mathcal S_{H,R}
 =\mathcal S_{B+H\circ L,\,R\circ L}.
\end{equation}
The class \(\mathcal C\) is closed under nonnegative linear
combinations, by adding representing measures, and under composition
by Theorem~\ref{thm:composition}.
Thus the admissible operators are closed under finite composition.
With the substitution series fixed, independent convolution
intensities may be retained in each factor. We obtain joint coefficientwise positivity by
Theorem~\ref{thm:substitution}.
Only finite coefficient identities are involved.

\subsection{Moment inputs for component exponents}

The next construction turns any Stieltjes moment sequence
into admissible component exponents. Its specializations organize
the tree and forest applications below.

\begin{theorem}[Moment inputs and universal component families]
\label{thm:moment-input-families}
Let \(\boldsymbol b=(b_n)_{n\ge0}\) be a real Stieltjes moment
sequence, represented by a nonnegative measure \(\nu\) on
\([0,\infty)\) with all moments finite. Define the formal series
\begin{equation}\label{eq:moment-input-exponents}
 \begin{aligned}
 I_{\boldsymbol b}(t)&=\sum_{n\ge0}b_n\frac{t^{n+1}}{(n+1)!},\\
 J_{\boldsymbol b}(t)&=\sum_{n\ge1}b_n\frac{t^n}{n!},\\
 O_{\boldsymbol b}(t)&=t\sum_{n\ge0}b_nt^n.
 \end{aligned}
\end{equation}
Then we have:
\begin{enumerate}
\renewcommand{\labelenumi}{\textup{(\roman{enumi})}}
\item All three formal series in (\ref{eq:moment-input-exponents}) belong to \(\mathcal C\);
\item For Stieltjes sequences \(\boldsymbol b^{(\ell)}\), choose
\(K_\ell\in\{I_{\boldsymbol b^{(\ell)}},
J_{\boldsymbol b^{(\ell)}},O_{\boldsymbol b^{(\ell)}}\}\).
For fixed \(q_\ell\ge0\), \(S\in\mathcal C\), and any initial
Stieltjes moment series \(A(t)\), the exponential coefficients of
\begin{equation}\label{eq:moment-input-universal}
 A(S(t))\exp\!\left(\sum_{\ell=1}^m z_\ell K_\ell(t)\right)
          \prod_{\ell=1}^m(1-q_\ell K_\ell(t))^{-y_\ell}
\end{equation}
are coefficientwise Hankel-totally positive jointly in the
 indeterminates \(\by,\boldsymbol z\).
\end{enumerate}

\end{theorem}

\begin{proof}
The shifted moments of \(I_{\boldsymbol b}\) and
\(J_{\boldsymbol b}\) are represented by \(\nu\) and
\(x\nu(\dif x)\), respectively. Both measures have all
moments finite. For \(O_{\boldsymbol b}\), applying Euler's
integral~\cite[Eq.~5.2.1]{DLMF-Gamma} to the gamma factor,
we obtain
\begin{equation}\label{eq:moment-input-gamma-product}
 (n+1)!b_n
 =\int_0^\infty\!\int_0^\infty
       (gx)^n g e^{-g}\,\dif g\,\nu(\dif x).
\end{equation}
The image of \(\nu(\dif x)g e^{-g}\dif g\) under
\((x,g)\mapsto xg\) therefore represents the third shifted
sequence. This proves all three admissibility assertions.

By Proposition~\ref{powerapp:prop:ordinary-log} and Theorem~\ref{thm:composition}, respectively, we obtain \(-\log(1-q_\ell t)\in\mathcal C\) and \(-\log(1-q_\ell K_\ell)\in\mathcal C\).
Thus, applying Theorem~\ref{thm:substitution} to these exponents,
\(K_\ell\), and substitution \(S\), we obtain positivity for
\eqref{eq:moment-input-universal}. 
\end{proof}

\begingroup
\color{black}
\hypersetup{allcolors=black}
\section{Functional digraphs and colored components}\label{sec:digraph}

We apply Theorems~\ref{thm:lifting} and~\ref{thm:hamburger}
to functional digraphs, iterated tree series, and Gaussian
perturbations. In Sections~\ref{sec:digraph}--\ref{sec:applications},
an \emph{initial Stieltjes moment series} means a formal series
\(A(t)=\sum_{n\ge0}a_nt^n/n!\) whose real coefficient sequence
is Stieltjes; its total mass \(a_0\) need not be one. Hamburger
initial sequences are specified explicitly. Unless stated
otherwise, structural parameters are fixed, and coefficientwise
positivity refers to the indicated independent indeterminates.

\subsection{Connected-component moments and colored functional digraphs}

To verify the moment hypothesis for the connected coefficients
\(c_n(q)\) in \eqref{eq:connected-egf}, we use the second assertion
of Theorem~\ref{prop:connected-stieltjes}.
Writing the rooted-forest numbers as
\begin{equation}\label{eq:display-030}
 f_{n,k}=\binom{n-1}{k-1}n^{n-k}\quad(1\le k\le n),
\end{equation}
we obtain the coefficient identity
\begin{equation}\label{eq:connected-coefficients}
 c_n(q)=\sum_{k=1}^n(k-1)!f_{n,k}q^k;
\end{equation}
as in the expansion used in
\cite[Proposition~6.18(iii)]{Zhu2024}.
By the component interpretation, we obtain \eqref{eq:colored-egf}
with exponent \(\sum_\ell y_\ell C_{q_\ell}(t)\).
Comparing this identity with the colored partition formula
\eqref{eq:bell-partitions}, we identify the cumulants as
\(\kappa_n^{(\ell)}=c_n(q_\ell)\).

\begin{proof}[Proof of Theorem~\ref{thm:digraph-scalar}]
Applying Theorem~\ref{prop:connected-stieltjes} at each fixed
\(q_\ell\ge0\), we obtain the Stieltjes sequences
\((c_{n+1}(q_\ell))_{n\ge0}\). We identify these as the
once-shifted cumulants in \eqref{eq:colored-egf}.
Applying Theorem~\ref{thm:lifting},
\textup{(i)}\(\Rightarrow\)\textup{(ii)}, to these inputs,
we prove Theorem~\ref{thm:digraph}.
\end{proof}

We may introduce a common component parameter by substituting
\(y_1=\cdots=y_m=y\), since this substitution preserves
coefficientwise nonnegativity.
Setting \(q_1=\cdots=q_m=q\), we also obtain the identity
\begin{equation}\label{eq:display-031}
 \Psi_n((q,\ldots,q);\by)=\psi_n(q,y_1+\cdots+y_m).
\end{equation}
To illustrate the multivariate conclusion, put
\begin{equation}\label{eq:display-032}
 U=\sum_{\ell=1}^{m}q_\ell y_\ell,\qquad
 V=\sum_{\ell=1}^{m}(2q_\ell+q_\ell^2)y_\ell.
\end{equation}
Then \(\Psi_0=1\), \(\Psi_1=U\), \(\Psi_2=U^2+V\), and
\begin{equation}\label{eq:display-033}
 \det\begin{pmatrix}\Psi_0&\Psi_1\\\Psi_1&\Psi_2\end{pmatrix}
 =V\in\R_{\ge0}[\by].
\end{equation}
For \(m=1\), \(q_1=1\), we recover the original polynomials
\begin{equation}\label{eq:display-034}
 \psi_1(y)=y,\qquad
 \psi_2(y)=y^2+3y,\qquad
 \psi_3(y)=y^3+9y^2+17y,
\end{equation}
and
\begin{equation}\label{eq:display-035}
 \psi_4(y)=y^4+18y^3+95y^2+142y,\qquad
 \det\begin{pmatrix}\psi_1(y)&\psi_3(y)\\
                     \psi_2(y)&\psi_4(y)\end{pmatrix}
 =y^2(6y^2+51y+91).
\end{equation}

\subsection{General connected series and iterated tree functions}

The following logarithmic composition extends the functional-digraph
construction to arbitrary exponents in \(\mathcal C\) in \eqref{eq:admissible-class}.

\begin{corollary}[Logarithmic and iterated-tree families]\label{thm:logarithmic}
Fix \(L_1,\ldots,L_m\in\mathcal C\) and \(q_1,\ldots,q_m\ge0\).
Let \(A(t)\) be any initial Stieltjes moment series.
Define \(F_n(\by)\) by
\begin{equation}\label{eq:general-logarithmic-egf}
 \sum_{n\ge0}F_n(\by)\frac{t^n}{n!}
 =A(t)\prod_{\ell=1}^{m}(1-q_\ell L_\ell(t))^{-y_\ell}.
\end{equation}
Then \((F_n(\by))_{n\ge0}\) is coefficientwise Hankel-totally positive
jointly in \(y_1,\ldots,y_m\).
\end{corollary}

\begin{proof}
Put \(G_q(w)=-\log(1-qw)\). For \(q>0\), applying Euler's gamma
integral~\cite[Eq.~5.2.1]{DLMF-Gamma} to the shifted coefficients,
we obtain
\begin{equation}\label{eq:display-036}
 G_q^{(n+1)}(0)=n!q^{n+1}
       =\int_0^\infty x^n e^{-x/q}\,\dif x
       \qquad(n\ge0),
\end{equation}
Thus \(G_q\in\mathcal C\); for \(q=0\), its shifted sequence
is represented by the zero measure. Applying Theorem~\ref{thm:composition} to
\(G_{q_\ell}\) and \(L_\ell\), we obtain
\(G_{q_\ell}\circ L_\ell\in\mathcal C\).
Applying Theorem~\ref{thm:lifting} to these exponents and \(A\),
we obtain the claimed positivity for \eqref{eq:general-logarithmic-egf}.
\end{proof}

To include iterated Cayley tree functions in this construction, we
first verify their admissibility.

\begin{proposition}\label{prop:tree-admissible}
The tree function \(T\) belongs to \(\mathcal C\).
Consequently, \(T^{\circ d}\in\mathcal C\) for every integer \(d\ge1\).
\end{proposition}

\begin{proof}
For every \(n\ge1\), expanding \eqref{eq:connected-coefficients}
at zero, we obtain
\begin{equation}\label{eq:display-037}
 \lim_{q\downarrow0}\frac{c_n(q)}q
 =f_{n,1}=n^{n-1}.
\end{equation}
For each \(q>0\), we know from Theorem~\ref{prop:connected-stieltjes}
that \((c_{n+1}(q)/q)_{n\ge0}\) is a Stieltjes moment sequence.
Passing to the limit in each finite Hankel minor, we obtain a
nonnegative minor for \(((n+1)^n)_{n\ge0}\). Hence this sequence
is a Stieltjes moment sequence.
Since \(T^{(n+1)}(0)=(n+1)^n\), we obtain \(T\in\mathcal C\).
Applying Theorem~\ref{thm:composition} successively to these
exponents, we obtain the assertion for all iterates.
For the Stieltjes moment property of \(((n+1)^n)\), see also
\cite{Sokal2023,Zhu2024}.
\end{proof}

Taking \(L_\ell=T^{\circ d_\ell}\) for \(d_\ell\ge1\), in
Corollary~\ref{thm:logarithmic}, we obtain the iterated-tree families.
Taking \(A=1\) and \(d_\ell=1\), we recover
Theorem~\ref{thm:digraph}. As before, we may identify
variables by substitution.

For example, take the initial measure
\(\tfrac12\delta_0+\tfrac12\delta_1\), with formal generating function
\begin{equation}\label{eq:display-038}
 A(t)=\frac{1+e^t}{2},
\end{equation}
and, applying Corollary~\ref{thm:logarithmic} to the two tree
exponents, we obtain coefficientwise Hankel total positivity for
\begin{equation}\label{eq:bernoulli-tree-example}
 \sum_{n\ge0}F_n(y,z)\frac{t^n}{n!}
 =\frac{1+e^t}{2}\,
       (1-T(T(t)))^{-y}(1-2T(t))^{-z}.
\end{equation}
For this family, the first moments and the first nontrivial leading
Hankel minor are
\begin{equation}\label{eq:display-039}
 \begin{aligned}
 F_0&=1,\qquad F_1=\tfrac12+y+2z,\\
 F_2&=(y+2z)^2+6y+10z+\tfrac12,\qquad
 F_0F_2-F_1^2=\tfrac14+5y+8z.
 \end{aligned}
\end{equation}

The initial moment series need not be the exponential of an
admissible exponent. For \(H(t)=\log((1+e^t)/2)\), we ahve
\begin{equation}\label{eq:display-040}
 H'(0)=\tfrac12,\qquad H''(0)=\tfrac14,\qquad H'''(0)=0,
 \qquad
 \det\begin{pmatrix}H'(0)&H''(0)\\H''(0)&H'''(0)\end{pmatrix}
 =-\tfrac1{16}.
\end{equation}
Thus \(H\notin\mathcal C\), although its exponential is an
allowed initial moment series in Theorem~\ref{thm:lifting}.

\subsection{Gaussian perturbations and centered component families}

We next use the real-line theorem to allow an arbitrary initial
Hamburger measure, real centering drifts, and an independent Gaussian
component in the functional-digraph construction.

\begin{corollary}[Centered and Gaussian-perturbed tree families]\label{cor:centered-digraph}
Fix integers \(d_\ell\ge1\), real numbers
\(q_\ell\ge0\), \(\beta_\ell\in\R\), and
\(\sigma^2\ge0\).  Let \(A(t)\) be the formal exponential
generating function of any initial Hamburger moment sequence
\((a_n)\).  Define \(G_n(\by,z)\) by
\begin{equation}\label{eq:centered-digraph-egf}
 \sum_{n\ge0}G_n(\by,z)\frac{t^n}{n!}
 =A(t)\exp\!\left(\frac{\sigma^2zt^2}{2}
                     -t\sum_{\ell=1}^{m}\beta_\ell y_\ell\right)
       \prod_{\ell=1}^{m}
          (1-q_\ell T^{\circ d_\ell}(t))^{-y_\ell}.
\end{equation}
For \(r\ge1\), we have
\begin{equation}\label{eq:display-041}
 H_r(\by,z):=\det[G_{i+j}(\by,z)]_{i,j=0}^{r-1}
          \in\R_{\ge0}[\by,z].
\end{equation}
With \(D=\binom r2\), its total degree is at most \(D\), and
\begin{equation}\label{eq:centered-digraph-top-degree}
 [H_r]_D=a_0^r\!\left(\prod_{j=0}^{r-1}j!\right)
       \left(\sigma^2z+
          \sum_{\ell=1}^{m}(2d_\ell q_\ell+q_\ell^2)y_\ell
       \right)^D.
\end{equation}
The determinant is independent of \(\beta_1,\ldots,\beta_m\).
\end{corollary}

\begin{proof}
Put \(C_\ell=-\log(1-q_\ell T^{\circ d_\ell})\).
Applying Corollary~\ref{thm:logarithmic} to the iterates from
Proposition~\ref{prop:tree-admissible}, we obtain \(C_\ell\in\mathcal C\).
If \(\nu_\ell\) represents its once-shifted cumulants, then
\(x\,\nu_\ell(\dif x)\) represents its twice-shifted cumulants.
Subtracting \(\beta_\ell\, t\) leaves these unchanged.
The Gaussian exponent has twice-shifted measure
\(\sigma^2\delta_0\). Applying Theorem~\ref{thm:hamburger}
to this exponent list and \(A\), we obtain \eqref{eq:display-041},
the degree bound, and drift invariance. Using the chain rule,
we obtain
\((T^{\circ d})'(0)=1\), \((T^{\circ d})''(0)=2d\), hence
\(C_\ell''(0)=2d_\ell q_\ell+q_\ell^2\).
Substituting these derivatives in \eqref{eq:hamburger-top-degree},
we obtain \eqref{eq:centered-digraph-top-degree}.
\end{proof}

Choosing \(\beta_\ell=q_\ell\) centers each added component family. So
the first moment of the output is then \(G_1=a_1\).
For \(\by\ge\boldsymbol0\) and \(z\ge0\), we also obtain a
Hamburger moment sequence. In the coefficientwise assertion, we
consider leading principal determinants, as in the real-line theorem.

\section{Implicit and autonomous equations}\label{sec:equations}

We will first characterize the implicit equations whose solutions
belong to \(\mathcal C\). The criteria
\cite[Theorem~4.2]{Zhu2026Riordan} and
\cite[Theorem~6.6]{Zhu2024} then supply moment inputs for
universal families with independent colors and powers of
prefactors.

\subsection{Implicit tree functions and an exact lifting criterion}
\label{sec:implicit-trees}

We now replace the Cayley tree equation (\ref{eq:tree-digraph-egf}) by a general implicit
equation. Let \(\Phi(u)\in\R[[u]]\) with \(\Phi(0)>0\), and define
the unique formal power series \(U(t)\in t\R[[t]]\) by
\begin{equation}\label{eq:implicit-tree-equation}
 U(t)=t\Phi(U(t)).
\end{equation}
When \(\Phi(u)=\sum_{k\ge0}w_ku^k/k!\) with \(w_k\ge0\), we
interpret \(U\) as the exponential generating function of labeled
rooted trees with weight \(w_k\) at a vertex with \(k\) children.
To compute the coefficients, we recall the classical Lagrange
inversion formula; see Stanley~\cite[Section~5.4]{Stanley}.

\begin{lemma}[Lagrange inversion]\label{lem:lagrange-source}
Under \eqref{eq:implicit-tree-equation}, for \(n\ge1\) and
\(H(u)\in\R[[u]]\), we have
\begin{equation}\label{eq:lagrange-source}
 [t^n]H(U(t))=\frac1n[u^{n-1}]H'(u)\Phi(u)^n.
\end{equation}
In particular, if \(U(t)=\sum_{n\ge1}\tau_nt^n/n!\), then
\begin{equation}\label{eq:implicit-tree-shifted-moments}
 \tau_{n+1}=n![u^n]\Phi(u)^{n+1}\qquad(n\ge0).
\end{equation}
\end{lemma}

Applying Theorem~\ref{thm:lifting} to the shifted coefficients (\ref{eq:implicit-tree-shifted-moments}), we obtain the following exact criterion for logarithmic
components of an implicit tree series.

\begin{theorem}[Implicit-tree lifting criterion]\label{thm:implicit-tree-criterion}
For \(U\) in \eqref{eq:implicit-tree-equation}, the following
conditions are equivalent:
\begin{enumerate}
\renewcommand{\labelenumi}{\textup{(\roman{enumi})}}
\item The sequence
\(\bigl(n![u^n]\Phi(u)^{n+1}\bigr)_{n\ge0}\) is a Stieltjes moment sequence.
\item For every fixed \(q>0\), the sequence $(\Psi_n^{U}(q,y))_{n}$ defined by
\begin{equation}\label{eq:implicit-component-family}
 \sum_{n\ge0}\Psi_n^{U}(q,y)\frac{t^n}{n!}
                 =(1-qU(t))^{-y}
\end{equation}
is coefficientwise Hankel-totally positive in \(y\).
\item For every \(m\ge1\), fixed \(q_1,\ldots,q_m\ge0\), and
real initial Stieltjes moment series \(A\), the sequence defined by
\begin{equation}\label{eq:implicit-colored-lifting}
 \sum_{n\ge0}F_n(\by)\frac{t^n}{n!}
       =A(t)\prod_{\ell=1}^m(1-q_\ell U(t))^{-y_\ell}
\end{equation}
is coefficientwise Hankel-totally positive in \(\by\).
\end{enumerate}
\end{theorem}

\begin{proof}
By \eqref{eq:implicit-tree-shifted-moments}, condition
\textup{(i)} is equivalent to \(U\in\mathcal C\).
Applying Corollary~\ref{thm:logarithmic} to \(U\), we obtain
\textup{(i)}\(\Rightarrow\)\textup{(iii)}. Taking \(A=1\)
and \(m=1\), we obtain
\textup{(iii)}\(\Rightarrow\)\textup{(ii)}.

Assume \textup{(ii)} holds and let \(C_q^U=-\log(1-qU)\).
Applying Theorem~\ref{thm:scalar-lifting} to \(C_q^U\) and
scaling the representing measure by \(1/q\), we obtain the
Stieltjes sequence \(((C_q^U)^{(n+1)}(0)/q)_{n\ge0}\)
for every \(q>0\). Expanding the formal logarithm, we have
\begin{equation}\label{eq:small-cycle-weight-limit}
 \frac{C_q^U(t)}q
 =U(t)+\sum_{k\ge2}\frac{q^{k-1}}kU(t)^k,
 \qquad
 \lim_{q\downarrow0}\frac{(C_q^U)^{(n+1)}(0)}q=\tau_{n+1}.
\end{equation}
Each finite Hankel minor is a polynomial in finitely many
entries, so its limit is nonnegative. Applying
Lemma~\ref{lem:moment-criteria} to the limiting sequence,
we obtain \textup{(i)}.
\end{proof}

To show why the moment condition is needed, we give an example
with nonnegative tree weights for which the positivity assertion fails.

\begin{rem}[Nonnegative branching weights alone]\label{rem:implicit-tree-counterexample}
Take \(\Phi(u)=1+u^2\). We obtain
\(U(t)=t+t^3+2t^5+\cdots\), and from
\eqref{eq:implicit-tree-shifted-moments} we have
\(\tau_1=1,\tau_2=0,\tau_3=6,\tau_4=0\). Hence
\begin{equation}\label{eq:implicit-tree-counterexample}
 \det\begin{pmatrix}\tau_2&\tau_3\\\tau_3&\tau_4\end{pmatrix}
 =-36.
\end{equation}
We can also verify the failure directly in \eqref{eq:implicit-component-family}:
for \(q=1\),
\begin{equation}\label{eq:implicit-component-counterexample}
 [y^2]\left(\Psi_2^{U}(1,y)\Psi_4^{U}(1,y)
                 -\Psi_3^{U}(1,y)^2\right)=-34.
\end{equation}
Thus we must check the Stieltjes condition in
Theorem~\ref{thm:implicit-tree-criterion}, even when all branching
weights are nonnegative.
\end{rem}

\subsection{Implicit equations with P\'olya frequency series}
\label{sec:pf-implicit-lifting}

A sufficient moment condition from implicit-equation
criterion was showed in \cite{Zhu2026Riordan}. For the related branched-continued-fraction approach,
see P\'etr\'eolle, Sokal and Zhu~\cite[Theorems~9.8 and~9.9]{PSZ2023}.
We distinguish ordinary coefficients in this input from the
exponential coefficients required by lifting.

For a real formal power series \(h(u)=\sum_{n\ge0}h_nu^n\), put
\(h_n=0\) for \(n<0\). We call \(h\) a \emph{P\'olya frequency
series}, or a \(\mathrm{PF}\) series, if
the Toeplitz matrix \([h_{i-j}]_{i,j\ge0}\) is totally positive.
In particular, \(h_n\ge0\). By the Cauchy--Binet formula, we
preserve this property when multiplying two such series, since the
corresponding lower triangular Toeplitz matrices multiply by ordinary
convolution.

For the ordinary-coefficient input, we use the real,
infinite-order case of the criterion
\cite[pp.~19--20, Theorems~4.2 and~4.2\(^{\prime}\),
Eq.~(4.12)]{Zhu2026Riordan}. We record its specialization and
coefficient convention explicitly as follows:

\begin{lemma}[Implicit-equation criterion]
\label{lem:zhu-implicit-source}
Let \(\varphi(u),\psi(u)\in\R[[u]]\) be
\(\mathrm{PF}\) series with
\(\varphi(0)\psi(0)>0\), and let
\(U(t)=t\varphi(U(t))\psi(U(t))\). Define
\begin{equation}\label{eq:zhu-implicit-ordinary}
 R_n(q):=[t^n]\frac{\psi(U(t))}{1-qU(t)}\qquad(n\ge0).
\end{equation}
Then \((R_n(q))_{n\ge0}\) is coefficientwise Hankel-totally
positive in \(q\). In particular, for $n\geq0$,
\([t^n]\psi(U(t))\) form a real Stieltjes moment
sequence \((b_n)_{n\ge0}\). 
\end{lemma}

Applying the lifting criterion to this moment input, we obtain
the following universal implicit-equation family.

\begin{theorem}[Universal lifting for implicit P\'olya frequency equations]
\label{thm:pf-implicit-lifting}
For \(1\le\ell\le m\), fix real \(\mathrm{PF}\) series
\(\varphi_\ell(u),\psi_\ell(u)\) with
\(\varphi_\ell(0)\psi_\ell(0)>0\), and define
\begin{equation}\label{eq:pf-implicit-equations}
 U_\ell(t)=t\varphi_\ell(U_\ell(t))\psi_\ell(U_\ell(t)).
\end{equation}
Then the following assertions hold.
\begin{enumerate}
\renewcommand{\labelenumi}{\textup{(\roman{enumi})}}
\item Each \(U_\ell\in\mathcal C\); equivalently, with
\(\Phi_\ell=\varphi_\ell\psi_\ell\), the sequence
\begin{equation}\label{eq:pf-implicit-shifted}
 \left(n![u^n]\Phi_\ell(u)^{n+1}\right)_{n\ge0}
 =\left(U_\ell^{(n+1)}(0)\right)_{n\ge0}
\end{equation}
is a Stieltjes moment sequence.
\item For any initial Stieltjes moment series \(A(t)\) and
fixed \(q_1,\ldots,q_m\ge0\), define
\begin{equation}\label{eq:pf-implicit-universal-family}
 \sum_{n\ge0}F_n(\bx,\by)\frac{t^n}{n!}
 =A(t)\prod_{\ell=1}^m
 \left\{\psi_\ell(U_\ell(t))e^{x_\ell U_\ell(t)}
                  (1-q_\ell U_\ell(t))^{-y_\ell}\right\}.
\end{equation}
Then \((F_n)\) is jointly coefficientwise Hankel-totally
positive in \(\bx,\by\).
We may omit any of the prefactors \(\psi_\ell(U_\ell(t))\)
and retain the same conclusion.
\end{enumerate}
\end{theorem}

\begin{proof}
\emph{Step 1: Ordinary shifted moments.}
Suppress the color index and put \(\Phi=\varphi\psi\).
Applying the Toeplitz product rule to \(\varphi\) and \(\psi\),
we obtain \(\Phi\in\mathrm{PF}\).
Applying Lemma~\ref{lem:zhu-implicit-source} to the factorization
\(\Phi=\Phi\cdot1\), we obtain the resolvent polynomials
\begin{equation}\label{eq:pf-implicit-resolvent-linear}
 Q_n(q):=[t^n](1-qU(t))^{-1},\qquad
 Q_n(q)=q[t^n]U(t)+O(q^2)\quad(n\ge1).
\end{equation}
Here \(O(q^2)\) denotes divisibility by \(q^2\).
For increasing nonnegative index lists \(I,J\) of length \(r\),
each entry of the shifted minor has zero constant term. Hence
\begin{equation}\label{eq:pf-implicit-linear-minor}
 [q^r]\det[Q_{i_a+j_b+1}(q)]_{a,b=1}^r
 =\det\bigl[[t^{i_a+j_b+1}]U(t)\bigr]_{a,b=1}^r\ge0.
\end{equation}
Thus \(d_n=[t^{n+1}]U(t)\) has nonnegative Hankel minors and,
by Lemma~\ref{lem:moment-criteria}, a Stieltjes measure \(\nu\).

\smallskip\noindent
\emph{Step 2: Exponential coefficients and the prefactor.}
Applying Euler's integral to the product measure, we obtain
\begin{equation}\label{eq:pf-implicit-factorial-moments}
 U^{(n+1)}(0)=(n+1)!d_n
 =\int_0^\infty\!\int_0^\infty
             (sx)^n s e^{-s}\,\dif s\,\nu(\dif x).
\end{equation}
Taking the image of \(\nu(\dif x)s e^{-s}\dif s\) under
\((x,s)\mapsto sx\), we represent the once-shifted cumulants
and obtain \(U\in\mathcal C\). Applying Lagrange inversion
to \(U=t\Phi(U)\), we obtain \eqref{eq:pf-implicit-shifted}.

Applying Lemma~\ref{lem:zhu-implicit-source} to the original
factorization at \(q=0\) yields a measure \(\sigma\)
representing \(b_n=[t^n]\psi(U(t))\). We then have
\begin{equation}\label{eq:pf-implicit-prefactor-moments}
 n!b_n=\int_0^\infty\!\int_0^\infty
                    (sx)^n e^{-s}\,\dif s\,\sigma(\dif x)
 \qquad(n\ge0).
\end{equation}
These are the exponential coefficients of the prefactor.
The distinct gamma weights in
\eqref{eq:pf-implicit-factorial-moments} and
\eqref{eq:pf-implicit-prefactor-moments} account for the
factors \((n+1)!\) and \(n!\), respectively.

\smallskip\noindent
\emph{Step 3: Initial moments and lifting.}
For moment series represented by \(\mu_1,\mu_2\), binomial
convolution is the image measure under addition:
\begin{equation}\label{eq:pf-prefactor-binomial-convolution}
 \sum_{k=0}^n\binom nk
     \left(\int_0^\infty x^k\,\mu_1(\dif x)\right)
     \left(\int_0^\infty z^{n-k}\,\mu_2(\dif z)\right)
 =\int_0^\infty\!\int_0^\infty(x+z)^n
                    \,\mu_1(\dif x)\mu_2(\dif z).
\end{equation}
In consequence,
\(B(t)=A(t)\prod_\ell\psi_\ell(U_\ell(t))\) is an initial
Stieltjes moment series. By Step~2 and
Corollary~\ref{thm:logarithmic}, both \(U_\ell\) and
\(-\log(1-q_\ell U_\ell)\) belong to \(\mathcal C\).
Applying Theorem~\ref{thm:lifting} to these \(2m\) exponents
and \(B\), we obtain \textup{(ii)}. Omitting prefactors simply
removes the corresponding factors from \(B\).
\end{proof}

For polynomial equation coefficients, we may also prescribe a
nonzero initial value. For the corresponding ordinary Riordan-array
results, see \cite[Theorem~6.1]{Zhu2026Riordan}. We now give the
universal cumulant-lifting counterpart.

\begin{corollary}[A nonzero polynomial initial value]
\label{cor:pf-implicit-nonzero}
Fix \(\kappa\ge0\) and
\(\rho(u)=\prod_{h=1}^d(\alpha_h+\beta_hu)\), where
 \(\alpha_h,\beta_h\ge0\), and
\(\rho(\kappa)>0\). Define the unique formal power series
\begin{equation}\label{eq:pf-nonzero-initial-equation}
 f(t)=\kappa+t\rho(f(t)),\qquad U(t)=f(t)-\kappa.
\end{equation}
Then \(U\in\mathcal C\). In addition, for every real initial Stieltjes moment
series \(A(t)\) and fixed \(q_1,\ldots,q_m\ge0\), the
exponential coefficients of
\begin{equation}\label{eq:pf-nonzero-initial-lifting}
 A(t)\exp\!\left(\sum_{\ell=1}^m x_\ell U(t)\right)
                     \prod_{\ell=1}^m(1-q_\ell U(t))^{-y_\ell}
\end{equation}
are coefficientwise Hankel-totally positive jointly in \(\bx,\by\).
\end{corollary}

\begin{proof}
Subtracting \(\kappa\) in
\eqref{eq:pf-nonzero-initial-equation}, we obtain
\begin{equation}\label{eq:pf-nonzero-shifted-equation}
 U(t)=t\widetilde\rho(U(t)),\qquad
 \widetilde\rho(u)=\prod_{h=1}^d
                         (\alpha_h+\kappa\beta_h+\beta_hu).
\end{equation}
The Toeplitz matrix of each nonnegative linear factor is bidiagonal
and totally positive, by Lemma~\ref{lem:matrix-tools}. By the
Toeplitz product rule, \(\widetilde\rho\) is
\(\mathrm{PF}\), and
\(\widetilde\rho(0)=\rho(\kappa)>0\).
Applying Theorem~\ref{thm:pf-implicit-lifting} to the factorization
\(\widetilde\rho=\widetilde\rho\cdot1\), we obtain the result.
\end{proof}

In this corollary, we substitute into a polynomial. For a general formal
power series \(\rho(u)\), substitution at a nonzero \(\kappa\)
requires a separate definition or suitable convergence hypotheses; see
\cite[Remark~1.2]{Zhu2026Riordan}.

\subsection{Autonomous differential equations with P\'olya frequency series}
\label{sec:pf-autonomous-lifting}

We next obtain an analogous universal result for autonomous
differential equations. As the known input, we use
\cite[Theorem~6.6 (iv) and~(iii)]{Zhu2024} with real coefficients.
If \(\overline V\) denotes the formal compositional inverse, then
\(1/\overline V'(u)=\varphi(u)\psi(u)\) is equivalent to
\(V'(t)=\varphi(V(t))\psi(V(t))\). We therefore restate the criterion in autonomous-equation form
and refer to the cited work for its proof.

\begin{lemma}[Autonomous-equation criterion]
\label{lem:zhu-autonomous-source}
Let \(\varphi(u),\psi(u)\in\R[[u]]\) be
\(\mathrm{PF}\) series with
\(\varphi(0)\psi(0)>0\). Let \(V(0)=0\) and
\begin{equation}\label{eq:zhu-autonomous-source}
 V'(t)=\varphi(V(t))\psi(V(t)).
\end{equation}
Then \((V^{(n+1)}(0))_{n\ge0}\) is a Stieltjes moment sequence.
Moreover, the exponential coefficients of
\(\varphi(V(t))e^{\lambda V(t)}\) form a coefficientwise
Hankel-totally positive sequence in \(\lambda\).
In particular, the exponential coefficients of
\(\varphi(V(t))\) form a real Stieltjes moment sequence.
\end{lemma}

Applying the autonomous criterion to the equations below,
we obtain the parallel universal family.

\begin{theorem}[Universal lifting for autonomous P\'olya frequency equations]
\label{thm:pf-autonomous-lifting}
For \(1\le\ell\le m\), fix real \(\mathrm{PF}\)
series \(\varphi_\ell(u),\psi_\ell(u)\) with
\(\varphi_\ell(0)\psi_\ell(0)>0\), and define the unique formal
power series \(V_\ell(t)\) by
\begin{equation}\label{eq:pf-autonomous-equations}
 V_\ell'(t)=\varphi_\ell(V_\ell(t))\psi_\ell(V_\ell(t)),
 \qquad V_\ell(0)=0.
\end{equation}
Then the following assertions hold:
\begin{enumerate}
\renewcommand{\labelenumi}{\textup{(\roman{enumi})}}
\item \(V_\ell\in\mathcal C\) for every \(\ell\).
\item For any initial Stieltjes moment series \(A(t)\) and
fixed \(q_1,\ldots,q_m\ge0\), define
\begin{equation}\label{eq:pf-autonomous-universal-family}
 \sum_{n\ge0}F_n(\bx,\by)\frac{t^n}{n!}
 =A(t)\prod_{\ell=1}^m
 \left\{\varphi_\ell(V_\ell(t))e^{x_\ell V_\ell(t)}
                  (1-q_\ell V_\ell(t))^{-y_\ell}\right\}.
\end{equation}
Then \((F_n)\) is coefficientwise Hankel-totally
positive jointly in \(\bx,\by\).
We may omit any of the prefactors \(\varphi_\ell(V_\ell(t))\)
and retain the same conclusion.
\end{enumerate}
\end{theorem}

\begin{proof}
Applying Lemma~\ref{lem:zhu-autonomous-source} to each equation,
we obtain \(V_\ell\in\mathcal C\) and a Stieltjes moment
series \(\varphi_\ell(V_\ell(t))\). By
\eqref{eq:pf-prefactor-binomial-convolution},
\(B(t)=A(t)\prod_\ell\varphi_\ell(V_\ell(t))\) is also a
Stieltjes moment series. Applying Corollary~\ref{thm:logarithmic} to \(V_\ell\),
we obtain \(-\log(1-q_\ell V_\ell)\in\mathcal C\).
Hence, by Theorem~\ref{thm:lifting},
we obtain the asserted joint positivity.
The same argument allows any prefactors to be omitted.
\end{proof}

\subsection{Independent powers of Riordan prefactors}
\label{powerapp:sec:riordan}

We now retain independent powers of admissible prefactors.
Applying the classical logarithm criterion in
Proposition~\ref{powerapp:prop:ordinary-log} to the ordinary moment
inputs below, we obtain the exponents needed for the implicit
case.

We apply Lemmas~\ref{lem:zhu-implicit-source} and
\ref{lem:zhu-autonomous-source} to positive linear factors in
either equation.

\begin{theorem}[Independent powers and universal lifting]
\label{powerapp:thm:powers}
Fix integers \(p,r\ge0\), \(c_1,\ldots,c_p\ge0\) and
\(d_1,\ldots,d_r\ge0\). Let \(\Gamma(u)\) be a real
\(\mathrm{PF}\) series with \(\Gamma(0)>0\), and put
\begin{equation}\label{powerapp:eq:factor-equation}
 \Phi(u)=\Gamma(u)\prod_{i=1}^p(1+c_i u).
\end{equation}
Let \(U(t)\) be the unique formal power series with \(U(0)=0\)
satisfying either \(U(t)=t\Phi(U(t))\) or
\(U'(t)=\Phi(U(t))\). Then the following hold.
\begin{enumerate}
\renewcommand{\labelenumi}{\textup{(\roman{enumi})}}
\item For \(1\le i\le p\) and \(1\le j\le r\), we have
\begin{equation}\label{powerapp:eq:admissible-logs}
 \left\{U(t),\, \log(1+c_iU(t)),\,
              -\log(1-d_jU(t))\right\}\subseteq\mathcal C.
\end{equation}
\item For any initial Stieltjes moment series \(A(t)\) and
independent intensities \(y\), \(\boldsymbol z=(z_1,\ldots,z_p)\),
\(\boldsymbol w=(w_1,\ldots,w_r)\), define
\begin{equation}\label{powerapp:eq:independent-powers}
 \sum_{n\ge0}F_n(y,\boldsymbol z,\boldsymbol w)\frac{t^n}{n!}
 =A(t)e^{yU(t)}\frac{\prod_{i=1}^p(1+c_iU(t))^{z_i}}{
                         \prod_{j=1}^r(1-d_jU(t))^{w_j}}.
\end{equation}
Then \((F_n)\) is coefficientwise Hankel-totally
positive jointly in \(y,\boldsymbol z,\boldsymbol w\). The conclusion
also holds for finite products from different equations with
independent indeterminates for each equation.
\end{enumerate}
\end{theorem}

\begin{proof}
If \(c_i=0\), the logarithm is zero. For \(c_i>0\), write
\begin{equation}\label{powerapp:eq:remove-linear-factor}
 \Psi_i(u)=\Gamma(u)\prod_{h\ne i}(1+c_hu),\qquad
 \Phi(u)=(1+c_i u)\Psi_i(u).
\end{equation}
Both factors are \(\mathrm{PF}\).
For the implicit equation, applying
Lemma~\ref{lem:zhu-implicit-source} to
\(\varphi=\Psi_i\), \(\psi(u)=1+c_i u\), with resolvent
parameter zero, we obtain the Stieltjes ordinary coefficients
of \(1+c_iU\). We have positive product \(\Gamma(0)\)
of the source constant terms and constant term one for
\(1+c_iU\). Applying Propostion~\ref{powerapp:prop:ordinary-log}
to this ordinary moment series, we obtain
\(\log(1+c_iU)\in\mathcal C\).

For the autonomous equation,
\begin{equation}\label{powerapp:eq:autonomous-log-derivative}
 \frac{\dif}{\dif t}\log(1+c_iU(t))
 =\frac{c_i\Phi(U(t))}{1+c_iU(t)}=c_i\Psi_i(U(t)).
\end{equation}
Applying Lemma~\ref{lem:zhu-autonomous-source} to the
prefactor \(\Psi_i\) gives the Stieltjes exponential
coefficients of \(\Psi_i(U)\). Multiplying the representing
measure by \(c_i\), we derive the required shifted moments.

Applying Theorem~\ref{thm:pf-implicit-lifting} or
\ref{thm:pf-autonomous-lifting} to the given equation, we obtain
\(U\in\mathcal C\). From Corollary~\ref{thm:logarithmic}, we have \(-\log(1-d_jU)\in\mathcal C\).
Rewrite the product in \eqref{powerapp:eq:independent-powers}
as $$A(t)\exp(yU+\sum_i z_i\log(1+c_iU)-\sum_j w_j\log(1-d_jU)).$$
Hence, by Theorem~\ref{thm:lifting}, we obtain
the conclusion including finite collections of equations with
independent indeterminates.
\end{proof}

For the underlying examples, coefficientwise Hankel total
positivity is already established: for multivariate Fuss--Narayana
polynomials, see P\'etr\'eolle, Sokal and Zhu~\cite[Theorem~10.1]{PSZ2023};
for the Riordan versions, see Zhu~\cite[Theorem~5.1]{Zhu2026Riordan};
for higher Eulerian and network families, see
\cite[Proposition~6.22 and Theorem~6.25]{Zhu2024}.
We now apply Theorem~\ref{powerapp:thm:powers} to record the
corresponding consequences with independent powers and arbitrary
initial moments.

\begin{corollary}[Fuss--Narayana--Riordan and higher Eulerian families]
\label{powerapp:cor:examples}
Let \(A(t)\) be any initial Stieltjes moment series.
The exponential coefficient sequences below are jointly
coefficientwise Hankel-totally positive in the indicated intensities.
\begin{enumerate}
\renewcommand{\labelenumi}{\textup{(\roman{enumi})}}
\item For fixed \(m\ge1\), \(c_0,\ldots,c_m,q\ge0\), let
\(U(t)=t\prod_{i=0}^m(1+c_iU(t))\). Then we may take
\begin{equation}\label{powerapp:eq:narayana}
 A(t)e^{yU(t)}(1-qU(t))^{-w}
                         \prod_{i=0}^m(1+c_iU(t))^{z_i},
 \quad\text{indeterminates}y,w,z_0,\ldots,z_m.
\end{equation}
\item For fixed integer \(r\ge1\) and \(x,q\ge0\), let \(V'(t)=(1+V(t))(1+xV(t))^r\) and \(V(0)=0\).
Then we may take
\begin{equation}\label{powerapp:eq:eulerian}
 A(t)e^{yV(t)}(1+V(t))^u(1+xV(t))^v(1-qV(t))^{-w},
 \quad\text{indeterminates }y,u,v,w.
\end{equation}
\item For fixed \(a,b\ge0\), let
\(S'(t)=(1+aS(t))/(1-bS(t))\) and \(S(0)=0\). Then we may take
\begin{equation}\label{powerapp:eq:networks}
 A(t)e^{yS(t)}(1+aS(t))^u(1-bS(t))^{-v},
 \quad\text{indeterminates }y,u,v.
\end{equation}
\end{enumerate}
\end{corollary}

\begin{proof}
For \textup{(i)}, applying Theorem~\ref{powerapp:thm:powers}
to \(\Gamma=1\) and the \(m+1\) factors \(1+c_i u\),
we obtain the asserted family. In \cite{Zhu2026Riordan}, the ordinary Riordan family uses the same equation with
prefactors \(\prod_{i=0}^{s-1}(1+c_iU)\), where \(s\) is an
integer, \(0\le s\le m+1\), and the empty product is one
\cite[pp.~23--24, Eqs.~(5.1)--(5.2) and
Theorem~5.1]{Zhu2026Riordan}. Its ordinary row generating series
is this prefactor divided by \(1-qU\). The independent powers
and arbitrary initial series \(A\) in our conclusion follow
from the admissible logarithms in
Theorem~\ref{powerapp:thm:powers}.

For \textup{(ii)}, use the factors \(1+u\) and \(r\) copies
of \(1+xu\), then identify the intensities of the repeated
factors. The source equation and prefactors, with powers
\(i\in\{0,\ldots,r\}\), \(j\in\{0,1\}\), occur in
\cite[Section~6.4, Eqs.~(6.145)--(6.147) and
Proposition~6.22]{Zhu2024}.

For \textup{(iii)}, take \(\Gamma(u)=(1-bu)^{-1}\), which
is \(\mathrm{PF}\), and the factor \(1+au\).
Applying Theorem~\ref{powerapp:thm:powers} to these factors,
we obtain the result. The source
prefactor powers are \(i,j\in\{0,1\}\)
\cite[Eq.~(6.161) and Theorem~6.25]{Zhu2024}.
The choices \(a=b=1\) and \(a=b=2\) give the series-parallel
network and nondegenerate fanout-free function series,
respectively \cite[Section~6.6]{Zhu2024}.
\end{proof}

Here \(c_i,x,a,b,q\) remain fixed. The extension concerns
independent powers and arbitrary initial moments; it makes no
joint coefficientwise assertion in these structural parameters.

\section{Tree moment inputs and universal component families}\label{sec:tree-applications}

\subsection{Fuss--Raney tree functions and colored components}
\label{sec:raney-trees}

For integer \(p\ge1\), we count rooted trees with \(p\)
ordered child slots at each vertex, allowing slots to be empty.
For a distinguished root, we may instead prescribe \(r\) slots.
Giving each edge weight \(b\), we obtain the series
\(tB_p(bt)^r\). We use the known Raney moments to extend
this model to real \(p\ge1\), \(0\le r\le p\), and then
apply lifting to obtain a universal colored component family.
For the related multivariate Fuss-Narayana positivity, see
\cite[Theorem~10.1]{PSZ2023}. Define the unique formal series
\(B_p(z)\in1+z\R[[z]]\) by
\begin{equation}\label{eq:fuss-catalan-series}
 B_p(z)=1+zB_p(z)^p.
\end{equation}
We use the generalized binomial expansion for real powers. Define
\(R_{p,r}(0)=1\) and
\begin{equation}\label{eq:raney-numbers}
 R_{p,r}(n)=\frac r{n!}\prod_{j=1}^{n-1}(pn+r-j)
\end{equation}
for \(n\ge1\).
The following representation is due to
M\l otkowski~\cite[Eqs.~(3.2)--(3.3) and Theorem~5.1]{Mlotkowski2010}.

\begin{lemma}[Classical Raney moment measures]\label{lem:raney-source}
For \(p\ge1\) and \(0\le r\le p\), there is a compactly supported
probability measure \(\mu_{p,r}\) on \([0,\infty)\) such that
\begin{equation}\label{eq:raney-moment-measure}
 \int_0^\infty  x^n\,d\mu_{p,r}(x)=R_{p,r}(n),\qquad
 B_p(z)^r=\sum_{n\ge0}R_{p,r}(n)z^n.
\end{equation}
For \(r=0\), we take \(\mu_{p,0}=\delta_0\).
\end{lemma}

To obtain the shifted exponential coefficients required for lifting
from these ordinary moment coefficients, we introduce an independent
gamma factor of shape two.
For \(\alpha>0\), we write
\(\Gamma(\alpha)=\int_0^\infty x^{\alpha-1}e^{-x}\,\dif x\).
For \(\alpha,c>0\), the notation \(\operatorname{Gamma}(\alpha,c)\)
means the probability distribution on \((0,\infty)\) with density
\begin{equation}\label{eq:gamma-distribution-convention}
 \frac{x^{\alpha-1}e^{-x/c}}{\Gamma(\alpha)c^\alpha}
 \qquad(x>0).
\end{equation}
Thus the second parameter is the scale, and unit scale means \(c=1\).

Fix \(p\ge1\), \(0\le r\le p\), and \(b\ge0\).
Applying Theorem~\ref{thm:moment-input-families} to the moments
\(b_n=b^nR_{p,r}(n)\) of a scaled Raney variable, we obtain
\begin{equation}\label{eq:raney-tree-definition}
 U_{p,r,b}(t):=tB_p(bt)^r\in\mathcal C.
\end{equation}
More precisely, if \(X\) has law \(\mu_{p,r}\) and \(G\) is an
independent gamma variable of shape two and unit scale, then
\begin{equation}\label{eq:raney-tree-shifted-measure}
 U_{p,r,b}^{(n+1)}(0)
       =(n+1)!b^nR_{p,r}(n)=\E(b\,GX)^n\qquad(n\ge0).
\end{equation}

\begin{corollary}[Fuss--Raney tree lifting]\label{thm:raney-tree-lifting}
For fixed \(p_\ell\ge1\), \(0\le r_\ell\le p_\ell\),
\(b_\ell,q_\ell\ge0\), and any initial Stieltjes moment series
\(A\), define
\begin{equation}\label{eq:raney-colored-components}
 \sum_{n\ge0}F_n(\by)\frac{t^n}{n!}
 =A(t)\prod_{\ell=1}^m
       (1-q_\ell U_{p_\ell,r_\ell,b_\ell}(t))^{-y_\ell}
\end{equation}
Then \((F_n)\) is jointly coefficientwise Hankel-totally
positive in \(\by\).
\end{corollary}

\begin{proof}
We have \(U_{p,r,b}=t\sum_{n\ge0}b_nt^n
=O_{\boldsymbol b}\) for the moments verified above.
Applying Theorem~\ref{thm:moment-input-families} to these
inputs, with zero exponential indeterminates, we obtain
positivity for \eqref{eq:raney-colored-components}. This includes
\(r=0\) or \(b=0\), when \(U_{p,r,b}=t\).
\end{proof}

For integer \(p\ge1\) and \(r\in\{0,\ldots,p\}\), we interpret
\(tB_p(bt)^r\) as a rooted tree with \(r\) distinguished child slots
at its root and \(p\) slots at every other vertex; we give each edge
weight \(b\). After labeling the vertices, we regard this series as an
exponential generating function. Using the labeled cycle construction
\cite[Section~5.1]{Stanley}, we obtain
\eqref{eq:raney-colored-components} for directed cycles decorated by
these rooted trees, with weight \(q_\ell\) at each cyclic vertex and
\(y_\ell\) at each component of color \(\ell\).

For nonintegral parameters, the construction uses the formal
series and moment representation
\eqref{eq:raney-tree-shifted-measure}.

To compute the connected coefficients explicitly, we expand the
logarithm. Applying \eqref{eq:lagrange-source} to
\(B_p(z)-1=z(1+B_p(z)-1)^p\), we obtain, for every \(s\ge0\)
and \(j\ge1\),
\begin{equation}\label{eq:raney-general-power-coefficients}
 [z^j]B_p(z)^s
 =\frac{s}{j}[u^{j-1}](1+u)^{pj+s-1}
 =R_{p,s}(j).
\end{equation}
For \(j=0\), both coefficients are one. Using
\eqref{eq:raney-general-power-coefficients} with \(s=rk\), we obtain
\begin{equation}\label{eq:raney-connected-coefficients}
 n![t^n]\bigl[-\log(1-qU_{p,r,b}(t))\bigr]
 =n!\sum_{k=1}^n\frac{q^k b^{n-k}}k
                         R_{p,rk}(n-k)\qquad(n\ge1).
\end{equation}
In \eqref{eq:raney-connected-coefficients}, we evaluate
\eqref{eq:raney-numbers} for the coefficient parameter \(rk\);
we use it only in the coefficient formula, without a moment
interpretation for the auxiliary parameter.

We recover two familiar implicit tree equations by selecting
subfamilies of \eqref{eq:raney-tree-definition}.

\begin{corollary}[Two implicit tree families]\label{cor:implicit-tree-subfamilies}
For \(b\ge0\), each of the following unique formal solutions belongs
to \(\mathcal C\):
\begin{align}
 U(t)&=t(1+b\,U(t))^p &(p\ge1),\label{eq:positive-slot-tree}\\
 V(t)&=t(1-b\,V(t))^{-a} &(a>0).\label{eq:plane-tree-family}
\end{align}
Their shifted coefficients are, respectively,
\begin{align}
 U^{(n+1)}(0)&=n!\binom{p(n+1)}n b^n,
                  \label{eq:positive-slot-moments}\\
 V^{(n+1)}(0)&=b^n\bigl(a(n+1)\bigr)_n
                  \qquad(n\ge0).\label{eq:plane-tree-moments}
\end{align}
\end{corollary}

\begin{proof}
For \eqref{eq:positive-slot-tree}, we take
\(U=tB_p(bt)^p\). From \eqref{eq:fuss-catalan-series}, we have
\(B_p(bt)=1+bU\), proving the equation. For
\eqref{eq:plane-tree-family}, take \(V=tB_{a+1}(bt)^a\).
Dividing \(B_{a+1}(bt)=1+btB_{a+1}(bt)^{a+1}\) by
\(B_{a+1}(bt)\), we obtain \(1-bV=B_{a+1}(bt)^{-1}\), and hence
the stated equation for \(V\). By
Corollary~\ref{thm:raney-tree-lifting}, both series are admissible.
Using \eqref{eq:implicit-tree-shifted-moments}, we obtain
\eqref{eq:positive-slot-moments} and \eqref{eq:plane-tree-moments}.
For \(a=1\), we obtain the usual plane-tree equation
\(V=t/(1-bV)\).
\end{proof}

Finally, we recover the original Cayley tree function as a limit of
finite-slot tree functions. Put \(U_p=U_{p,p,1/p}\). By
\eqref{eq:positive-slot-tree} and \eqref{eq:positive-slot-moments},
we have
\begin{equation}\label{eq:cayley-finite-slot-limit}
 U_p=t(1+U_p/p)^p,\qquad
 U_p^{(n)}(0)=\prod_{j=0}^{n-2}(n-j/p)
          \longrightarrow n^{n-1}\quad(n\ge1).
\end{equation}
We interpret an empty product as one. We therefore obtain coefficientwise convergence
\(U_p\to T\), where \(T=te^T\). We may also take
nonnegative finite sums and finite compositions of the series in
\eqref{eq:raney-tree-definition}, by adding their representing measures
and applying Theorem~\ref{thm:composition}. Using these series in
Corollary~\ref{thm:logarithmic}, or as substitution series in
Theorem~\ref{thm:substitution}, we obtain further colored tree families
with arbitrary initial Stieltjes moments.

\subsection{Integrated tree moment kernels}
\label{treeapp:sec:kernels}

We count sets of rooted-tree components and connected
functional-digraph components on a labeled vertex set, marking
the two component types by fixed weights \(a,b\ge0\).
The exponential coefficients of \(e^{aT+b\mathcal L}\)
are the resulting weighted counting sequence. Applying
Theorem~\ref{thm:lifting} to \(T,\mathcal L\), we obtain
its numerical Stieltjes property. We then integrate this moment
series and obtain the general component family below. We first
record the scale closure used in that construction.

\begin{lemma}[Positive combinations and argument scales]
\label{treeapp:lem:scales}
If \(K_1,\ldots,K_s\in\mathcal C\) and
\(\lambda_j,c_j\ge0\) are fixed real numbers, then
\begin{equation}\label{treeapp:eq:scale-closure}
 \sum_{j=1}^s\lambda_j K_j(c_jt)\in\mathcal C.
\end{equation}
\end{lemma}
\begin{proof}
Let \(\nu_j\) represent the once-shifted coefficients of \(K_j\).
Differentiating the scaled series, we obtain
\begin{equation}\label{treeapp:eq:scale-moments}
 \left.\frac{\dif^{n+1}}{\dif t^{n+1}}
       \sum_{j=1}^s\lambda_jK_j(c_jt)\right|_{t=0}
 =\sum_{j=1}^s\lambda_jc_j
            \int_0^\infty(c_jx)^n\,\nu_j(\dif x).
\end{equation}
For each \(j\), take the image of
\(\lambda_jc_j\nu_j\) under \(x\mapsto c_jx\), and add these
positive measures. We thereby obtain the moments in
\eqref{treeapp:eq:scale-moments}. When \(c_j=0\), we take the
zero measure, consistently with \(K_j(0)=0\). This proves
\eqref{treeapp:eq:scale-closure}.
\end{proof}

Define
\begin{equation}\label{treeapp:eq:connected-log}
 \mathcal L(t):=-\log(1-T(t)).
\end{equation}
By Proposition~\ref{prop:tree-admissible} and
Corollary~\ref{thm:logarithmic}, we have
\(T,\mathcal L\in\mathcal C\).
We use these two exponents to construct the following family of integrated kernels.

\begin{theorem}[Integrated tree kernels and component lifting]
\label{treeapp:thm:kernels}
For fixed \(a,b\ge0\), define the formal primitive
\begin{equation}\label{treeapp:eq:kernel-integrals}
 \mathcal K_{a,b}(t):=\int_0^t
          \frac{e^{aT(s)}}{(1-T(s))^b}\,\dif s
       =\int_0^{T(t)}e^{(a-1)u}(1-u)^{1-b}\,\dif u.
\end{equation}
Then \(\mathcal K_{a,b}\in\mathcal C\).
More generally, fix \(a_\ell,b_\ell,c_\ell\ge0\) for
\(1\le\ell\le m\), and let
\(A(t)\) be any initial Stieltjes moment series. Define
\begin{equation}\label{treeapp:eq:kernel-lifting}
 \sum_{n\ge0}F_n(\by)\frac{t^n}{n!}
 =A(t)\exp\!\left(\sum_{\ell=1}^m
             y_\ell \mathcal K_{a_\ell,b_\ell}(c_\ell t)\right).
\end{equation}
Then \((F_n)\) is jointly coefficientwise Hankel-totally
positive in \(\by\).
\end{theorem}

\begin{proof}
Define \(m_n(a,b)\) by
\begin{equation}\label{treeapp:eq:kernel-derivative}
 \mathcal K_{a,b}'(t)=e^{aT(t)+b\mathcal L(t)}
            =\sum_{n\ge0}m_n(a,b)\frac{t^n}{n!}.
\end{equation}
Since \(T,\mathcal L\in\mathcal C\),
applying Theorem~\ref{thm:lifting} \textup{(iii)} to these
exponents at fixed \(a,b\ge0\), we obtain the Stieltjes
sequence \((m_n(a,b))\). Applying
Theorem~\ref{thm:moment-input-families} to this sequence,
we obtain \(\mathcal K_{a,b}=I_{\boldsymbol m}\in\mathcal C\).

Making the formal change of variables \(s=ue^{-u}\), with
\(\dif s=(1-u)e^{-u}\dif u\), we obtain the second expression
in \eqref{treeapp:eq:kernel-integrals}; both primitives have
zero constant term. Applying Lemma~\ref{treeapp:lem:scales} to the argument scales
and Theorem~\ref{thm:lifting} to the resulting exponents,
we obtain the scaled component family.
\end{proof}

Several familiar connected series occur in this family:
\begin{equation}\label{treeapp:eq:kernel-specializations}
 \mathcal K_{1,0}=T-\tfrac12T^2,\qquad \mathcal K_{1,1}=T,\qquad
 \mathcal K_{1,2}=\mathcal L,\qquad \mathcal K_{1,3}=\frac{T}{1-T}.
\end{equation}

These are the series for unrooted labeled trees, rooted trees,
connected functional digraphs, and nonempty functional digraphs,
respectively; the last equals \((1-T)^{-1}-1\).

\subsection{Exponential tree composites and unrooted forests}
\label{treeapp:sec:forests}

For rooted forests, Sokal introduced a weight \(1+w\) for each
root descent, namely an edge from a root to a child with a smaller
label \cite[Section~6.7]{Sokal2023}. For that model, we express the connected series as an exponential
composite of \(T\). Ding and Zeng obtained the same coefficient
identity in their study of Ramanujan's series
\cite[Eq.~(1.3)]{DingZeng2024}. For unrooted trees, we use the
classical relation recorded by Josuat-Verg\`es
\cite[Section~2]{JV2015}. We collect these known identities in the following lemma.

\begin{lemma}[Classical rooted and unrooted tree identities]
\label{treeapp:lem:forest-sources}
Define \(B_w(t)=(e^{wT(t)}-1)/w\) for \(w\ne0\), and
\(B_0(t)=T(t)\). Then
\begin{equation}\label{treeapp:eq:abel-series}
 B_w(t)=\sum_{n\ge1}(n+w)^{n-1}\frac{t^n}{n!},\qquad
 B_w'(t)=\frac{e^{(w+1)T(t)}}{1-T(t)}.
\end{equation}
Also, the unrooted-tree series satisfies
\begin{equation}\label{treeapp:eq:unrooted-series}
 U(t):=T(t)-\tfrac12T(t)^2
      =\sum_{n\ge1}n^{n-2}\frac{t^n}{n!},\qquad U'(t)=e^{T(t)}.
\end{equation}
The coefficient identities use \(0^0=1\) when needed.
\end{lemma}

We now extend these two forest models by allowing arbitrary initial
moments and independent component colors.

\begin{corollary}[Universal forest component families]
\label{treeapp:cor:forests}
For every fixed \(w\ge-1\), we have \(B_w\in\mathcal C\),
and \(U\in\mathcal C\). Consequently, for fixed
\(w_\ell\ge-1\), \(c_\ell,d_\ell\ge0\), and any real initial
Stieltjes moment series \(A(t)\), the exponential coefficients of
\begin{equation}\label{treeapp:eq:forest-colors}
 A(t)\exp\!\left(\sum_{\ell=1}^m
         \{y_\ell B_{w_\ell}(c_\ell t)+z_\ell U(d_\ell t)\}\right)
\end{equation}
are coefficientwise Hankel-totally positive jointly in
\(\by,\boldsymbol z\).
In particular, for fixed \(w\ge-1\), the polynomials defined by
\begin{equation}\label{treeapp:eq:sokal-forest-egf}
 e^{xB_w(t)}=\sum_{n\ge0}F_n^\sharp(x,w)\frac{t^n}{n!}
\end{equation}
are coefficientwise Hankel-totally positive in \(x\).
\end{corollary}

\begin{proof}
Using the derivative identities in
Lemma~\ref{treeapp:lem:forest-sources} and the zero constant
terms, we obtain
\begin{equation}\label{treeapp:eq:forest-kernel-identification}
 B_w=\mathcal K_{w+1,1}\quad(w\ge-1),\qquad U=\mathcal K_{1,0}.
\end{equation}
This also holds at \(w=0\). By Theorem~\ref{treeapp:thm:kernels},
we derive admissibility. Applying
Lemma~\ref{treeapp:lem:scales} to the scales and
Theorem~\ref{thm:lifting} to the resulting exponents, we obtain
\eqref{treeapp:eq:forest-colors}. Taking \(A=1\), one color,
and unit scale, we have \eqref{treeapp:eq:sokal-forest-egf}.
\end{proof}

Sokal's joint assertions in \((x,w)\) and in \((x,v)\),
with \(w=v-1\), remain open
\cite[Conjectures~6.11(b) and~6.12(b)]{Sokal2023}.
The corollary fixes \(w\) and does not imply either assertion.

For an explicit check, put \(v=w+1\). From
\(B_w^{(n+1)}(0)=(n+w+1)^n\), we obtain
\begin{equation}\label{treeapp:eq:abel-minor}
 \det[F_{i+j}^\sharp(x,v-1)]_{i,j=0}^2
 =2x^3(v+1)^3+x^2(2v^3+12v^2+22v+11).
\end{equation}

For unrooted forests, using the usual labeled set construction, we obtain
\begin{equation}\label{treeapp:eq:unrooted-forest-egf}
 e^{yU(t)}=\sum_{n\ge0}\mathcal F_n(y)\frac{t^n}{n!}.
\end{equation}
Here we mark the number of tree components in \(\mathcal F_n(y)\)
by the exponent of the component variable. For example,
\begin{equation}\label{treeapp:eq:unrooted-initial}
 \mathcal F_1=y,\quad \mathcal F_2=y^2+y,\quad
 \mathcal F_3=y^3+3y^2+3y,\quad
 \mathcal F_4=y^4+6y^3+15y^2+16y.
\end{equation}
Applying Corollary~\ref{treeapp:cor:forests} to this family, we
obtain coefficientwise Hankel-total positivity.

\subsection{Shifted tree functions and Ramanujan--Shor and Greg trees}
\label{treeapp:sec:shifted}

Sokal~\cite[Theorem~1.3]{Sokal2023} proved coefficientwise
Hankel total positivity for the underlying rooted-forest
polynomials. We begin with their counted objects and equation,
verify the shifted moment condition from a tilted tree measure,
and then extend the fixed-weight model to arbitrary initial
moments and independent component intensities.

Ramanujan--Shor polynomials refine rooted labeled trees by proper
and improper edges. For an edge directed from a parent to a child, we call the edge
improper if the parent label exceeds the smallest label in the
child's subtree, and proper otherwise.
Zeng developed this refinement in connection with Ramanujan's
series \cite{Zeng1999}. From Sokal's generating function with two
edge weights \cite[Proposition~3.1(d) and Corollary~3.2]{Sokal2023},
we obtain the differential equation below. To connect this equation
with our moment criterion, we first study translations of the tree
function.

For this purpose, we use the principal real continuation of \(T\)
to \(( -\infty,1/e)\): it is the inverse of
\(u\mapsto ue^{-u}\) on \(u<1\). In the next lemma, we evaluate this analytic continuation at
nonzero arguments.

\begin{lemma}[Tilted tree moment measures]
\label{treeapp:lem:tilted-tree}
Let \(\nu\) be a probability measure representing
\(((n+1)^n)_{n\ge0}\). For every fixed real \(c<1/e\), the
Taylor series
\begin{equation}\label{treeapp:eq:shifted-tree}
 \widehat T_c(s):=T(c+s)-T(c)
\end{equation}
belongs to \(\mathcal C\), with shifted moment measure
\(e^{cx}\nu(\dif x)\).
\end{lemma}

\begin{proof}
To expand about a nonzero center, we must justify differentiating
the moment representation and show that the resulting integrals
are finite. We prove these assertions in the following claim.

\begin{claim}\label{treeapp:claim:tilted-identity}
For \(c<1/e\) and \(n\ge0\), we have
\begin{equation}\label{treeapp:eq:tilted-identity}
 T^{(n+1)}(c)=\int_0^\infty x^ne^{cx}\,\nu(\dif x)<\infty.
\end{equation}
\end{claim}
\begin{proof}
For \(0\le c<1/e\), applying Tonelli's theorem to the
nonnegative exponential expansion, we obtain
\begin{equation}\label{treeapp:eq:positive-center}
 \int_0^\infty e^{cx}\,\nu(\dif x)
 =\sum_{k\ge0}\frac{c^k}{k!}(k+1)^k=T'(c)<\infty.
\end{equation}
From the defining coefficients, we know that the series on the
right has radius \(1/e\) and is the derivative of \(T\). For any fixed \(c<1/e\), choose
\(r\in(\max\{c,0\},1/e)\). For every \(n\ge0\), there is
a finite constant \(C_{n,c,r}\) such that
\begin{equation}\label{treeapp:eq:tilted-bound}
 x^ne^{cx}\le C_{n,c,r}e^{rx}\qquad(x\ge0).
\end{equation}
By \eqref{treeapp:eq:positive-center}, we have an integrable bound.
Using this estimate on a neighborhood of \(c\), we may
differentiate under the integral any finite number of times. More
explicitly, for \(|h|<r-c\), the exponential expansion around
\(c\) is dominated in absolute value by \(e^{rx}\).
We therefore obtain a real analytic function throughout
\((-\infty,1/e)\).

On the other hand, \(\frac{\dif}{\dif u}(ue^{-u})=(1-u)e^{-u}>0\)
for \(u<1\). We therefore obtain a real analytic inverse \(T\),
with real analytic derivative, on the interval of centers under
consideration. By the identity theorem for real
analytic functions, we extend \eqref{treeapp:eq:positive-center}
to every \(c<1/e\). Differentiating this identity under the
integral, justified by \eqref{treeapp:eq:tilted-bound}, we obtain
\eqref{treeapp:eq:tilted-identity}.
\end{proof}

Since \(\widehat T_c^{(n+1)}(0)=T^{(n+1)}(c)\), we obtain
a finite positive measure representing all the shifted
coefficients by \eqref{treeapp:eq:tilted-identity}. Hence
\(\widehat T_c\in\mathcal C\), as asserted.
\end{proof}

We now state the known rooted-tree equation and its centered
tree-function form. In \cite[Eq.~(3.18)]{Sokal2023}, the two
parameters are denoted by \(y,z\). We use \(\eta,\zeta\)
to distinguish them from the component intensities.

\begin{lemma}[Sokal's rooted-tree equation]
\label{treeapp:lem:ramanujan-source}
For \(\eta,\zeta\ge0\), the rooted-tree series
\(R_{\eta,\zeta}(t)\) satisfies
\begin{equation}\label{treeapp:eq:ramanujan-ode}
 R_{\eta,\zeta}'(t)=
       \frac{e^{\zeta R_{\eta,\zeta}(t)}}{1-\eta R_{\eta,\zeta}(t)},
 \qquad R_{\eta,\zeta}(0)=0.
\end{equation}
For \(\eta,\zeta>0\), let
\begin{equation}\label{treeapp:eq:ramanujan-centers}
 a=1-\zeta/\eta<1,\qquad c=ae^{-a}<1/e,\qquad
 b=(\zeta^2/\eta)e^{-a}>0.
\end{equation}
Then
\begin{equation}\label{treeapp:eq:ramanujan-shift}
 R_{\eta,\zeta}(t)=\frac{T(c+bt)-a}{\zeta}
                         =\frac1\zeta\widehat T_c(bt).
\end{equation}
\end{lemma}

We now use these admissible series in the lifting theorem with
arbitrary initial moments and independent component colors.

\begin{corollary}[Rooted-tree component lifting]
\label{treeapp:cor:ramanujan-lifting}
For every fixed \(\eta,\zeta\ge0\), we have
\(R_{\eta,\zeta}\in\mathcal C\).
For independently chosen fixed \(\eta_\ell,\zeta_\ell,c_\ell\ge0\)
and every real initial Stieltjes moment series \(A(t)\), the
exponential coefficients of
\begin{equation}\label{treeapp:eq:ramanujan-lifting}
 A(t)\exp\!\left(\sum_{\ell=1}^m
           y_\ell R_{\eta_\ell,\zeta_\ell}(c_\ell t)\right)
\end{equation}
are coefficientwise Hankel-totally positive jointly in \(\by\).
\end{corollary}

\begin{proof}
For \(\eta,\zeta>0\), combine
\eqref{treeapp:eq:ramanujan-shift} with
Lemmas~\ref{treeapp:lem:tilted-tree} and
\ref{treeapp:lem:scales}. On the boundary, integration gives
\begin{equation}\label{treeapp:eq:ramanujan-boundaries}
 R_{0,\zeta}(t)=-\frac{\log(1-\zeta t)}{\zeta}\quad(\zeta>0),
 \qquad R_{\eta,0}(t)=\frac{1-\sqrt{1-2\eta t}}{\eta}\quad(\eta>0).
\end{equation}
The shifted coefficients are \(n!\zeta^n\) and
\((1/2)_n(2\eta)^n\), the moments of
\(\operatorname{Gamma}(1,\zeta)\) and
\(\operatorname{Gamma}(1/2,2\eta)\), respectively.
At \(\eta=\zeta=0\), \(R_{0,0}=t\) has shifted measure
\(\delta_0\). Hence every stated exponent is admissible.
Applying scale closure to these exponents and
Theorem~\ref{thm:lifting} to the initial series and scaled
exponents, we obtain the conclusion.
\end{proof}

In a rooted Greg tree, we allow additional unlabeled branching
vertices alongside the vertices labeled \(1,\ldots,n\). Every unlabeled
vertex has degree at least three, except that an unlabeled root may
have degree two. We regard the trees as non-plane and identify
them by isomorphisms preserving the root and all labels.
Josuat-Verg\`es \cite[Definitions~3.1--3.2 and Theorem~4.2]{JV2015}
shows that their generating polynomial, with weight \(q\) for
each unlabeled vertex, is
\(\mathfrak G_n(q)=R_{1+q,1}^{(n)}(0)\).
For the first few polynomials, we have
\begin{equation}\label{treeapp:eq:greg-initial}
 \mathfrak G_1=1,\quad \mathfrak G_2=q+2,\quad
 \mathfrak G_3=3q^2+10q+9,\quad
 \mathfrak G_4=15q^3+70q^2+113q+64.
\end{equation}
Thus, from \eqref{treeapp:eq:ramanujan-lifting}, we also obtain
Greg component families for independently fixed \(q_\ell\ge0\).
We retain the component intensities as indeterminates and use the
established scalar edge-weight positivity as an input.

\subsection{Ward moment inputs and phylogenetic component families}
\label{treeapp:sec:phylogenetic}

We count non-plane rooted trees with \(n+1\) labeled leaves
and unlabeled internal vertices, each having at least two
children. Giving each internal vertex weight \(u\), we obtain
the Ward polynomial \(W_n(u)\). Elvey Price and
Sokal~\cite[Section~1, item~4, pp.~11--12, and Appendix~A]{EPS2020}
related these trees to augmented perfect matchings and proved
coefficientwise Hankel total positivity in \(u\); see also
\cite[Section~6.5]{Zhu2024}.
We recall the generating identity from
\cite[Eqs.~(1.17)--(1.19)]{EPS2020}.

\begin{lemma}[Ward counting series and its moment input]
\label{treeapp:lem:ward-source}
We have
\begin{equation}\label{treeapp:eq:phylo-equation}
 H_u(t):=\sum_{n\ge0}W_n(u)\frac{t^{n+1}}{(n+1)!}
 =t+u\bigl(e^{H_u(t)}-1-H_u(t)\bigr).
\end{equation}
Thus \(H_u^{(n+1)}(0)=W_n(u)\). For each fixed \(u\ge0\),
\((W_n(u))_{n\ge0}\) is a Stieltjes moment sequence.
\end{lemma}

Applying the real moment criterion to the known Ward sequence
at a fixed \(u\), we choose a measure \(\nu_u\) with moments
\(W_n(u)\). We therefore have the shifted moment representation
\(H_u^{(n+1)}(0)=\int_0^\infty x^n\nu_u(\dif x)\).
Applying the primitive moment construction to this input,
we obtain \(H_u\in\mathcal C\). We now retain independent
component intensities, logarithmic powers, an admissible
substitution, and arbitrary initial moments in one general result.

\begin{theorem}[Universal lifting of Ward moment inputs]
\label{treeapp:thm:phylo-lifting}
Fix \(u_\ell,c_\ell,\theta_\ell\ge0\), put
\(H_\ell(t)=H_{u_\ell}(c_\ell t)\), and choose
\(S\in\mathcal C\). For any initial Stieltjes moment series
\(A\), the exponential coefficients of
\begin{equation}\label{treeapp:eq:ward-universal-family}
 A(S(t))\exp\!\left(\sum_{\ell=1}^m y_\ell H_\ell(t)\right)
       \prod_{\ell=1}^m(1-\theta_\ell H_\ell(t))^{-z_\ell}
\end{equation}
are jointly coefficientwise Hankel-totally positive in
\(\by,\boldsymbol z\).
\end{theorem}

\begin{proof}
By Lemma~\ref{treeapp:lem:scales}, we derive the shifted Stieltjes inputs of
\(H_\ell\). Applying Theorem~\ref{thm:moment-input-families}
to these inputs, with substitution \(S\) and logarithmic
parameters \(\theta_\ell\), we obtain the assertion.
\end{proof}

Taking \(S=t\) and \(\boldsymbol z=\boldsymbol0\), we have
the exponential component family
\begin{equation}\label{treeapp:eq:phylo-lifting}
 A(t)\exp\!\left(\sum_{\ell=1}^m
                         y_\ell H_{u_\ell}(c_\ell t)\right)
\end{equation}
At \(A=1\), we count forests of independently weighted
phylogenetic components; we mark their colors by \(\by\).
In the one-color case at unit scale, we recover
\(W_n(u)=[y]F_{n+1}(y;u)\). Our new coefficientwise
conclusion concerns the component intensities; we keep the
internal-vertex weights fixed. We therefore distinguish this
universal extension from the cited positivity in the structural
variable \(u\).

We also identify the component series explicitly. For fixed
\(u\ge0\), put
\begin{equation}\label{treeapp:eq:phylo-parameters}
 q=\frac{u}{1+u},\qquad d=\frac1{1+u},\qquad c=qe^{-q}<1/e.
\end{equation}
We claim that
\begin{equation}\label{treeapp:eq:phylo-composition}
 H_u(t)=dt-q+T(ce^{dt})
       =dt+\widehat T_c\bigl(c(e^{dt}-1)\bigr).
\end{equation}
To verify the formula, set \(Z(t)=T(ce^{dt})\) and
\(J(t)=dt-q+Z(t)\), using the Taylor expansion at \(c\).
Since \(T(c)=q\), we have \(J(0)=0\). By \(T=te^T\),
we have
\begin{equation}\label{treeapp:eq:phylo-verification}
 Z(t)=ce^{dt+Z(t)}=qe^{J(t)},\qquad
 (1+u)J(t)=t-u+ue^{J(t)}.
\end{equation}
Rearranging the second identity, we deduce the defining equation
for \(H_u\). Comparing coefficients successively, we obtain
\(J=H_u\), because \(e^J-1-J\) starts in degree two.
At \(u=0\), we obtain \(H_0=t\), with shifted measure
\(\delta_0\).

For example, at unit internal-vertex weight we have
\begin{equation}\label{treeapp:eq:phylo-unit}
 H_1(t)=\frac{t-1}{2}
       +T\!\left(\tfrac12e^{(t-1)/2}\right)
       =t+\frac{t^2}{2!}+4\frac{t^3}{3!}+26\frac{t^4}{4!}+\cdots.
\end{equation}
Combining the unrooted-tree and Ward inputs, we obtain
\begin{equation}\label{treeapp:eq:mixed-forest}
 e^{yU(t)+zH_1(t)}=\sum_{n\ge0}M_n(y,z)\frac{t^n}{n!},
\end{equation}
with joint coefficientwise Hankel total positivity in \(y,z\).
The first leading determinants are
\begin{equation}\label{treeapp:eq:mixed-minors}
 \det[M_{i+j}]_{i,j=0}^1=y+z,\qquad
 \det[M_{i+j}]_{i,j=0}^2=2(y+z)^3+7y^2+18yz+10z^2.
\end{equation}
We may include arbitrary initial moments in this mixed family by
applying Theorem~\ref{thm:lifting} to \(U\) and \(H_1\).

\subsection{Gamma damping of tree coefficients and an excluded composite}
\label{treeapp:sec:damping}

The family \(T_p(t)=\sum_{n\ge1}n^{n-p}t^n/n!\) includes
rooted and unrooted tree series at \(p=1,2\)
\cite[Section~2]{JV2015}. Using a gamma damping factor,
we verify the shifted moment condition for every real \(p\ge1\)
and obtain a general component family by lifting.

\begin{proposition}[Gamma-damped tree exponents]
\label{treeapp:prop:damped}
For every real \(p\ge1\), the formal power series
\begin{equation}\label{treeapp:eq:damped-definition}
 T_p(t):=\sum_{n\ge1}n^{n-p}\frac{t^n}{n!}
\end{equation}
belongs to \(\mathcal C\). For fixed \(p_\ell\ge1\) and
\(c_\ell\ge0\), the exponential coefficients of
\begin{equation}\label{treeapp:eq:damped-lifting}
 A(t)\exp\!\left(\sum_{\ell=1}^m y_\ell T_{p_\ell}(c_\ell t)\right)
\end{equation}
are coefficientwise Hankel-totally positive jointly in \(\by\)
for every real initial Stieltjes moment series \(A(t)\).
\end{proposition}

\begin{proof}
Let \(X\) have moments \((n+1)^n\). For \(p>1\), take an
independent \(G\sim\operatorname{Gamma}(p-1,1)\) and put
\(Q=e^{-G}\). Applying Euler's integral to the gamma density,
we obtain
\begin{equation}\label{treeapp:eq:damping-moments}
 \E Q^n=\frac1{\Gamma(p-1)}
       \int_0^\infty s^{p-2}e^{-(n+1)s}\,\dif s
       =(n+1)^{-(p-1)}.
\end{equation}
Using independence, we then have
\begin{equation}\label{treeapp:eq:damped-shifted-moments}
 \E(XQ)^n=(n+1)^n(n+1)^{-(p-1)}
          =(n+1)^{n+1-p}=T_p^{(n+1)}(0).
\end{equation}
All moments are finite since \(0<Q\le1\). For \(p=1\),
take \(Q=1\). Thus \(T_p\in\mathcal C\). Applying scale closure
and then Theorem~\ref{thm:lifting} to these exponents yields the component assertion.
\end{proof}

Nonnegative coefficients alone do not imply admissibility.

For \(K(t)=T(t)^2\), the first two cumulants are
\(\kappa_1=0\) and \(\kappa_2=2\). A positive measure with
zeroth moment zero is the zero measure, so its first moment cannot
equal two. Hence \(T^2\notin\mathcal C\).
Equivalently, defining \(P_n(y)\) by
\begin{equation}\label{treeapp:eq:excluded-egf}
 e^{yT(t)^2}=\sum_{n\ge0}P_n(y)\frac{t^n}{n!},
\end{equation}
we have
\begin{equation}\label{treeapp:eq:excluded-minor}
 P_1(y)P_3(y)-P_2(y)^2=-4y^2.
\end{equation}
The shifted moment condition is therefore essential.

\section{Polynomial moment inputs and further lifted families}\label{sec:applications}

We will obtain general coefficientwise Hankel-positive families
from the shifted moment conditions of concrete counting sequences.
For the Schl\"afli-Gessel-Seo polynomials polynomials, we first prove that condition and obtain
fixed-shape positivity. For Bessel, Lah, and Laguerre polynomials,
we recover known Hankel positivity by a different proof and allow
arbitrary initial moments, independent colors, and independently
chosen fixed shape parameters. We specify these extensions below
and attribute the earlier polynomial results at each specialization.
We also obtain lognormal jump families from their classical moments;
we use formal series even when exponential moments are infinite.

\subsection{Schl\"afli-Gessel-Seo polynomials}

For \(n\ge1\), let
\begin{equation}\label{eq:sgs-polynomial-definition}
 P_0(x;a,b)=1,\qquad
 P_n(x;a,b)=x\prod_{h=1}^{n-1}\bigl(x+ha+(n-h)b\bigr).
\end{equation}
We first recall the combinatorial interpretation
\cite[Section~6.1, Eqs.~(6.10)--(6.11)]{Sokal2023}.
For a rooted labeled forest \(F\), call a vertex proper if
every other vertex in its descendant subtree has a larger label.
Let \(k(F)\) and \(p(F)\) be its numbers of components and
proper vertices. We obtain \(P_n(x;a,b)\) by giving \(F\)
weight \(x^{k(F)}a^{n-p(F)}b^{p(F)-k(F)}\).
Sokal conjectured joint coefficientwise Hankel total positivity
in \(x,a,b\). We verify the shifted moment condition below
for every fixed \(a,b\ge0\), and then obtain a new universal
coefficientwise result in the component intensities.

Writing \(P_{n,k}(a,b)=[x^k]P_n(x;a,b)\), we recall the
Hankel assertions of the original conjecture before giving the
conclusion obtained from real cumulant lifting.

\begin{conjecture}[Sokal, Conjecture~6.1(b),(c) in~\cite{Sokal2023}]
\label{conj:sgs-joint}
The sequence \((P_n(x;a,b))_{n\ge0}\) is coefficientwise
Hankel-totally positive jointly in \(x,a,b\). Moreover, the sequence
\((P_{n+1,1}(a,b))_{n\ge0}\) is coefficientwise Hankel-totally
positive jointly in \(a,b\).
\end{conjecture}

For fixed \(a,b\ge0\), we verify the shifted moment
condition and apply Theorem~\ref{thm:lifting}. The source
identities below follow from
\cite[Section~6.1, Eqs.~(6.1)--(6.4)]{Sokal2023}; extracting
the coefficient of \(x\), we obtain the final identity.

\begin{lemma}[Sokal's generating-function identities]
\label{lem:sgs-source}
Let \(G_{a,b}(t)\in t\R[[t]]\) be the compositional inverse of
\begin{equation}\label{eq:sgs-inverse-definition}
 g_{a,b}(u)=
 \begin{cases}
 \dfrac{e^{-bu}-e^{-au}}{a-b},&a\ne b,\\[5pt]
 ue^{-au},&a=b.
 \end{cases}
\end{equation}
Then \(g_{a,b}(0)=0\), \(g_{a,b}'(0)=1\), and
\begin{equation}\label{eq:sgs-known-egf}
 \sum_{n\ge0}P_n(x;a,b)\frac{t^n}{n!}=e^{xG_{a,b}(t)}.
\end{equation}
Consequently, if \(c_n(a,b):=P_{n+1,1}(a,b)\), then
\begin{equation}\label{eq:sgs-shifted-coefficients}
 c_n(a,b)=G_{a,b}^{(n+1)}(0)
 =\prod_{h=1}^{n}\bigl(ha+(n+1-h)b\bigr)\qquad(n\ge0).
\end{equation}
\end{lemma}

To construct a positive moment measure, we use the classical
Bernstein--Widder theorem \cite[Theorem~1.4]{SSV2012}.
We call a function \(h:(0,\infty)\to\R\)
\emph{completely monotone} if it is infinitely differentiable and
\((-1)^nh^{(n)}(s)\ge0\) for every \(n\ge0\) and \(s>0\).

\begin{lemma}[Bernstein--Widder representation]
\label{lem:bernstein-widder-source}
A function \(h\) is completely monotone if and only if
\(h(s)=\int_0^\infty e^{-sz}\,\nu(\dif z)\) for a positive
measure \(\nu\) on \([0,\infty)\), with finite integral for
every \(s>0\). If \(h(0+)<\infty\), then
\(\nu([0,\infty))=h(0+)\).
\end{lemma}

Using this representation and an explicit calculation of
derivative signs, we construct the shifted moment measure
required for lifting.

\begin{lemma}[Shifted SGS moments at fixed shape parameters]
\label{lem:sgs-shifted-moments}
For fixed \(a,b\ge0\), there is a probability measure
\(\nu_{a,b}\) on \([0,\infty)\) such that
\begin{equation}\label{eq:sgs-shifted-moment-measure}
 c_n(a,b)=\int_0^\infty z^n\,\nu_{a,b}(\dif z)
 \qquad(n\ge0).
\end{equation}
In particular, \(G_{a,b}\in\mathcal C\).
\end{lemma}

\begin{proof}
\emph{Step 1: A global inverse.}
By symmetry, assume \(a\ge b\ge0\). For \(u\ge0\), put
\begin{equation}\label{eq:sgs-eta-definition}
 \eta(u):=-g_{a,b}(-u)=
 \begin{cases}
 \dfrac{e^{au}-e^{bu}}{a-b},&a>b,\\[5pt]
 ue^{au},&a=b.
 \end{cases}
\end{equation}
Since \(\eta(0)=0\), \(\eta'(u)>0\), and
\(\eta(u)\to\infty\), it has a smooth inverse
\(\xi:[0,\infty)\to[0,\infty)\). The inverse is analytic
near zero because \(\eta'(0)=1\), and
\begin{equation}\label{eq:sgs-local-global-inverse}
 \eta(\xi(s))=s,\quad \xi(0)=0,\quad
 G_{a,b}(-s)=-\xi(s).
\end{equation}
The last identity is used only near zero.

\smallskip\noindent
\emph{Step 2: Complete monotonicity of the reciprocal derivative.}
Put \(R=1/\eta'\). If \(a>b\), let \(d=a-b>0\). Then
\begin{equation}\label{eq:sgs-reciprocal-expansion}
 R(u)=\frac{d}{ae^{au}-be^{bu}}
 =\frac da\,\frac{e^{-au}}{1-(b/a)e^{-du}}
 =\frac da\sum_{k\ge0}(b/a)^k e^{-(a+kd)u}.
\end{equation}
For \(b=0\), only the \(k=0\) term remains. For \(b>0\),
every differentiated series converges uniformly on \([0,\infty)\),
since \(\sum_{k\ge0}(b/a)^k(a+kd)^n<\infty\).
Differentiating term by term, we obtain \((-1)^nR^{(n)}\ge0\).
If \(a=b=c\),
\begin{equation}\label{eq:sgs-reciprocal-diagonal}
 R(u)=\frac{e^{-cu}}{1+cu}.
\end{equation}
Both factors are completely monotone, including at \(c=0\),
so their product is also completely monotone by Leibniz's rule.

\smallskip\noindent
\emph{Step 3: Complete monotonicity of the inverse derivative.}
Differentiating the inverse relation, we obtain
\(\xi'=R\circ\xi\). Define
\begin{equation}\label{eq:sgs-sign-recursion}
 H_0(u)=R(u),\qquad H_{n+1}(u)=-R(u)H_n'(u).
\end{equation}
Complete monotonicity is preserved by \(h\mapsto-h'\) and
by products, since
\begin{equation}\label{eq:sgs-cm-product}
 (-1)^n(h_1h_2)^{(n)}(u)
 =\sum_{k=0}^n\binom nk
       (-1)^kh_1^{(k)}(u)(-1)^{n-k}h_2^{(n-k)}(u)\ge0.
\end{equation}
Inductively, every \(H_n\) is completely monotone and
nonnegative. Using the chain rule, we obtain
\begin{equation}\label{eq:sgs-inverse-derivative-sign}
 (-1)^n\frac{\dif^n}{\dif s^n}\xi'(s)
 =H_n(\xi(s))\ge0\qquad(n\ge0,s>0).
\end{equation}
Differentiating and changing the sign, we obtain
\(-H_n'(\xi(s))R(\xi(s))=H_{n+1}(\xi(s))\).
Thus \(\xi'\) is completely monotone.

\smallskip\noindent
\emph{Step 4: Shifted moment representation.}
Since \(\xi'(0)=1\), applying
Lemma~\ref{lem:bernstein-widder-source} to \(\xi'\),
we obtain a probability measure \(\nu_{a,b}\) with
\begin{equation}\label{eq:sgs-derivative-laplace}
 \xi'(s)=\int_0^\infty e^{-sz}\,\nu_{a,b}(\dif z)
 \qquad(s>0).
\end{equation}
For fixed \(s>0\), boundedness of \(z^ne^{-sz/2}\) and
finiteness of the measure justify differentiation:
\begin{equation}\label{eq:sgs-laplace-differentiation}
 (-1)^n(\xi')^{(n)}(s)
 =\int_0^\infty z^ne^{-sz}\,\nu_{a,b}(\dif z).
\end{equation}
Near zero, \(\xi'(s)=G_{a,b}'(-s)\). As \(s\downarrow0\),
the left side tends to \(G_{a,b}^{(n+1)}(0)=c_n(a,b)\).
Applying monotone convergence to the integrals on the right,
we obtain \eqref{eq:sgs-shifted-moment-measure}, including
finiteness of every moment.
\end{proof}

Having verified the moment condition, we apply the lifting theorem
to obtain the following general result. We allow different fixed
shape parameters in each color and an arbitrary initial moment
series; the original one-color family is a specialization.

\begin{theorem}[Universal lifting for fixed SGS shape parameters]
\label{thm:sgs-fixed-lifting}
Fix \(a_\ell,b_\ell\ge0\), and let \(A(t)\) be any initial
Stieltjes moment series. Define
\begin{equation}\label{eq:sgs-multicolor-lifting}
 \sum_{n\ge0}F_n(\by)\frac{t^n}{n!}
 =A(t)\exp\!\left(\sum_{\ell=1}^m
                         y_\ell G_{a_\ell,b_\ell}(t)\right).
\end{equation}
Then \((F_n)\) is coefficientwise Hankel-totally
positive in \(\by\). In particular, for fixed \(a,b\ge0\),
\((P_n(x;a,b))\) is coefficientwise Hankel-totally positive
in \(x\), and is Stieltjes at every numerical \(x\ge0\).
\end{theorem}

\begin{proof}
Applying Lemma~\ref{lem:sgs-shifted-moments} to each fixed
pair \((a_\ell,b_\ell)\), we obtain admissible exponents
\(G_{a_\ell,b_\ell}\). Applying
Theorem~\ref{thm:lifting} to these exponents and \(A\), we
obtain the joint coefficientwise assertion. Taking \(A=1\),
one color, and \(y_1=x\), we obtain the assertion for \(P_n\)
from \eqref{eq:sgs-known-egf}; taking nonnegative numerical
\(x\), we obtain its Stieltjes property.
\end{proof}

We have obtained a fixed-shape Hankel consequence of
Conjecture~\ref{conj:sgs-joint}. The joint assertion in all
three variables \(x,a,b\) remains open; our new generality
concerns independent component intensities and initial moments.

\subsection{Bessel and Lah polynomials from gamma moment inputs}
\label{sec:gamma-families}

We begin with two counting models. We count forests of non-plane
full binary trees with leaves labeled by \([n]\), giving each
tree component weight \(x\). A component with \(j\) leaves has
\((2j-3)!!\) realizations, with the one-leaf tree included.
Using the labeled set construction, we obtain the reverse Bessel
polynomials \(B_n(x)\). We also count partitions of \([n]\)
into linearly ordered blocks, giving each block weight \(x\);
we obtain the Lah polynomials \(L_n(x)\). Set \(B_0=L_0=1\).
\begin{align}
 \sum_{n\ge0}B_n(x)\frac{t^n}{n!}
 &=\exp\!\left(x(1-\sqrt{1-2t})\right),
 & B_n(x)&=\sum_{k=1}^{n}
 \frac{(2n-k-1)!\,x^k}{2^{n-k}(n-k)!(k-1)!},\label{eq:bessel-moments}\\
 \sum_{n\ge0}L_n(x)\frac{t^n}{n!}
 &=\exp\!\left(\frac{xt}{1-t}\right),
 & L_n(x)&=\sum_{k=1}^{n}\frac{n!}{k!}
                         \binom{n-1}{k-1}x^k.\label{eq:lah-moments}
\end{align}
The finite formulas hold for \(n\ge1\).
In \cite[Propositions~4.4 and~4.7]{Zhu2021}, we proved
coefficientwise Hankel total positivity of these polynomial
sequences. We will obtain both conclusions below from gamma moment
inputs and a universal lifting result, with a different proof.

For the binary-tree exponent \(1-\sqrt{1-2t}\), the shifted
cumulants are \((2n-1)!!=2^n(1/2)_n\). For the ordered-block
exponent \(t/(1-t)\), they are \((n+1)!=(2)_n\).
We represent these sequences by the gamma laws
\(\operatorname{Gamma}(1/2,2)\) and
\(\operatorname{Gamma}(2,1)\), respectively.
We now use the whole gamma moment family to obtain a common
extension of the two models.

For \(p\in\R\) and \(c\ge0\), define
\begin{equation}\label{eq:gamma-exponent}
 K_{p,c}(t)=\sum_{n\ge1}(p)_{n-1}c^{n-1}\frac{t^n}{n!},
 \qquad K_{p,c}'(t)=(1-ct)^{-p}.
\end{equation}
Here \(K_{p,0}=K_{0,c}=t\). For \(c>0\), we have
\begin{equation}\label{eq:display-042}
 K_{p,c}(t)=
 \begin{cases}
 \displaystyle\frac{1-(1-c\,t)^{1-p}}{c(1-p)},&p\ne1,\\[5pt]
 \displaystyle-\frac{\log(1-c\,t)}c,&p=1.
 \end{cases}
\end{equation}
For \(p,c>0\), take the gamma measure
\begin{equation}\label{eq:display-043}
 \nu_{p,c}(\dif x)
 =\frac{x^{p-1}e^{-x/c}}{\Gamma(p)c^p}\,\dif x
 \qquad(x>0)
\end{equation}
Applying Euler's integral and the gamma-ratio formula
\cite[Eqs.~5.2.1 and~5.2.5]{DLMF-Gamma} to this measure,
we obtain \(\int_0^\infty x^n\nu_{p,c}(\dif x)=c^n(p)_n\).
Thus we have verified the once-shifted moment condition for
\(K_{p,c}\) before applying the lifting theorem.

\begin{corollary}[Universal gamma component families]
\label{thm:gamma-exponents}
Fix \(p_\ell,c_\ell\ge0\). For any initial Stieltjes moment
series \(A\), define
\begin{equation}\label{eq:gamma-family-egf}
 \sum_{n\ge0}P_n(\by)\frac{t^n}{n!}
 =A(t)\exp\!\left(\sum_{\ell=1}^{m}
                             y_\ell K_{p_\ell,c_\ell}(t)\right)
\end{equation}
Then \((P_n)\) is coefficientwise Hankel-totally
positive in \(\by\). For \(c>0\), the condition \(p\ge0\)
is necessary for \(K_{p,c}\in\mathcal C\).
\end{corollary}

\begin{proof}
Using the gamma measures above, we have
\(K_{p,c}=I_{\boldsymbol b}\), with \(b_n=c^n(p)_n\).
At \(p=0\) or \(c=0\), use the shifted measure \(\delta_0\).
Applying Theorem~\ref{thm:lifting} to these exponents and \(A(t)\),
we obtain the joint coefficientwise assertion.
For \(p<0<c\), we have second cumulant \(pc<0\), so we
cannot have \(K_{p,c}\in\mathcal C\).
\end{proof}

We allow independent shapes and scales, arbitrary initial moments,
and independent component intensities in this result. Taking
\((p,c)=(1/2,2)\) or \((2,1)\), we recover the two classical
polynomial families. We may also retain a common scale as an
indeterminate.

\begin{corollary}[A common scale as an indeterminate]
\label{cor:common-scale}
Fix \(p_1,\ldots,p_m\ge0\), and let \(A\) be any initial
Stieltjes moment series. The exponential coefficients of
\begin{equation}\label{eq:common-scale-egf}
 A(ct)\exp\!\left(\sum_{\ell=1}^{m}y_\ell K_{p_\ell,c}(t)\right)
\end{equation}
are jointly coefficientwise Hankel-totally positive in
\(c,\by\).
\end{corollary}

\begin{proof}
Let \(R_n(\by)\) be the coefficients at unit scale.
For \(c>0\), using \(K_{p,c}(t)=c^{-1}K_{p,1}(ct)\),
we obtain coefficients \(c^nR_n(\by/c)\). For increasing
index lists \(I,J\), let \(S=\sum_a i_a+\sum_b j_b\).
Factoring the row and column scales, we obtain
\begin{equation}\label{eq:display-045}
 c^S\det[R_{i_a+j_b}(\by/c)]_{a,b=1}^{r}.
\end{equation}
Applying Corollary~\ref{thm:gamma-exponents} to the unit-scale
family, we obtain a determinant with nonnegative coefficients.
Since its degree is at most \(S\), we multiply each monomial
\(\by^{\boldsymbol\alpha}\) by a nonnegative power
\(c^{S-|\boldsymbol\alpha|}\). Thus we obtain a polynomial
identity valid also at \(c=0\), with the initial factor
\(A(ct)\) specified in the statement.
\end{proof}

\begin{corollary}[Recovery of the Bessel and Lah results]
\label{cor:bessel-lah}
The sequences \((B_n(x))\) and \((L_n(x))\) defined in
\eqref{eq:bessel-moments}--\eqref{eq:lah-moments} are
coefficientwise Hankel-totally positive in \(x\).
\end{corollary}

\begin{proof}
Applying Corollary~\ref{thm:gamma-exponents} to one-component
families with \(A(t)=1\) and indeterminate \(x\), we obtain the two assertions
from \(K_{1/2,2}=1-\sqrt{1-2t}\) and
\(K_{2,1}=t/(1-t)\), respectively.
\end{proof}

For numerical \(x>0\), we may also represent \(B_n(x)\) by
an inverse Gaussian law with mean \(x\) and shape \(x^2\)
\cite{GinerSmyth2016}. We represent \(L_n(x)\) by a
compound-Poisson sum with rate \(x\) and unit exponential jumps.
More generally, for \(p>1,c>0\), we write
\begin{equation}\label{eq:display-044}
 K_{p,c}(t)=
 \frac{(1-ct)^{-(p-1)}-1}{c(p-1)}.
\end{equation}
We represent \(e^{yK_{p,c}}\) by gamma jumps of shape \(p-1\)
and scale \(c\), with Poisson rate \(y/[c(p-1)]\).
At \(p=1\), we use a gamma law of shape \(y/c\), including
its zero-shape degenerate case.

For generalized Lah polynomials, their Hankel total-positivity result was proved
in \cite[Proposition~4.7]{Zhu2021} and
\cite[Section~6.8, Eq.~(6.174) and Proposition~6.30]{Zhu2024}.
For fixed \(a,b,c>0\) and independent \(d,q,\lambda\), put
\begin{equation}\label{powerapp:eq:lah-real}
 f(t)=\frac cb\bigl((1-abt)^{-1/a}-1\bigr),\qquad
 A(t)(1-abt)^{-d}\exp((q+\lambda)f(t)).
\end{equation}
We represent the shifted cumulants of \(f\), namely
\(c(1+1/a)_n(ab)^n\), by \(c\) times the gamma law of shape
\(1+1/a\) and scale \(ab\). Applying
Corollary~\ref{thm:gamma-exponents} to \(-\log(1-abt)\),
we also obtain its admissibility. Applying
Theorem~\ref{thm:lifting} to these two exponents, and replacing
the intensity of \(f\) by \(q+\lambda\), we obtain joint
coefficientwise Hankel total positivity in \(d,q,\lambda\).
We allow continuous \(a,b,c>0\) and arbitrary initial moments;
we do not require the restrictions \(a,ad\in\N\) and
\(ad\le a+1\) used in the cited production-matrix result.

\subsection{Univariate Laguerre polynomials and universal component families}
\label{lagapp:sec:laguerre}

We count assemblies of permutation cycles and linearly ordered
blocks on \([n]\), with weights \(\beta\) per cycle and \(x\)
per ordered block. Using the classical exponential formula, we
obtain unsigned Laguerre polynomials with the normalization
\begin{equation}\label{lagapp:eq:normalization}
 \Lambda_n(x;\beta):=n!L_n^{(\beta-1)}(-x).
\end{equation}
For \(\beta>0\), \(L_n^{(\beta-1)}(z)\) is orthogonal on
\([0,\infty)\) with weight \(z^{\beta-1}e^{-z}\).
In \cite[Proposition~4.14]{Zhu2021} and
\cite[Proposition~6.10]{Zhu2024}, scalar coefficientwise Hankel total positivity was provedd the joint parameter
case using \cite[Theorem~5.2]{Zhu2024}; see also the
rook-polynomial interpretation in the latter work.
We retain the Laguerre parameter as a component indeterminate and
obtain the general independently colored extension below.
The normalization, finite formula and generating function are
classical~\cite[Table~18.3.1 and Eqs.~18.5.12,
18.12.13]{DLMF-Laguerre}.
At \(\beta=0\), we interpret the finite expression by polynomial
continuation. We first record these coefficient identities.

\begin{lemma}[Classical unsigned Laguerre formulas]
\label{lagapp:lem:classical}
For \(n\ge0\), we have
\begin{equation}\label{lagapp:eq:finite-formula}
 \Lambda_n(x;\beta)
   =\sum_{k=0}^n\binom nk(\beta+k)_{n-k}x^k,
 \qquad
 \sum_{n\ge0}\Lambda_n(x;\beta)\frac{t^n}{n!}
   =(1-t)^{-\beta}\exp\!\left(\frac{xt}{1-t}\right).
\end{equation}
For every fixed \(\beta\ge0\), the sequence
\((\Lambda_n(x;\beta))_{n\ge0}\) is coefficientwise
Hankel-totally positive in \(x\), as in
\cite[Proposition~6.10(ii)]{Zhu2024}.
\end{lemma}

We first verify the two shifted gamma moment conditions.
Applying lifting to these inputs, we obtain arbitrary initial
moments, independent scales, and separate intensities for cycles
and ordered blocks. Taking one color, unit scale, and \(A=1\),
we recover the known Laguerre coefficientwise positivity by
this different proof.

For fixed \(c\ge0\), put
\begin{equation}\label{lagapp:eq:component-exponents}
 U_c(t):=-\log(1-ct),\qquad V_c(t):=\frac{t}{1-ct}.
\end{equation}
Their shifted coefficients are
\begin{equation}\label{lagapp:eq:shifted-coefficients}
 U_c^{(n+1)}(0)=n!c^{n+1},\qquad
 V_c^{(n+1)}(0)=(n+1)!c^n\qquad(n\ge0).
\end{equation}
For \(c>0\), applying Euler's integral to these shifted coefficients, we obtain
\begin{equation}\label{lagapp:eq:shifted-measures}
 n!c^{n+1}=\int_0^\infty s^ne^{-s/c}\,\dif s,\qquad
 (n+1)!c^n=\int_0^\infty s^n\frac{s e^{-s/c}}{c^2}\,\dif s.
\end{equation}
At \(c=0\), use the zero measure for \(U_0=0\) and
\(\delta_0\) for \(V_0=t\). Thus \(U_c,V_c\in\mathcal C\).

\begin{theorem}[Universal Laguerre component lifting]
\label{lagapp:thm:universal}
Fix \(c_1,\ldots,c_m\ge0\), and let
\(A(t)\) be any initial Stieltjes moment series. Introduce independent indeterminates
\(\boldsymbol\beta=(\beta_1,\ldots,\beta_m)\) and
\(\bx=(x_1,\ldots,x_m)\), and define
\begin{equation}\label{lagapp:eq:colored-family}
 \sum_{n\ge0}F_n(\boldsymbol\beta,\bx)\frac{t^n}{n!}
 =A(t)\prod_{\ell=1}^m(1-c_\ell t)^{-\beta_\ell}
       \exp\!\left(\sum_{\ell=1}^m
                              \frac{x_\ell t}{1-c_\ell t}\right).
\end{equation}
Then \((F_n)\) is coefficientwise Hankel-totally
positive jointly in \(\boldsymbol\beta,\bx\).
In particular, for \(m=1\) and \(c_1=1\), we obtain
\begin{equation}\label{lagapp:eq:initial-convolution}
 F_n(\beta,x)=\sum_{k=0}^n\binom nk a_{n-k}\Lambda_k(x;\beta),
\end{equation}
which has the same joint positivity in \(\beta,x\).
\end{theorem}

\begin{proof}
Applying Theorem~\ref{thm:lifting} to
\(U_{c_\ell},V_{c_\ell}\), with  indeterminates
\(\beta_\ell,x_\ell\) and initial series \(A\),
we obtain joint coefficientwise positivity for
\eqref{lagapp:eq:colored-family}. At \(m=1,c_1=1\),
multiplying \eqref{lagapp:eq:finite-formula} by \(A(t)\)
and extracting \(n![t^n]\), we obtain
\eqref{lagapp:eq:initial-convolution}.
\end{proof}

For numerical \(x,\beta\ge0\), we represent the pure family
at unit scale by the sum of an independent gamma variable of
shape \(\beta\) and a compound-Poisson variable with rate \(x\)
and exponential jumps of mean one. At \(\beta=0\), the gamma
summand is zero. We thus identify the logarithmic exponent with the gamma
component and the rational exponent with the Lah component
already considered in Subsection~\ref{sec:gamma-families}.
The first three polynomials are
\begin{equation}\label{lagapp:eq:initial-polynomials}
 \Lambda_0=1,\qquad \Lambda_1=x+\beta,\qquad
 \Lambda_2=x^2+2(\beta+1)x+\beta(\beta+1).
\end{equation}
For example, we obtain the explicit positive minor
\begin{equation}\label{lagapp:eq:first-minor}
 \Lambda_0\Lambda_2-\Lambda_1^2=\beta+2x.
\end{equation}
For the next leading minor, we obtain
\begin{equation}\label{lagapp:eq:second-minor}
 \det[\Lambda_{i+j}(x;\beta)]_{i,j=0}^2
 =2\beta^2(\beta+1)+12\beta(\beta+1)x
                   +12(2\beta+1)x^2+16x^3.
\end{equation}

We may also combine the Laguerre components with the tree models
above. In particular, for any real initial Stieltjes moment series
\(A(t)\), the exponential coefficients of
\begin{equation}\label{lagapp:eq:mixed-tree-family}
 A(t)(1-t)^{-\beta}
       \exp\!\left(\frac{xt}{1-t}+yU(t)+zH_1(t)\right)
\end{equation}
are coefficientwise Hankel-totally positive jointly in
\(\beta,x,y,z\). Applying Theorem~\ref{thm:lifting} to
\(U_1,V_1,U,H_1\), for which we have already verified the moment
hypotheses, we obtain the mixed example.

\subsection{Set partitions and classical moment transforms}

We count partitions whose blocks carry successive layers of
nonempty sets, with singleton layers allowed. At depth one we
obtain the ordinary Bell polynomials, whose coefficientwise
Hankel positivity was already studied in
\cite[Section~3]{Wang-Zhu}. Starting from the moment input
\(1,1,\ldots\), represented by \(\delta_1\), we obtain
admissible iterated exponents by composition. Applying lifting,
we extend the one-depth model to independent depths and colors
with arbitrary initial moments.

Define
\begin{equation}\label{eq:hierarchy-exponents}
 E_0(t)=t,\qquad E_{d+1}(t)=\exp(E_d(t))-1\quad(d\ge0).
\end{equation}

Iterating the construction of nonempty sets, we obtain the following families.

\begin{corollary}[Hierarchical Bell families]\label{thm:hierarchical-bell}
Fix positive integers \(d_1,\ldots,d_m\).
For any initial Stieltjes moment series \(A(t)\), the exponential coefficients of
\begin{equation}\label{eq:hierarchy-egf}
 A(t)\exp\!\left(\sum_{\ell=1}^{m}y_\ell E_{d_\ell}(t)\right)
\end{equation}
are coefficientwise Hankel-totally positive in \(\by\).
\end{corollary}

\begin{proof}
We represent the shifted cumulants of \(t\) and \(e^t-1\)
by \(\delta_0\) and \(\delta_1\), respectively. Applying
Theorem~\ref{thm:composition} successively to
\(E_{d+1}=(e^t-1)\circ E_d\), we obtain
\(E_d\in\mathcal C\). Applying Theorem~\ref{thm:lifting}
to the selected depths and \(A\), we obtain the conclusion.
\end{proof}

For \(A=1\), we use the exponential formula to interpret a
color-\(\ell\) component as \(d_\ell\) nested layers of nonempty
sets, allowing singleton layers.
At depth \(1\), we obtain ordinary set partitions and Bell polynomials,
while depth \(2\) groups nonempty blocks into nonempty collections.
For this standard iteration of the nonempty-set construction, see \cite{Comtet1974} for the exponential formula.
Differentiating \eqref{eq:hierarchy-exponents} at zero, we obtain
\begin{equation}\label{eq:display-046}
 E_d'(0)=1,\qquad E_d''(0)=d,\qquad
 E_d'''(0)=\frac{d(3d-1)}2.
\end{equation}

At the step from \(d\) to \(d+1\), the second and third
derivatives increase by \(1\) and \(3d+1\), respectively.

Writing \(Y=\sum_\ell y_\ell\) and
\(V=\sum_\ell d_\ell y_\ell\), the first pure-family moments are
\begin{equation}\label{eq:display-047}
 M_1=Y,\quad M_2=Y^2+V,\quad
 M_3=Y^3+3YV+
                 \sum_\ell\frac{d_\ell(3d_\ell-1)}2y_\ell.
\end{equation}
In particular, \(M_0M_2-M_1^2=V\).
For an arbitrary initial Stieltjes moment sequence, the component of
degree \(\binom r2\) in its leading determinant of order \(r\) is
\begin{equation}\label{eq:display-048}
 a_0^r\!\left(\prod_{j=0}^{r-1}j!\right)
             \left(\sum_\ell d_\ell y_\ell\right)^{\binom r2},
\end{equation}
by Theorem~\ref{thm:hamburger}.

For a formal power series \(L\) with \(L(0)=0\), put
\begin{equation}\label{eq:display-049}
 R_{n,k}(L)=\frac{n!}{k!}[t^n]L(t)^k.
\end{equation}
Formal substitution transforms an initial sequence into
\(b_n=\sum_{k=0}^{n}R_{n,k}(L)a_k\).

For the criterion based on row moment sequences, see
Wang and Zhu~\cite[Theorem~4.2]{Wang-Zhu}. By their Corollary~4.4(i),
we obtain the second-kind Stirling case. The first-kind case is the gamma
moment family \((s)_n\) at a nonnegative argument. For the Bessel and Lah cases, we may also use
\cite[Propositions~4.4(v) and~4.7(vi)]{Zhu2021}. The coefficient
arrays are the classical ones in those sources.

\begin{lemma}[Classical Stieltjes-preserving transforms]\label{cor:classical-transforms}
Each transform in the following table sends every Stieltjes moment
sequence to a Stieltjes moment sequence:
\begin{equation}\label{eq:display-050}
 \renewcommand{\arraystretch}{1.45}
 \begin{array}{c|c|c}
 \text{transform}&L(t)&R_{n,k}(L)\ \ (1\le k\le n)\\ \hline
 \text{Stirling, second kind}&e^t-1&\genfrac{\{}{\}}{0pt}{}{n}{k}\\
 \text{unsigned Stirling, first kind}&-\log(1-t)&
                                           \genfrac{[}{]}{0pt}{}{n}{k}\\
 \text{Lah}&t/(1-t)&\dfrac{n!}{k!}\binom{n-1}{k-1}\\
 \text{Bessel}&1-\sqrt{1-2t}&
                  \dfrac{(2n-k-1)!}{2^{n-k}(n-k)!(k-1)!}
 \end{array}
\end{equation}
In all four cases, \(R_{0,0}=1\) and \(R_{n,0}=0\) for \(n>0\).
\end{lemma}

For each $L(t)$ in four examples above, it is easy to obtain the corresponding shifted cumulant sequence
\((\kappa_{n+1})_{n\ge0}\) is a Stieltjes moment sequence.
After any finite composition of these known transforms, we may
convolve with the family defined by
\begin{equation}\label{eq:display-051}
 \exp\!\left(\sum_{\ell=1}^m y_\ell K_\ell(t)\right),
 \qquad K_\ell\in\mathcal C.
\end{equation}
By Theorem~\ref{thm:substitution}, we obtain Hankel minors with nonnegative
coefficients jointly in the indeterminates. We thus combine the four
known numerical preservation results with universal cumulant lifting.

For example, for any initial Stieltjes moment series \(A\) and fixed
\(q\ge0\), the coefficients of
\begin{equation}\label{eq:hybrid-application}
 A(e^t-1)(1-qT(T(t)))^{-y}
 \exp\!\left(x(1-\sqrt{1-2t})+
               \frac{zt}{1-t}+w(e^{e^t-1}-1)\right)
\end{equation}
are coefficientwise Hankel-totally positive jointly in \(x,z,w,y\).

Here the Stirling substitution and the Bessel, Lah,
depth-\(2\) Bell, and iterated-tree exponents all belong to
\(\mathcal C\). Applying Theorem~\ref{thm:substitution} to
these exponents and the Stirling substitution, we obtain
the asserted joint positivity.

\subsection{Compound-Poisson moment families and lognormal jumps}

Let \(\rho\) be any probability measure on \([0,\infty)\)
with all moments finite, and write \(b_n=\int x^n\rho(\dif x)\).
The jump exponent \(J_{\boldsymbol b}(t)=\sum_{n\ge1}b_nt^n/n!\)
in Theorem~\ref{thm:moment-input-families} has once-shifted
moments represented by \(x\rho(\dif x)\). Hence arbitrary
independently chosen jump laws give universal component families
with arbitrary initial Stieltjes moments. At numerical intensities,
these are the classical compound-Poisson moment families
\cite[Chapter~3]{LastPenrose2017}.

For any jump moment sequence \(b_n\), we can also count
set partitions with block weight \(b_{|B|}\) and a component
weight for each block. We specialize these weighted counts to
the known lognormal moments. Verifying the shifted moment
condition, we obtain coefficientwise Hankel positivity in the
Poisson-rate variables, including divergent moment series.

\begin{proposition}[Lognormal jumps]\label{prop:lognormal-jumps}
Fix \(\mu_\ell\in\R\) and \(\sigma_\ell^2\ge0\), and put
\begin{equation}\label{eq:lognormal-cumulants}
 K_\ell(t)=\sum_{n\ge1}
       \exp\!\left(n\mu_\ell+\frac{n^2\sigma_\ell^2}{2}\right)
                              \frac{t^n}{n!}.
\end{equation}
For any initial Stieltjes moment series \(A(t)\), the series
\(A(t)\exp(\sum_\ell y_\ell K_\ell(t))\) has an exponential
coefficient sequence that is jointly coefficientwise Hankel-totally
positive in \(\by\).
Whenever \(\sigma_\ell^2>0\), the series \(K_\ell\) has radius of
convergence zero.
\end{proposition}

\begin{proof}
Put \(\sigma_\ell=\sqrt{\sigma_\ell^2}\) and let
\(X_\ell=\exp(\mu_\ell+\sigma_\ell Z)\), where \(Z\)
is standard normal. Integrating against its density, we obtain
\begin{equation}\label{eq:display-052}
 \E X_\ell^n
 =\exp\!\left(n\mu_\ell+\frac{n^2\sigma_\ell^2}{2}\right).
\end{equation}
Thus \(b_n^{(\ell)}=\E X_\ell^n\) is Stieltjes and
\(K_\ell=J_{\boldsymbol b^{(\ell)}}\).
Applying Theorem~\ref{thm:moment-input-families} to these
jump moment inputs, with \(S=t\) and zero logarithmic
intensities, we obtain the positivity assertion.
For \(\sigma_\ell^2>0\), applying Stirling's formula to
the exponential coefficient, we have
\begin{equation}\label{eq:display-053}
 \left(\frac{\exp(n\mu_\ell+n^2\sigma_\ell^2/2)}{n!}\right)^{1/n}
 \sim\frac{e}{n}\exp(\mu_\ell+n\sigma_\ell^2/2)
 \longrightarrow\infty,
\end{equation}
so the radius of convergence is zero.
\end{proof}

For numerical \(\by\ge\boldsymbol0\), the pure-family coefficients
are the moments of
\begin{equation}\label{eq:display-054}
 S_{\by}=\sum_{\ell=1}^{m}\sum_{j=1}^{N_\ell}X_{\ell,j},
 \qquad N_\ell\sim\operatorname{Poisson}(y_\ell),
\end{equation}
where, for each \(\ell\), the variables \(X_{\ell,j}\) are
independent copies of the lognormal variable \(X_\ell\) constructed above.
We take all counts and jumps to be independent.
Conditioning on the counts, or using the set-partition formula,
we verify finiteness of every moment and recover the same formal
coefficients. Applying Theorem~\ref{thm:lifting} to the lognormal jump
exponents, we obtain nonnegativity of all mixed coefficients
in the Poisson rates in every Hankel minor. This remains valid when a lognormal jump of
positive variance has positive rate, although the moment generating
function is then infinite at every positive argument.
Indeed, the event of exactly one such jump and no other jumps has
positive probability, and \(\E e^{tX_\ell}=\infty\) for \(t>0\):
in the normal integral, we compare the term
\(t e^{\mu_\ell+\sigma_\ell z}\) with the quadratic normal
exponent and see that the former eventually dominates.
We therefore need the formulation in terms of formal power series
and finite moments for this application.

\subsection{Poisson differences and Gaussian perturbations}

The difference of independent Poisson counts has the classical
Skellam distribution \cite{GanKolaczyk2017}. Using its
twice-shifted moment measures, we obtain a general
Hamburger lifting result with arbitrary initial moments
and an independent Gaussian component.

\begin{corollary}[Poisson differences with an arbitrary initial measure]\label{cor:skellam}
Let \(A(t)=\sum_{n\ge0}a_nt^n/n!\) be any initial
Hamburger moment series, and fix \(\sigma^2\ge0\). Define
\begin{equation}\label{eq:skellam-egf}
 \sum_{n\ge0}S_n(u,v,w)\frac{t^n}{n!}
 =A(t)\exp\!\left(u(e^t-1)+v(e^{-t}-1)
                              +\frac{\sigma^2wt^2}{2}\right).
\end{equation}
For \(r\ge1\), the leading determinant
\(\det[S_{i+j}(u,v,w)]_{i,j=0}^{r-1}\) belongs to
\(\R_{\ge0}[u,v,w]\), has degree at most \(D=\binom r2\),
and has homogeneous part of degree \(D\) equal to
\begin{equation}\label{eq:skellam-top-degree}
 a_0^r\!\left(\prod_{j=0}^{r-1}j!\right)
                         (u+v+\sigma^2w)^D.
\end{equation}
Any additional common polynomial drift leaves these determinants
unchanged.
\end{corollary}

\begin{proof}
The twice-shifted cumulants of \(e^t-1\), \(e^{-t}-1\),
and \(\sigma^2t^2/2\) are represented by \(\delta_1\),
\(\delta_{-1}\), and \(\sigma^2\delta_0\), respectively.
Applying Theorem~\ref{thm:hamburger} to these exponents,
we obtain coefficientwise positivity, the degree bound, and
\eqref{eq:skellam-top-degree}, using second cumulants
\(1,1,\sigma^2\). Applying
Remark~\ref{rem:common-drifts} to the output family,
we obtain drift invariance.
\end{proof}

For \(A=1\) and numerical nonnegative parameters, we obtain the
moment sequence of \(N_u-N_v+\sigma\sqrt w\,Z\), where
\(\sigma=\sqrt{\sigma^2}\), \(N_u\sim\operatorname{Poisson}(u)\),
\(N_v\sim\operatorname{Poisson}(v)\), and \(Z\) is standard normal.
We take these three variables to be independent.
To see the positivity explicitly, put \(V=u+v+\sigma^2w\).
The first two nontrivial leading determinants are
\begin{align}
 \det[S_{i+j}]_{i,j=0}^{1}&=V,\label{eq:skellam-h2}\\
 \det[S_{i+j}]_{i,j=0}^{2}
 &=2V^3+4uv+\sigma^2w(u+v).\label{eq:skellam-h3}
\end{align}
Writing \(\kappa_j\) for the \(j\)-th cumulant of the pure
family in \eqref{eq:skellam-egf}, these determinants are
\(\kappa_2\) and \(2\kappa_2^3+\kappa_2\kappa_4-\kappa_3^2\).
Here \(\kappa_2=V\), \(\kappa_3=u-v\), and \(\kappa_4=u+v\).
After cancelling the signed odd-cumulant contribution, we obtain
exactly \eqref{eq:skellam-h3}.

The coefficientwise assertion concerns leading determinants,
as in Theorem~\ref{thm:hamburger}.

\subsection{Further consequences for lifted Hankel matrices}
\label{sec:lifted-consequences}

We conclude the applications with three standard consequences of
the coefficientwise Hankel total positivity proved above. For the first,
we use the classical Sylvester identity
\cite[Theorem~1.1]{Konvalinka2007}. For the second, we use the
standard bidiagonal and Pascal factorizations; see
\cite{Pinkus2010} and Layman's binomial-invariance formula
\cite[Theorem~1]{Layman2001}. For the third, we use a bidiagonal
factorization of the moment matrices appearing in classical
Stieltjes--Pad\'e approximation \cite{Simon1998}.
We include the deductions to explain how we use these classical
tools for the output sequences of the lifting theorem.

\begin{corollary}[Standard consequences for the lifted sequences]
\label{cor:lifted-matrix-consequences}
Let \((F_n(\by))_{n\ge0}\) be an output of
Theorem~\ref{thm:lifting} under its equivalent conditions.
Then the following assertions hold.
\begin{enumerate}
\renewcommand{\labelenumi}{\textup{(\roman{enumi})}}
\item The centered kernel
\begin{equation}\label{eq:lifted-centered-kernel}
 K_{i,j}(\by):=F_0(\by)F_{i+j}(\by)-F_i(\by)F_j(\by)
 \qquad(i,j\ge1)
\end{equation}
is coefficientwise totally positive in \(\by\).
For increasing positive index lists \(I,J\) of length \(r\),
we have
\begin{equation}\label{eq:lifted-centered-minors}
 \det K[I,J]=F_0^{r-1}
       \det[F_{u+v}]_{u\in\{0\}\cup I,\,v\in\{0\}\cup J}.
\end{equation}
\item Fix \(h\ge0\) and introduce independent indeterminates
\(c_1,\ldots,c_h,u,v\). Define
\begin{equation}\label{eq:lifted-weight-polynomial}
 P(z):=\prod_{\ell=1}^h(z+c_\ell)=\sum_{\ell=0}^h p_\ell z^\ell,
 \qquad \widetilde F_n:=\sum_{\ell=0}^h p_\ell F_{n+\ell},
\end{equation}
and
\begin{equation}\label{eq:lifted-affine-sequence}
 G_n:=\sum_{k=0}^n\binom nk u^{n-k}v^k\widetilde F_k.
\end{equation}
Both \((\widetilde F_n)\) and \((G_n)\) are coefficientwise
Hankel-totally positive jointly in the original and new parameters.
Moreover,
\begin{equation}\label{eq:lifted-affine-leading}
 \det[G_{i+j}]_{i,j=0}^{r-1}
 =v^{r(r-1)}\det[\widetilde F_{i+j}]_{i,j=0}^{r-1}.
\end{equation}
\item For an independent indeterminate \(s\), \(h\in\{0,1\}\),
and \(r\ge1\), the polynomial
\begin{equation}\label{eq:lifted-pade-denominators}
 \det[F_{i+j+h}(\by)+sF_{i+j+h+1}(\by)]_{i,j=0}^{r-1}
\end{equation}
has nonnegative coefficients jointly in \(s,\by\).
\end{enumerate}
\end{corollary}

\begin{proof}
For \textup{(i)}, applying Sylvester's identity to the Hankel
matrix \(H=[F_{i+j}]_{i,j\ge0}\), with the fixed row and column
index both equal to zero, we obtain
\eqref{eq:lifted-centered-minors}. Since this is a polynomial
identity, we may also use it when \(F_0=0\). By
Theorem~\ref{thm:lifting}, every factor on the right has
nonnegative coefficients, and we obtain the required positivity.

For \textup{(ii)}, define lower bidiagonal matrices
\(T(c)_{i,j}=c\mathbf1_{\{i=j\}}+\mathbf1_{\{i=j+1\}}\).
Writing \(H(F)=[F_{i+j}]_{i,j\ge0}\), we have
\begin{equation}\label{eq:lifted-weight-factorization}
 H(\widetilde F)=H(F)\prod_{\ell=1}^hT(c_\ell).
\end{equation}
For every weight factor, we have coefficientwise total positivity.
For the
affine change, put \(B_{n,k}=\binom nk u^{n-k}v^k\) when
\(k\le n\), and \(B_{n,k}=0\) otherwise. By the binomial
identity, we obtain
\begin{equation}\label{eq:lifted-affine-factorization}
 H(G)=BH(\widetilde F)B^{\mathsf T}.
\end{equation}
For positive numerical \(u,v\), every minor of \(B\) is a
minor of the Pascal matrix multiplied by the row and column factors
\(u^{n}\) and \((v/u)^k\). For a nonzero minor, we can match row and column indices within
the lower triangular support. We therefore obtain a factor
\(u^{\sum n-\sum k}\) with nonnegative exponent. Thus the
scaling gives a polynomial identity with nonnegative coefficients,
valid also at \(u=0\) or \(v=0\), and we conclude that \(B\)
is coefficientwise totally positive.
Applying the finite Cauchy--Binet formula to
\eqref{eq:lifted-weight-factorization} and
\eqref{eq:lifted-affine-factorization}, we obtain the asserted
Hankel total positivity. For the leading \(r\times r\) block,
we use the triangular matrix \(B\), with diagonal
\(1,v,\ldots,v^{r-1}\). Taking determinants, we obtain
\eqref{eq:lifted-affine-leading}.

For \textup{(iii)}, let \(H^{(h)}=[F_{i+j+h}]_{i,j\ge0}\),
and let \(L\) be the down-shift matrix,
\(L_{i,j}=\mathbf1_{\{i=j+1\}}\). The required lower bidiagonal
matrix is \(B(s):=I+sL\), with diagonal entries one and
subdiagonal entries \(s\). Indeed,
\begin{equation}\label{eq:lifted-pade-factorization}
 [H^{(h)}B(s)]_{i,j}
 =F_{i+j+h}+sF_{i+j+h+1}.
\end{equation}
Thus \eqref{eq:lifted-pade-denominators} is the determinant of
the leading \(r\times r\) block of \(H^{(h)}(I+sL)\).
Every minor of \(H^{(h)}\) is a minor of \(H(F)\), so it has
nonnegative coefficients by Theorem~\ref{thm:lifting}.

The matrix \(I+sL\) is coefficientwise totally positive.
For the first \(r\) columns, its nonzero entries are confined
to rows \(0,\ldots,r\). Hence the Cauchy--Binet expansion for
the leading \(r\times r\) block is a finite sum over
\(r\)-element subsets of \(\{0,\ldots,r\}\), and every term
has nonnegative coefficients. This proves \textup{(iii)},
including the specialization \(s=0\).
\end{proof}

When \(F_0=1\) and the parameters have numerical nonnegative
values, we may interpret \eqref{eq:lifted-centered-kernel} as a
covariance kernel of powers of a nonnegative random variable.
In \eqref{eq:lifted-pade-denominators}, we obtain coefficientwise
positivity for the moment-matrix denominators used in classical
Stieltjes--Pad\'e bounds whenever those approximants are defined.
To establish convergence of these bounds, one needs additional
moment-determinacy hypotheses. Here we assert the displayed
coefficientwise positivity.

\section{Some final remarks }\label{sec:scope}

For colored functional digraphs, we keep the cycle weights
\(q_1,\ldots,q_m\) fixed and nonnegative. In
Theorem~\ref{thm:digraph}, we obtain joint coefficientwise
positivity in the component intensities \(y_1,\ldots,y_m\).
It remains to establish joint coefficientwise positivity in both
\(\bq\) and \(\by\).

We likewise keep the real equation coefficients fixed in
Theorems~\ref{thm:pf-implicit-lifting} and
\ref{thm:pf-autonomous-lifting}. We retain all independent component indeterminates in their
conclusions. To retain the equation parameters as well, we would
need to extend the lifting argument beyond the real moment-sequence
hypotheses used here.

For the Schl\"afli-Gessel-Seo family, we prove the shifted
moment condition for every fixed \(a,b\ge0\) in
Lemma~\ref{lem:sgs-shifted-moments}. Applying lifting to those
inputs, we obtain coefficientwise Hankel total positivity in
\(x\), arbitrary initial moments, and independent colors in
Theorem~\ref{thm:sgs-fixed-lifting}.
We leave the joint \(x,a,b\) assertion of
Conjecture~\ref{conj:sgs-joint} unresolved. Nonnegative minors for
every numerical specialization of \(a,b\) do not imply nonnegative
coefficients in those parameters.

Throughout the lifting arguments, we use formal power series
and finite moments. In the lognormal example, we obtain positivity
even though the moment generating series has zero radius of
convergence. To deduce uniqueness of a representing measure,
one would need additional hypotheses beyond the positivity
established here.

For the bivariate forest model in
Corollary~\ref{treeapp:cor:forests}, we keep \(w\ge-1\) fixed.
The corresponding joint assertions in \((x,w)\) and in the shifted
variables \((x,w+1)\), discussed by Sokal in
\cite[Conjectures~6.11(b) and~6.12(b)]{Sokal2023}, are not proved
by fixing that parameter. We use the same convention for
structural parameters in the further tree component families.
In the Laguerre application, we can express the shape parameter
as the intensity of a logarithmic exponent and hence retain it
as an indeterminate in Theorem~\ref{lagapp:thm:universal}.

\color{black}
\hypersetup{allcolors=black}

\endgroup

\begin{thebibliography}{99}
\small
\setlength{\itemsep}{4pt}
\setlength{\parskip}{0pt}
\bibitem{Andreief1886}
C.~Andr\'eief,
Note sur une relation entre les int\'egrales d\'efinies des produits
des fonctions,
\emph{M\'em. Soc. Sci. Phys. Nat. Bordeaux} (3) \textbf{2} (1886), 1--14.

\bibitem{BGKP2016}
A.~Belton, D.~Guillot, A.~Khare and M.~Putinar,
\href{https://doi.org/10.1016/j.aim.2016.04.016}
{Matrix positivity preservers in fixed dimension. I},
\emph{Adv. Math.} \textbf{298} (2016), 325--368.

\bibitem{BGKP2022}
A.~Belton, D.~Guillot, A.~Khare and M.~Putinar,
\href{https://doi.org/10.4171/JEMS/1145}
{Moment-sequence transforms},
\emph{J. Eur. Math. Soc.} \textbf{24} (2022), no.~9, 3109--3160.

\bibitem{Bre89}
F. ~Brenti, Unimodal, log-concave, and P\'olya frequency sequences in
combinatorics, \emph{Mem. Amer. Math. Soc.} \textbf{413} (1989).

\bibitem{ChenSokal2024}
X.~Chen and A.~D.~Sokal,
Total positivity of some polynomial matrices that enumerate labeled
trees and forests. II. Rooted labeled trees and partial functional
digraphs,
\emph{Adv. Appl. Math.} \textbf{157} (2024), 102703.

\bibitem{Comtet1974}
L.~Comtet,
\emph{Advanced Combinatorics},
D.~Reidel, Dordrecht, 1974.

\bibitem{DDPS23}
B. Deb, A. Dyachenko, M. Pétréolle, A. D. Sokal, Lattice paths and branched continued fractions. III. Generalizations of the Laguerre, rook and Lah polynomials, arXiv:2312.11081 (2023).
\bibitem{DingZeng2024}
M.-J.~Ding and J.~Zeng,
\href{https://doi.org/10.1016/j.disc.2024.113953}
{Proof of an explicit formula for a series from Ramanujan's Notebooks
via tree functions},
\emph{Discrete Math.} \textbf{347} (2024), no.~6, 113953.



\bibitem{EPS2020}
A.~Elvey Price and A.~D.~Sokal,
\href{https://doi.org/10.37236/9571}
{Phylogenetic trees, augmented perfect matchings, and a Thron-type
continued fraction (T-fraction) for the Ward polynomials},
\emph{Electron. J. Combin.} \textbf{27} (2020), no.~4, Paper~4.6, 36~pp.


\bibitem{Flajolet1980}
P.~Flajolet,
\href{https://algo.inria.fr/flajolet/Publications/Flajolet80b.pdf}
{Combinatorial aspects of continued fractions},
\emph{Discrete Math.} \textbf{32} (1980), 125--161.

\bibitem{FominZelevinsky1999}
S.~Fomin and A.~Zelevinsky,
\href{https://doi.org/10.1090/S0894-0347-99-00295-7}
{Double Bruhat cells and total positivity},
\emph{J. Amer. Math. Soc.} \textbf{12} (1999), no.~2, 335--380.

\bibitem{GanKolaczyk2017}
H.~L.~Gan and E.~D.~Kolaczyk,
\href{https://arxiv.org/abs/1708.04018}
{Approximation of the difference of two Poisson-like counts with Skellam},
arXiv:1708.04018 (2017).

\bibitem{GinerSmyth2016}
G.~Giner and G.~K.~Smyth,
\href{https://journal.r-project.org/articles/RJ-2016-024/}
{statmod: Probability calculations for the inverse Gaussian distribution},
\emph{The R Journal} \textbf{8} (2016), no.~1, 339--351.

\bibitem{Goldie1967}
C.~Goldie,
\href{https://doi.org/10.1017/S0305004100042225}
{A class of infinitely divisible random variables},
\emph{Proc. Cambridge Philos. Soc.} \textbf{63} (1967), 1141--1143.

\bibitem{HornJohnson2013}
R.~A.~Horn and C.~R.~Johnson,
\href{https://doi.org/10.1017/CBO9781139020411}{\emph{Matrix Analysis}},
2nd ed., Cambridge University Press, Cambridge, 2013.

\bibitem{JohnsonKempKotz2005}
N.~L.~Johnson, A.~W.~Kemp and S.~Kotz,
\href{https://doi.org/10.1002/0471715816}{\emph{Univariate Discrete Distributions}},
3rd ed., Wiley, Hoboken, NJ, 2005.


\bibitem{JV2015}
M.~Josuat-Verg\`es,
\href{https://doi.org/10.1007/s11139-015-9684-x}
{Derivatives of the tree function},
\emph{Ramanujan J.} \textbf{38} (2015), no.~1, 1--15;
\href{https://arxiv.org/abs/1310.7531}{arXiv:1310.7531}.


\bibitem{Kar68}
S.~Karlin,
\emph{Total Positivity}, Vol.~I,
Stanford University Press, Stanford, CA, 1968.

\bibitem{KodamaWilliams2014}
Y.~Kodama and L.~K.~Williams,
\href{https://doi.org/10.1007/s00222-014-0506-3}
{KP solitons and total positivity for the Grassmannian},
\emph{Invent. Math.} \textbf{198} (2014), no.~3, 637--699.

\bibitem{Konvalinka2007}
M.~Konvalinka,
\href{https://doi.org/10.37236/960}
{Non-commutative Sylvester's determinantal identity},
\emph{Electron. J. Combin.} \textbf{14} (2007), no.~1, Research Paper~42.

\bibitem{LastPenrose2011}
G.~Last and M.~D.~Penrose,
\href{https://arxiv.org/abs/0909.3205}{Poisson process Fock space
representation, chaos expansion and covariance inequalities},
\emph{Probab. Theory Related Fields} \textbf{150} (2011), 663--690.

\bibitem{LastPenrose2017}
G.~Last and M.~Penrose,
\emph{Lectures on the Poisson Process},
Cambridge University Press, Cambridge, 2017.

\bibitem{Layman2001}
J.~W.~Layman,
\href{https://cs.uwaterloo.ca/journals/JIS/VOL4/LAYMAN/hankel.pdf}
{The Hankel transform and some of its properties},
\emph{J. Integer Seq.} \textbf{4} (2001), Article~01.1.5.

\bibitem{Lusztig1994}
G.~Lusztig,
Total positivity in reductive groups,
in \emph{Lie Theory and Geometry},
J.-L.~Brylinski, R.~Brylinski, V.~Guillemin and V.~Kac (eds.),
Progress in Mathematics \textbf{123},
Birkh\"auser, Boston, MA, 1994, pp.~531--568.

\bibitem{Macdonald1995}
I.~G.~Macdonald,
\emph{Symmetric Functions and Hall Polynomials},
2nd ed., Oxford University Press, Oxford, 1995.

\bibitem{Mlotkowski2010}
W.~M\l otkowski,
\href{https://doi.org/10.4171/DM/318}
{Fuss-Catalan numbers in noncommutative probability},
\emph{Doc. Math.} \textbf{15} (2010), 939--955.

\bibitem{DLMF-Gamma}
NIST Digital Library of Mathematical Functions,
\emph{Gamma function: integrals and rising factorials},
Eqs.~5.2.1, 5.2.5 and~5.9.1,
\href{https://dlmf.nist.gov/5.2}{dlmf.nist.gov/5.2} and
\href{https://dlmf.nist.gov/5.9}{5.9}.


\bibitem{DLMF-Laguerre}
NIST Digital Library of Mathematical Functions,
\emph{Laguerre polynomials: orthogonality, finite formulas and generating functions},
Table~18.3.1, Eqs.~18.5.12 and~18.12.13,
\href{https://dlmf.nist.gov/18.3\#T1}{dlmf.nist.gov/18.3\#T1},
\href{https://dlmf.nist.gov/18.5\#E12}{18.5\#E12} and
\href{https://dlmf.nist.gov/18.12\#E13}{18.12\#E13}.


\bibitem{DLMF-Norms}
NIST Digital Library of Mathematical Functions,
\emph{Orthogonality, norms and moment determinants},
Section~18.2(ix), Tables~18.3.1 and~18.19.1,
\href{https://dlmf.nist.gov/18.2\#ix}{dlmf.nist.gov/18.2(ix)},
\href{https://dlmf.nist.gov/18.3\#T1}{18.3.1} and
\href{https://dlmf.nist.gov/18.19\#T1}{18.19.1}.

\bibitem{PS19}
M. P\'{e}tr\'{e}olle, A.D. Sokal, Lattice paths and branched
continued fractions. II. Multivariate Lah polynomials and Lah
symmetric functions, European J. Combin. 92, Article 103235.

\bibitem{PSZ2023}
M.~P\'etr\'eolle, A.~D.~Sokal and B.-X.~Zhu,
Lattice paths and branched continued fractions: An infinite sequence
of generalizations of the Stieltjes--Rogers and Thron--Rogers
polynomials, with coefficientwise Hankel-total positivity,
\emph{Mem. Amer. Math. Soc.} \textbf{291} (2023), no.~1450, 154~pp.

\bibitem{Pinkus2010}
A.~Pinkus,
\emph{Totally Positive Matrices},
Cambridge University Press, Cambridge, 2010.

\bibitem{PS64}
G. P\'olya, G. Szeg\"o,
Problems and Theorems in Analysis, Vol. II,
Springer-Verlag, Berlin, 1976.

\bibitem{Rietsch2003}
K.~Rietsch,
\href{https://doi.org/10.1090/S0894-0347-02-00412-5}
{Totally positive Toeplitz matrices and quantum cohomology of
partial flag varieties},
\emph{J. Amer. Math. Soc.} \textbf{16} (2003), no.~2, 363--392.

\bibitem{SSV2012}
R.~L.~Schilling, R.~Song and Z.~Vondra\v cek,
\href{https://doi.org/10.1515/9783110269338}
{\emph{Bernstein Functions: Theory and Applications}},
2nd ed., de Gruyter Studies in Mathematics \textbf{37},
Walter de Gruyter, Berlin, 2012.

\bibitem{Schmudgen2017}
K.~Schm\"udgen,
\href{https://doi.org/10.1007/978-3-319-64546-9}{\emph{The Moment Problem}},
Graduate Texts in Mathematics \textbf{277}, Springer, Cham, 2017.

\bibitem{Simon1998}
B.~Simon,
\href{https://arxiv.org/abs/math-ph/9906008}
{The classical moment problem as a self-adjoint finite difference operator},
\emph{Adv. Math.} \textbf{137} (1998), 82--203.

\bibitem{Sokal_totalpos}  
A.D. Sokal, Coefficientwise total positivity
   (via continued fractions) for some Hankel matrices of combinatorial
   polynomials, in preparation.
   
\bibitem{Sokal2023}
A.~D.~Sokal,
Total positivity of some polynomial matrices that enumerate labeled
trees and forests. I. Forests of rooted labeled trees,
\emph{Monatsh. Math.} \textbf{200} (2023), 389--452.


\bibitem{Speed1983}
T.~P.~Speed,
\href{https://doi.org/10.1111/j.1467-842X.1983.tb00391.x}
{Cumulants and partition lattices},
\emph{Austral. J. Statist.} \textbf{25} (1983), no.~2, 378--388.


\bibitem{Stanley}
R.~P.~Stanley, \emph{Enumerative Combinatorics, Volume 2},
Cambridge University Press, 1999, Sections~5.1, 5.4, 7.10 and~7.15.

\bibitem{Wang-Zhu}
Y.~Wang and B.-X.~Zhu,
\emph{Log-convex and Stieltjes moment sequences}, Adv. Appl. Math.
\textbf{81} (2016), 115--127.


\bibitem{Zeng1999}
J.~Zeng,
\href{https://doi.org/10.1023/A:1009809224933}
{A Ramanujan sequence that refines the Cayley formula for trees},
\emph{Ramanujan J.} \textbf{3} (1999), no.~1, 45--54.

\bibitem{Zhu20}
B.-X. Zhu, A generalized Eulerian triangle from staircase tableaux
and tree-like tableaux, J. Combin. Theory Ser. A 172 (2020), Article
105206, 32 pp.

\bibitem{Zhu21}
B.-X. Zhu, Stieltjes moment properties and continued fractions from
combinatorial triangles, Adv. in Appl. Math. 130 (2021), Article
102232, 33 pp.
\bibitem{Zhu212}
B.-X. Zhu, On a Stirling-Whitney-Riordan triangle, J. Algebraic
Combin. 54 (2021) 999--1019.

\bibitem{Zhu2021}
B.-X.~Zhu,
\href{https://doi.org/10.1137/20M1379952}
{Total positivity from the exponential Riordan arrays},
\emph{SIAM J. Discrete Math.} \textbf{35} (2021), 2971--3003.

\bibitem{Zhu2024}
B.-X.~Zhu,
\href{https://arxiv.org/abs/2202.03793}{Coefficientwise Hankel-total
positivity of the row-generating polynomials for the output matrices
of certain production matrices},
\emph{Adv. Math.} \textbf{447} (2024), 109695.

\bibitem{Zhu2026Riordan}
B.-X.~Zhu,
\href{https://doi.org/10.1112/jlms.70450}
{Coefficientwise total positivity and Stieltjes moment properties from
Riordan arrays},
\emph{J. Lond. Math. Soc.} (2) \textbf{113} (2026), no.~2, e70450.

\end{thebibliography}
\end{document}